\documentclass[12pt]{article}
\usepackage{amsfonts}
\usepackage{amsmath}
\usepackage{amssymb}
\usepackage{graphicx}

\usepackage{amsfonts}
\usepackage{amsmath}
\newtheorem{theorem}{Theorem}[section]
\newtheorem{lemma}[theorem]{Lemma}
\newtheorem{corollary}[theorem]{Corollary}

\newcommand{\HS}[3]{\left(\frac{#1,#2}{#3}\right)}

\newcommand{\IR}{\mathbb R}
\newcommand{\IQ}{\mathbb Q}
\newcommand{\IZ}{\mathbb Z}
\newcommand{\IC}{\mathbb C}
\newcommand{\IH}{\mathbb H}

\newcommand{\G}{\Gamma}

\def\IH{{\Bbb H}} 
 
\def\IR{{\Bbb R}} 
 
\def\IS{{\Bbb S}}

\def\IZ{{\Bbb Z}} 
 
\def\IC{\Bbb C} 
\def\ID{{\Bbb D}}
\def\tr{{\rm tr}} 
\def\oC{\hat{\IC}} 
\def\tr{\mbox{\rm{tr\,}}}

\def\PSL{\mbox{\rm{PSL}}}

\def\em{\it}

\title{The $(4,p)$-arithmetic hyperbolic lattices, $p\geq 2$,  in three dimensions. }
\author{G.J. Martin, K. Salehi and Y. Yamashita \thanks{Research of all authors supported in part by the NZ Marsden Fund.}}
\date{}
\begin{document}

\maketitle

\begin{abstract}
We identify the finitely many arithmetic lattices $\Gamma$ in the orientation preserving isometry group of hyperbolic $3$-space $\IH^3$ generated by an element of order $4$ and and element of order $p\geq 2$.  Thus $\Gamma$ has a presentation of the form
\[ \Gamma\cong\langle f,g: f^4=g^p=w(f,g)=\cdots=1 \rangle \]  
We find that necessarily $p\in \{2,3,4,5,6,\infty\}$,  where $p=\infty$ denotes that $g$ is a parabolic element,  the total degree of the invariant trace field $k\Gamma=\IQ(\{\tr^2(h):h\in\Gamma\})$ is at most $4$, and each orbifold is either a two bridge link of slope $r/s$ surgered with $(4,0)$, $(p,0)$ Dehn surgery (possibly a two bridge knot if $p=4$) or a Heckoid group with slope $r/s$ and $w(f,g)=(w_{r/s})^r$ with $r\in \{1,2,3,4\}$.  We give a discrete and faithful representation in $PSL(2,\IC)$ for each group and identify the associated number theoretic data.
\end{abstract}

\section{Introduction}

In this paper we continue our long running programme to identify (up to conjugacy) all the finitely many arithmetic lattices $\Gamma$ in the group of orientation preserving isometries of hyperbolic $3$-space $Isom^+(\IH^3)\cong PSL(2,\IC)$ generated by two elements of finite order $p$ and $q$,  we also allow $p=\infty$ or $q=\infty$ to mean that a generator is parabolic.  In \cite{MM} it is proved that there are only finitely many such lattices.  In fact it is widely expected that there are only finitely many two generator arithmetic lattices in total in $PSL(2,\IC)$,  though this remains unknown.  There are infinitely many distinct three generator arithmetic lattices. 

It is well known that for such a lattice as $\Gamma$ above, $\IH^3/\Gamma$ is geometrically the $3$-sphere with a marked trivalent graph of singular points -- the vertices of the graph corresponding to the finite spherical triangle groups and dihedral groups. The possible graphs for the cases at hand can be seen in Figure 1. These precisely follow the descriptions of all lattices in $PSL(2,\IC)$ generated by two parabolic elements as proved by Aimi, Lee,Sakai, and Sakuma, and also Akiyoshi,  Ohshika, Parker,  Sakuma and  Yoshida, \cite{ALSS,AOPSY}, who resolved a conjecture of Agol to this effect.  Our results suggest that this conjecture (modified in the obvious way as per Figure 1.) is valid for all groups generated by two elements of finite order which is not freely generated by those two elements.

\medskip

In two dimensions a lattice in $Isom^+(\IH^2)$ generated by two elements of finite order is necessarily a triangle group. Takeuchi \cite[Theorem 3]{Take} has identified the $82$ arithmetic Fuchsian triangle groups.  Here is a summary of what is known in dimension three. 
 
 \begin{itemize}
\item There are exactly $4$ arithmetic lattices generated by two parabolic elements. These are all knot and link complements and so torsion free. The explicit groups can be found in \cite{GMM}.
\item There are exactly $14$ arithmetic lattices generated by a parabolic element and an elliptic of order $2\leq p < \infty$.  Necessarily $p\in \{2,3,4,6\}$. Six of these groups have $p=2$, three have $p=3$, and three have $p=4$. The explicit groups can be found in \cite{CMMO}.
\item There are exactly $16$ arithmetic lattices generated by two elements of finite orders $p$ and $q$ with $p,q \geq  6$. In each case $p=q$,  twelve of these groups have $p=q=6$,  and two have $p=q=12$ and one each for $p=q=8$ and $p=q=10$. The explicit groups can be found in \cite{MM} 
\end{itemize}
Most of these groups are orbifold surgeries on two bridge knots and links of ``small'' slope.  Here. we will meet examples of Dehn surgeries on knots of more than 12 crossings.
 
 \medskip
 
 From this information the reader can see that there are basically four cases left to deal with in order to compete the identification we seek. Those are when one generator is elliptic of order $q=2,3,4,5$ and $2\leq p < \infty$.  In fact a discrete sugroup of $Isom^+(\IH^3)$ generated by two elements of the same finite order (say $p$) admits a $\IZ_2$ extension to a group generated by an element of order two and an element of order $p$.  As a finite extension remains an arithmetic lattice,  there really are now only three cases to deal with,   $q\in\{3,4,5\}$ and here we assume that $q=4$ and $2\leq p\leq \infty$.  Below is a summary of our main result.
 
 \begin{theorem}\label{mainthm} Let $\Gamma=\langle f,g \rangle$ be an arithmetic lattice generated by $f$ of order $4$ and $g$ of finite order $p\geq 2$.  Then $p\in \{2,3,4,5,6\}$ and the degree of the invariant trace field 
 \[ k\Gamma=\IQ(\{\tr^2(h):h\in\Gamma\}) \]
  is at most $4$.
 \begin{enumerate}
 \item If $p=2$,  then there are fifty four groups.
 \item If $p=3$,  then there are ten groups.
 \item If $p=4$,  then there are twenty seven groups.
  \item If $p=5$,  then there is one group.
\item If $p=6$,  then there are five groups.    
 \end{enumerate}
 \end{theorem}

 \medskip
 
 The groups appearing in Theorem \ref{mainthm} all fall into the pattern of Heckoid groups as described in \cite{ALSS,AOPSY}.  The singular graph of $\IH^3/\Gamma$ is one of the following.
 
 \scalebox{0.5}{\includegraphics[viewport=-40 520 570 800]{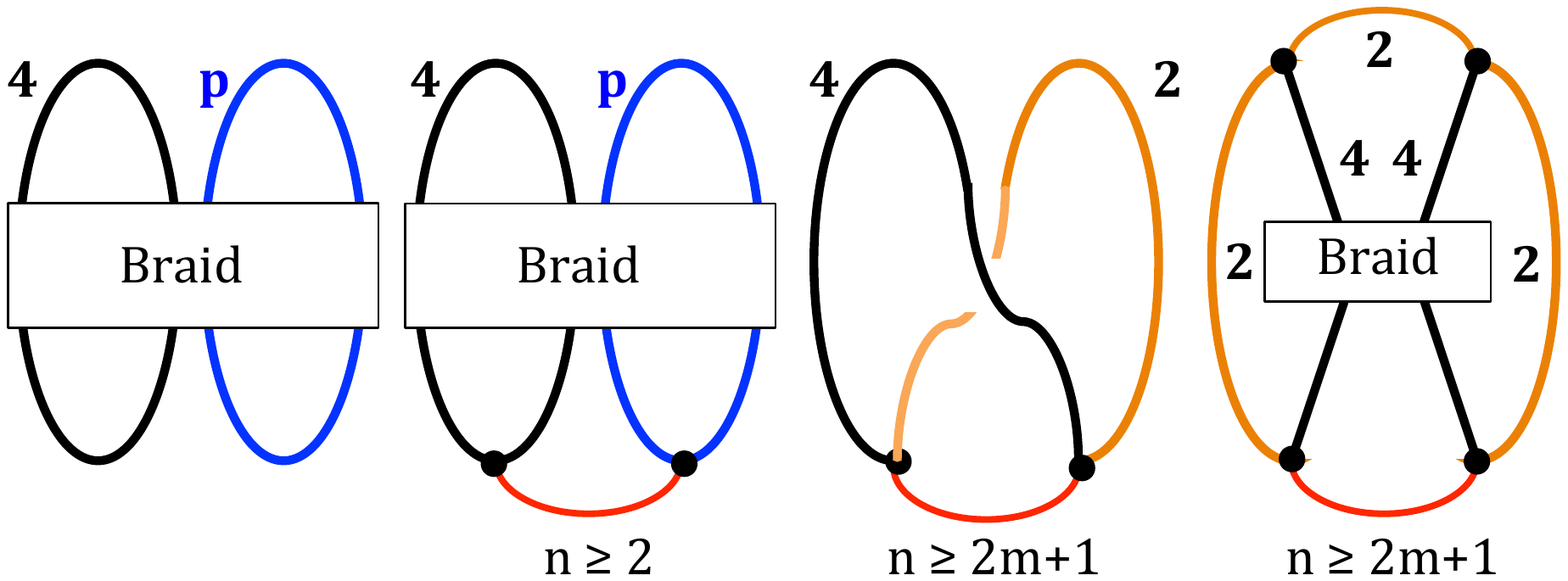}}\\
 {\bf Figure 1.} {\em The Heckoid groups.  When $p\neq 4$ the braid must have an even number of crossings in the first example.  Here $m\geq 1$.}
   
  \bigskip
  
 In the tables that follow we set
 \begin{equation}
 \gamma =\gamma(f,g) = \tr [f,g]-2.
 \end{equation}
 and give the minimal polynomial for $\gamma$.  With the orders of the generators,  this completely determines the arithmetic structure of the group,  the invariant trace field $k\Gamma=\IQ(\gamma)$ and the associated quaternion algebra -- we discuss these things below.
 
 Then $\gamma$ is an invariant of the Nielsen class of the generating pair $\{f,g\}$. With the order of the generators,  which we shall always assume are primitive (that is their trace should be $\pm 2\cos \big(\frac{\pi}{p}\big)$) and $\gamma$ it is straightforward to construct a discrete and faithful representation of the group in $PSL(2,\IC)$, \cite{GM1}. We give this. representation below at (\ref{XY}).  Note that the very difficult general problem of proving discreteness is overcome by the arthimetic criterion as described in \cite{GMMR} and \cite{MMpq}.  We recall those results,  and in particular the identification theorem,  below at Theorem \ref{idthm}.
 
For each group $\Gamma$, as noted,  we give the minimal polynomial for $\gamma$ and the discriminant of the invariant trace field $k\Gamma=\IQ(\gamma)$.    Next the column $(Farey,order)$ indicates a particular word in the group $\langle f,g \rangle$ has a particular order (or is the identity).  Briefly this word is found as follows.  The simple closed curves on the $4$ times punctured sphere  $\IS_4$  separating one pair of points from another are enumerated by their slope,  coming from the Farey expansion.  The enumeration of these simple closed curves informs the deformation theory of Keen-Series-Maskit \cite{KS,KS2,KSM} which relates this slope (and the geodesic in the homotopy class) to a bending deformation along this geodesic of $\IS_4$ which terminates on the boundary of moduli space as the length of this curve shrinks to zero (and so the associated word becomes parabolic). This bending locus is called a pleating ray.   Bending further creates cone manifolds some of which are discrete lattices when the cone angle becomes $2\pi/n$ for a positive integer $n$.  The recent results of \cite{ALSS,AOPSY} show this process to describe all the discrete and faithful representations of groups generated by two parabolic elements which are not free.  There is of course a strong connection with the Schubert normal form of a two-bridge knot or link, \cite{BZ}.  If one takes a two-bridge link,  say with Schubert normal form $r/s$, and performs orbifold $(p,0)$ and $(q,0)$ Dehn surgery on the components of the link,  and if the result is hyperbolic,  then the orbifold fundamental group of the space so formed is a discrete subgroup of $PSL(2,\IC)$ generated by elliptic elements of order $p$ and $q$,  say $f$ and $g$ respectively - the images of the two meridians around each link component (or the images of the two parabolic generators of the group). The word associated with the slope remains the identity.  However the value $\gamma(f,g)$ may not always lie on the pleating ray (we give more details below) but possibly on the pleating ray $1-\frac{r}{s}$ - this is simply because of the way we chose to enumerate slopes between.  This then also provides the Conway notation for the knot of link complement via the continued fraction expansion of the slope $r/s$.  There is one further issue here and that is that the exact choice of generators in $PSL(2,\IC)$ matters since we classify all Nielson pairs and in some cases there may be different pairs generating the same group and obviously they will have different relations. In the description of the group (in terms of the pleating ray $\gamma$ lies on) we will have the following setup:

 If $X,Y\in PSL(2,\IC)$,  representing $f$ and $g$,  are given by 
\begin{equation}\label{XY} X= \left[\begin{array}{cc} \zeta & 1\\ 0 &\bar\zeta \end{array}\right], Y= \left[\begin{array}{cc} \eta & 0\\ \mu &\bar\eta \end{array}\right],\; \zeta = \cos\frac{\pi}{p}+i\sin\frac{\pi}{p},\;\; \eta = \cos\frac{\pi}{q}+i\sin\frac{\pi}{q},\end{equation}
we compute
\begin{equation}\label{mug}
\gamma(f,g) = \mu (\mu -4 \sin \frac{\pi }{p}\; \sin  \frac{\pi }{q}).
\end{equation}
Thus two different choices for $\mu$ lead to the same group up to M\"obius conjugacy.  One of these values of $\mu$ lies on the pleating ray and gives the ``presentation'' we choose - invariably we choose
\[ \mu = 1-\sqrt{1+\gamma},\quad \Im m(\mu)>0.\]
 We make this choice so as to give the simplest presentation of the underlying orbifold,  though often it does not matter and both choices of $\mu$ lead to conjugate groups.

There is a simple recipe for  moving from a slope to a word and we give the algorithm for this in the appendix.  For instance the first few Farey slopes (greater then $\frac{1}{2}$) are 
\[ \left\{\frac{1}{2},\frac{4}{7},\frac{3}{5},\frac{5}{8},\frac{2}{3},\frac{5}{7},\frac{3}{4},\frac{4}{5},1\right\} \]
and in the same order the words are,  with $x=X^{-1}$ and $y=Y^{-1}$,
\[ \{ XYxy ,\;  XYxyXYxYXyxYXy ,\;  XYxyXyxYXy,  \;  XYxyXyxYxyXYxYXy, \]
\[  XYxYXy,\;  XYxYXyXyxYxyXy,\;  XYxYxyXy,\; XYxYxYXyXy,\; Xy\}\]

Now let us describe how to read the tables by example. If we look at the $7^{th}$ entry of the $(4,3)$-arithmetic lattices table below we have $X=f$ has order $4$ and $Y=g$ has order $3$.   We see $\gamma(f,g)\approx-2.55804+1.34707 i$,  the complex root with positive real part of its minimal polynomial $z^4+6z^3+13z^2+8z+1$ (there is only ever one conjugate pair of complex roots). The pair $(3/5,2)$ indicates the word $w_{3/5} = XYxyXyxYXy$ is elliptic of order two when $X,Y$ have the form at (\ref{XY})  and $\mu$ has one of the two values solving (\ref{mug}),  in this case
\begin{eqnarray*} \tr w_{3/5} &=& -\mu ^5+\sqrt{5 \sqrt{3}+38} \mu ^4-\left(2 \sqrt{3}+15\right) \mu ^3+\frac{\left(15 \sqrt{3}+7\right) \mu ^2}{\sqrt{2}}\\ && -\left(\sqrt{3}+10\right) \mu +\sqrt{\sqrt{3}+2}\end{eqnarray*}
which is $0$ when $\mu=1.79696 \pm 1.17706 i$. The other possibility for $\mu$ so that (\ref{mug}) holds is $0.652527 - 1.17706 i$ and for that choice $w_{3/5}$ is loxodromic. With minor additional work this gives us a presentation
\[ \langle f,g:f^4=g^3=(fgf^{-1}g^{-1}fg^{-1}f^{-1}gfg^{-1})^2=1 \rangle \]
This group is therefore an arithmetic generalised triangle group in the sense of \cite{HMR},  however it is not identified in that paper as they restrict themselves largely to the case $w_{1/2}=[X,Y]$ is elliptic.   

Next,  we have also given a co-volume approximation in some cases.  This is obtained from adapting the Poincar\'e subroutine in Week's programme Snap to our setting of groups generated by two elliptic elements.  This was previously done by  Cooper \cite{Cooper} in his PhD thesis.  This gives an approximation to the volume of $\IH^3/\langle f,g\rangle$.  However,  this approximation is enough to give the precise index of the the group $\langle f,g\rangle$ in its maximal order and the latter has an explicit volume formula due to Borel \cite{Bo} if further refinement is needed.
 
 \medskip
 
We expect that the topological classification (the structure of the singular locus and the observation that $\gamma(f,g)$ lies on an extended pleating ray) suggested by our data from the arithmetic groups  is also generally true for groups generated by two elements of finite order which are not freely generated.  

\medskip

A major technical advance used in this article is provided by the development of the Keen-Series-Maskit  deformation theory   from the case of the Riley slice (specifically \cite{KS}), to the more general case of groups generated by two elements of finite order,  \cite{EMS1,EMS2}. Using this we are able to further give defined neighbourhoods of pleating rays which lie entirely in these quasiconformal deformation spaces.  This allows us to set up a practical algorithm, which we describe below and implement here,  to decide whether or not a two generator discrete group -- identified solely by arithmetic criteria -- is free on the given generating pair.  And if it is not free on the two generators, then to  identify a nontrivial word in the group.  This process is guaranteed to work if the discrete group is geometrically finite and certain conjectures are true,  and of course the groups we seek will fall into this category,  though the groups we test may not {\em apriori} be geometrically finite.  This gives us an effective computational description of certain moduli spaces akin to the Riley Slice \cite{KS}. We expect that any lattice we might find comes from a bending of a word of some slope that has become elliptic - it's just a matter of finding them.  Since we additionally expect (hope) that the slope is not too big (denominator small), the possibilities might be few and so we simply search for them. However we will see that in fact run into technical difficulties using this approach in a couple of cases where the slops get up to $101/120$ and some guesses are needed to limit the search space (there are $4387$ slopes less than this and we need to precisely generate a polynomial of degree depending on the numerator - and when the degree is large these span pages of course.  Here we are lucky because it turns out that the slopes we seek (as we found out with a lot of experimentation) do not have large integers in their continued fraction expansion. Note that $101/120$ has Conway notation $[6,3,5,1]$ and represents a link of about 15 crossings (we do not know what the crossing number is) and similarly with $7/102$ (Conway $[3,1,1,14]$) probably with about $19$ crossings.  However we do know that $29/51$ below has Conway notation $[7,3,1,1]$ and does have exactly $12$ crossings ($\#762$ in Rolfsen's tables) - it is a surprising group discovered by Zhang \cite{Zhang} in her PhD thesis.  It is quite remarkable that the invariant trace field is at most quartic for all of the above groups.

\medskip

We recall that we require that $f$ and $g$ are primitive elliptic elements,  that is $\tr^2(f)=4\cos^2 \frac{\pi}{4}$ and $\tr^2(g)=4\cos^2 \frac{\pi}{p}$.  This now identifies the Heckoid group from which a complete presentation can be determined if required.  However some of these can be quite long and so we have not given them.

 \subsection{The five $(4,6)$ arithmetic hyperbolic lattices.} 
 \begin{center}
\begin{tabular}{|c|c|c|c|c|} \hline
No. & polynomial & (Farey,order) & $\gamma$ & Discr. $k\Gamma$ \\ \hline
1   & $z+1$ & $Tet(4,6;3)$ & $-1$ & $ -12 $ \\ \hline
2 & $z^2 +6$ & $(7/10,1)$  &  $\sqrt{6} i$ & $-24$ \\  \hline
3 & $z^3 +3z^2 +6z+2$ &  $(5/8,1)$  &  $-1.2980+ 1.8073i$& $-216$  \\ \hline
4 & $ z^3 +5z^2 +9z+3 $ & $(7/12,1)$  & $-2.2873+ 1.3499i$   & $-204$ \\ \hline
5 & $ z^4+8z^2+6z+1$ &   $(17/24,1)$ & $0.3620+2.8764i$& $-2412$  \\ \hline
\end{tabular} \\ 
\end{center}

 \subsection{The (4,5) arithmetic hyperbolic lattice.}  In the case $p=5$ there is one arithmetic lattice and it is cocompact,  $Tet(4,5;3)$.  It has $\gamma=-1$.  This group has presentation 
 \[ \langle x,y:x^4=y^5 =[x,y]^3=(x[yx^{-1}])^2=(x[y,x^{-1}]y)^2=[y,x^{-1}]y)^2=1 \rangle \]
 and is the orientation preserving subgroup of the reflection group with Coxeter diagram $4-3-5$. The number field $k\Gamma$ has discriminant $-400$.
 
  \subsection{The ten $(4,3)$ arithmetic hyperbolic lattices.}
 \begin{center}
\begin{tabular}{|c|c|c|c|c|} \hline
No. & polynomial & (Farey , order)& $\gamma$ & Disc. $k\Gamma$ \\ \hline
1 & $z+2$ & $GT(4,3;2)$ & $ -2 $ & $-24 $ \\ \hline
2 & $z^2 +4$ &  $ (3/4;2)$ & $2i$ & $ -4 $ \\ \hline
3 & $z^2+4 z+6$ &  $ (3/4;1)$ & $-2+\sqrt{2}i$ & $ -8 $ \\ \hline
- & $z^3-z^2+z+1$ &  $ (4/5;4)^*$ & $0.77184 + 1.11514 i$& $-44$ \\ \hline
4 & $z^3 +3z^2 +5z+1$ & $(2/3;2)$  & $-1.38546 + 1.56388 i$& $-140 $  \\ \hline 
5& $ z^3+3z^2+9z+9$ &  $ (17/24;2)$ & $-0.83626+2.46585 i$& $ -108 $ \\ \hline
6 & $ z^3+7z^2+17z+13$ &   $(7/12;1)$ & $-2.77184 + 1.11514 i$& $-44 $ \\ \hline
- & $ z^4+2z^3+3z^2+4z+1$ & $(3/4;3)^*$  & $-0.10176+1.4711 i$& $ -976$ \\ \hline
- & $ z^4-2z^3-z^2+6z+1$ &  $(6/7;3)^*$ & $1.78298 +1.08419 i$& $ -1424 $ \\ \hline 
7 & $ z^4+6z^3+13z^2+8z+1$ & $(3/5;2)$  & $-2.55804+1.34707 i$& $ -3376 $ \\ \hline
8 & $ z^4+6z^3+13z^2+10z+1$ & $(3/5;4)$ \&  $(8/13;3)$ & $-2.20711+0.97831 i $& $ -448$  \\ \hline
9 & $ z^4+4z^3+10z^2+12z+4$ & $(7/10;1)$  & $-1+2.05817 i $& $ -400$ \\ \hline 
10 & $ z^4+8z^3+22z^2+22z+6$ &  $(9/16;1)$  & $-3.11438+0.83097 i $  & $ -2096 $  \\ \hline 
\end{tabular} \\ 
\end{center}  The group $\#8$ above is in fact the index two subgroup of the tetrahedral reflection group $T_5[2,3,3;2,3,4]$, \cite[pg. 228 \& Table 1]{MR2}.  The three groups marked with an asterisk are infinite covolume web-groups with a hyperideal vertex,  or of finite index in the truncated tetrahedral reflection groups enumerated in \cite{CM}. 
  
\subsection{The fifty-four  $(4,2)$ arithmetic hyperbolic lattices.}

 \subsubsection{Real, quadratic.}
 \begin{center}
\begin{tabular}{|c|c|c|c|c|c|} \hline
No. & polynomial & (Farey , order)& $\gamma$ & Disc. $k\Gamma$ & covolume\\ \hline
{\bf Real} & $z+2$ & $GT(4,3;2)$ & $-2$ & $-24$ & \\ \hline
{\bf Quad.} & &&& \\ \hline
2 & $z^2 +2z+2$ &  $ (3/4;2) $ & $-1+ i$ & $ -4$&$ 0.457 $ \\ \hline
3 & $z^2 +z+1$ &  $ (4/5;2) $ & $-0.5 + 0.86602 i$ & $ -3$ & $ 0.253$ \\ \hline
4 & $z^2 +3z+3$ &  $ (7/10;1)$ & $-1.5 + 0.86602 i$& $-3$ & $ 0.253$\\ \hline
5 & $z^2 +1$ &  $ (5/6;2)$ & $i$ & $ -4$ & $ 0.915 $\\ \hline
6 & $z^2 +4z+5$ &  $ (2/3;4)$ & $-2+i$ & $-4 $& $ 0.915 $\\ \hline
- & $z^2 +2z+3$ &  $(3/4;3)$ & $-1+1.41421 i$& $ -8$& $ \infty $ \\ \hline
7 & $z^2 +z+2$ &  $(19/24;1)$  & $-0.5 - 1.32288 i$& $ -7$ & $ 1.332$ \\ \hline 
8 & $z^2 +3z+4$ &  $(17/24;1)$   & $-1.5 +1.32288 i$& $ -7$& $ 1.332 $ \\ \hline
9 & $z^2 +5z+7$ &  $(19/30;1) $  & $-2.5+ 0.86602 i$ & $-3 $ & $ 1.524$\\ \hline
10 & $z^2-z+1$ &  $(13/15;2)$  & $0.5+ 0.86602 i$ & $-3 $ & $ 1.524$ \\ \hline   
\end{tabular} \\ 
\end{center} 
\newpage

 \subsubsection{Cubic.}
\begin{center}
\begin{tabular}{|c|c|c|c|c|c|} \hline
No. & polynomial & (Farey , order)& $\gamma$ & Disc. $k\Gamma$ & covolume\\ \hline      
 11 & $z^3-2z+2$ &  $(43/48;1)$ & $0.884646 + 0.58974 i$ & $-76 $ & $ 1.985$\\ \hline 
 12 & $z^3+6z^2+10z+2$ &  $(39/40;4)$  & $-2.884646 + 0.58974 i$ & $-76 $ & $ 1.985 $\\ \hline  
13 & $z^3-z+1$ &  $(8/9;2)$ & $0.662359 - 0.56228 i$ & $-23 $ & $ 0.824$ \\ \hline  
14 & $z^3+6z^2+11z+5$ &  $(11/18;1)$ & $-2.66236 + 0.56228 i$ & $ -23 $ & $ 0.824$  \\ \hline  
15 & $z^3+z^2-z+1$ &  $(7/8;3)$ & $0.419643 - 0.60629 i$ & $ -44$ & $ 0.264$\\ \hline  
16 & $z^3+5z^2+7z+1$ &  $(5/8;3)$ & $-2.41964 - 0.60629 i$ & $-44 $ & $ 0.264$ \\ \hline  
17 & $z^3+z^2+1$ &  $(6/7;2)$   & $0.232786 + 0.79255 i$ & $-31 $ & $0.595$\\ \hline 
18 & $z^3+5z^2+8z+3$ &  $(9/14;1)$ & $-2.23279+ 0.79255 i$ & $ -31$ & $0.595$\\ \hline 
19 & $z^3+z^2+2$ &  $(17/24;1)$  & $0.34781+1.02885i$ & $-116$ & $ 2.232 $ \\ \hline 
20 & $z^3+5z^2+8z+2$ &  $(17/24;1)$ & $-2.34781 +1.02885 i$ & $-116 $ & $ 2.232 $ \\ \hline 
21 & $z^3+2z^2+z+1$ &    $(5/6;3)$  & $-0.122561 + 0.74486 i$ & $ -23$& $ 0.137$  \\ \hline 
22 & $z^3+ 4 z^2 + 5 z+1$ &    $(2/3;3)$   & $-1.87744 - 0.74486 i$ & $ -23$& $ 0.137$ \\ \hline 
23 & $z^3+z^2+z+2$ &    $(101/120;1)$ & $0.17660+1.20282 i $& $-83 $ & $ 3.308$\\ \hline
24 & $z^3+5z^2+9z+4$ &   $(79/120;1)$  & $-2.1766+ 1.20282i$ & $-83 $ & $ 3.308$ \\ \hline
25 & $z^3+2z^2+2z+3$ &    $(49/60;1)$ & $-0.09473 + 1.28374 i$ & $ -139$& $ 2.597  $ \\ \hline
26 & $z^3+4z^2+6z+1$ &    $(41/60;1)$  & $-1.90527 + 1.28374 i$ & $-139 $ & $ 2.597 $ \\ \hline
27 & $z^3+2z^2+2z+2$ &   $(13/16;1)$ & $-0.22815+1.11514 i$ & $ -44$ &  $0.793$ \\ \hline
28 & $z^3+4z^2+6z+2$ &   $(11/16;1)$ & $-1.77184 + 1.11514 i$ & $ -44$& $ 0.793 $ \\ \hline
29 & $z^3+z^2+2z+1$ &   $(17/21;2)$ & $-0.21508 + 1.30714 i$& $-23 $ & $ 0.824$ \\ \hline
30 & $z^3+5z^2+10z+7$ &   (29/42;1)  & $-1.78492+ 1.30714 i$ & $ -23 $& $ 0.824 $\\ \hline
31 & $z^3+2z^2+3z+1$ &   $(7/9;2)$, & $-0.78492 + 1.30714 i$ & $ -23$ & $ 0.824 $ \\ \hline
32 & $z^3+4z^2+7z+5$ &   $(13/18;1)$, & $-1.21508 + 1.30714 i$ & $-23 $& $ 0.824 $  \\ \hline 
33 & $z^3+3z^2+5z+4$ &   $ (31/40;1)$ & $-0.773301 + 1.46771 i$ & $ -59$ & $ 1.927 $\\ \hline
34 & $z^3+3z^2+5z+2$ &   $ (29/40;1)$ & $-1.2267 + 1.46771 i$ & $ -59$ & $ 1.927$ \\ \hline
35 & $z^3+3z^2+4z+3$ &  $(11/14;1)$ & $-0.658836-1.16154 i$ & $-31 $ & $0.593$  \\ \hline
36 & $z^3+3z^2+4z+1$ &   $(5/7;2)$ & $-1.34116 + 1.16154 i$ & $-31 $ & $0.593$ \\ \hline
\end{tabular} \\ \end{center} 

 \newpage
\subsubsection{Quartic.}
\begin{center}
\begin{tabular}{|c|c|c|c|c|c|} \hline
No. & polynomial & (Farey , order)& $\gamma$ & Disc. $k\Gamma$ & covol. \\ \hline
 \\ \hline
 37& $ z^4 -  4 z^2 +  z + 3 $ &    $(7/102;1)$ & $1.36778 +  0.23154 i $ & $ -731$  & $ 2.9003 $  \\ \hline
   38 & $ 1 + 15 z + 20 z^2 + 8 z^3 + z^4 $ &    $(29/51;2)$ & $-3.36778+  0.23154 i $ & $ -731$  & $ 2.9003 $  \\ \hline
 39& $1 + z - 2 z^2 + z^4$ &    $(77/85;2)$ & $1.00755 + 0.51311i $ & $ -283$  & $ 3.6673 $  \\ \hline
40 & $ 7 + 23 z + 22 z^2 + 8 z^3 + z^4 $ &    $(101/170;1)$ & $-3.00755  + 0.51311i $ & $ -283$  & $ 3.6673 $  \\ \hline
41 & $z^4+z^3-3z^2-z+3$ &    $(31.34;1)$ & $1.06115 + 0.38829 i$ & $ -491 $ & $ 1.542$ \\ \hline
42 & $z^4+7z^3+15z^2+9z+1$ & $(10/17;2) $ & $-3.06115+ 0.38829 i$ & $  -491$ & $ 1.542$\\ \hline
- & $z^4+4z^3+7z^2+6z+1$ & $(3/4;5) $   & $-1+ 1.27202 i$ & $ -400 $ & $\infty$ \\ \hline
- & $z^4+4z^3+8z^2+8z+2$ & $(3/4;4)  $ & $-1+ 1.55377 i$ & $ -1024 $ & $\infty$ \\ \hline
- & $z^4+4z^3+8z^2+8z+1$ & $(3/4;6)  $ & $-1+1.65289 i$ & $ -4608 $ & $\infty$  \\ \hline
43 & $z^4+5z^3+11z^2+11z+3$ & $(40/51;2)$ & $-1.42057+ 1.45743 i$ & $ -731$ & $ 2.889 $\\ \hline
44 & $z^4+3z^3+5z^2+5z+1$ & $(40/51;2)$ & $-0.57943+ 1.45743 i$& $ -731$ & $ 2.889$ \\ \hline
45 & $z^4+2z^3+2z^2+3z+1$ & $(14/17;2)$ & $0.03640 +1.21238 i$ & $-491 $ & $ 1.540 $\\ \hline
46 & $z^4+6z^3+14z^2+13z+3$ &  $(23/34;1)$  & $-2.03640 +1.21238 i$ & $ -491$& $ 1.540$ \\ \hline
47 & $z^4+6z^3+13z^2+10z+1$ & $(13/20;1)$ & $-2.20711 +0.97831 i$ & $ -448$ & $ 1.028 $ \\ \hline
48 & $z^4+2z^3+z^2+2z+1$ & $(17/20;1)$  & $0.20711 +0.97831 i$ & $ -448$ & $ 1.028 $ \\ \hline
49 & $z^4+z^3-z^2+z+1$ & $(15/17;2)$    & $0.65138 + 0.75874 i$ & $-507 $ & $ 1.624 $ \\ \hline
50 & $z^4+7z^3+17z^2+15z+3$ &  $(21/34;1)$   & $-2.65139 + 0.75874 i$ & $ -507$& $ 1.624$  \\ \hline
51 & $ z^4+3 z^3+4 z^2+4 z+1 $ &  (4/5;3) &	$ -0.40630 +  1.19616 i $ &$-331$ &  $  0.447  $\\ \hline
52 & $  z^4+5 z^3+10 z^2+8 z+1 $ & (7/10;3) &	 $ -1.59369 +  1.19616 i $   & $-331$&$ 0.447 $ \\ \hline
53 & $ z^4 + z^3 - 2 z^2+ 1 $ &  (9/10;3) &	$  0.78810 +  0.40136 i $ &$-283$ & $ 0.347$ \\ \hline
54 & $  z^4+7 z^3+16 z^2+12 z+1 $ & (3/5;3) &	 $ -2.78810 + 0.40135i$   & $-283$&$ 0.347 $\\ \hline
\end{tabular} 
\end{center} 
\newpage

 \subsection{The twenty-seven $(4,4)$ arithmetic hyperbolic lattices.}
  
 \begin{center}
\begin{tabular}{|c|c|c|c|c|c|} \hline
No. & polynomial & (Farey,order)& $\gamma$ & Disc. $k\Gamma$ & covol.\\ \hline
 \\ \hline 
1 & $z+2$ &  $GT(4,4;2)$ & $-2$ & $-4$ &$0.91596$ \\ \hline
-& $z+3$ &  $Tet(4,4;3)$ & $-3$ & $-8$ & $\infty$ \\ \hline
-& $1 + 3 z + z^2$ & $Tet(4,4;5)$ & $-2.61803$ & $ 5 $ & $\infty$ \\ \hline
-& $2 + 4 z + z^2$ & $GT(4,4;4)$ & $-3.41421$ & $ 8$ & $\infty$ \\ \hline
-& $1 + 4 z + z^2$ & $GT(4,4;6) $ & $-3.73205$ & $12$ & $\infty$  \\ \hline
 \\ \hline 
2 & $3 + 3 z + z^2$ &  $ (3/5,1) $  & $-1.5 + 0.866025 i$ & $ -3$ & $0.5067$ \\ \hline
3& $5 + 2 z + z^2$ &  $ (2/3;2) $ & $-1 + 2 i$ & $ -4$ & $1.8319 $ \\ \hline
4& $8 + 5 z + z^2$ &  $(7/12;1)$ & $-2.5 + 1.32288 i$ & $ -7$  & $2.6648$ \\ \hline 
5& $7 - z + z^2$ &    $(11/15;1)$  & $0.5 + 2.59808 i$ & $-3 $ & $3.0491$ \\ \hline

  \\ \hline
6& $3 + 8  z + 5 z^2 + z^3$ &   $(4/7;1)$ & $-2.23278+0.79255 i$ & $ -31$  & $1.1878$ \\ \hline 
7& $4 + 8 z - 4 z^2 + z^3$ &  $(19/24;1)$ & $2.20409 + 2.22291 i$ & $-76$  & $3.9710$ \\ \hline 
8& $5 + 3 z - 2 z^2 + z^3$ &  $(7/9;1)$ & $1.44728 +1.86942 i$ & $-23 $  & $1.6481$ \\ \hline  
9& $1 + 3 z - z^2 + z^3$ &  $(3/4;3)$ & $0.647799 + 1.72143 i$ & $-44$ & $0.5295$  \\ \hline  
10& $3 + 4 z + z^2 + z^3 $ &  $(5/7;1)$  & $-0.108378 + 1.95409 i$ & $-31$ & $1.1900$ \\ \hline 
11& $4 + 8 z + z^2 + z^3$ &  $(17/24;1)$ & $-0.241944 + 2.7734 i $& $-116$ & $4.4654$  \\ \hline 
12& $1 + 3 z + 2 z^2 + z^3$ &    $(2/3,3)$  & $-0.78492 + 1.30714 i$ & $ -23$ & $0.2746$  \\ \hline 
13& $8 + 11 z + 3 z^2 + z^3$ &    $(41/60;1)$ & $-1.06238 + 2.83049 i$& $-83$ & $6.6173$  \\ \hline
14& $3 + 10 z + 4 z^2 + z^3$ &    $(19/30;1)$ & $-1.82848 + 2.32426 i$ & $ -139$  & $5.1943$ \\ \hline
15& $4 + 8 z + 4 z^2 + z^3$ &   $(5/8;1)$ & $-1.6478 + 1.72143 i $& $-44$  & $0.1322 $ \\ \hline
16& $7 + 12 z + 5 z^2 + z^3$ &   $(5/9;1)$ & $-2.09252 + 2.052 i$& $-23 $  & $4.2441$  \\ \hline
17& $8 + 15 z + 7 z^2 + z^3$ &   $(11/20;1)$ & $-3.10278 + 0.665457 i$ & $ -59$  & $3.8557$ \\ \hline
18& $3 + 8  z + 5 z^2 + z^3$ &   $(4/7;1)$ & $-2.23278+0.79255 i$ & $ -31$  & $1.1878$ \\ \hline 
\end{tabular} 
\end{center} 
 \begin{center}
\begin{tabular}{|c|c|c|c|c|c|} \hline
No. & polynomial & (Farey,order)& $\gamma$ & Disc. $k\Gamma$ & covol.\\ \hline
 \\ \hline 

19& $z^4-8 z^3+12 z^2+23 z+3$ &    $(44/51;1) $ & $4.55279 + 1.09648 i$ & $-731$  & $5.8007 $  \\ \hline
20& $7 + 15 z + 4 z^2 - 4 z^3 + z^4$ &    $(69/85;1) $ & $2.76698 + 2.06021 i$ & $-283$  & $7.3346 $  \\ \hline
21& $3 + 13 z + 5 z^2 - 5 z^3 + z^4$ &    $(14/17;1) $ & $3.09755 + 1.60068 i$ & $-491$  & $3.0852$  \\ \hline
22& $3 + 13 z + 17 z^2 + 7 z^3 + z^4$ & $(29/51;1$ & $-2.94722 + 1.22589 i$ & $ -731$  & $5.7795$ \\ \hline
23& $3 + 11 z + 12 z^2 + 4 z^3 + z^4$ & $(11/17;1)$   & $-1.39573 + 2.51303 i$ &   $-491 $  & $3.0815$ \\ \hline
24& $1 + 6 z + 7 z^2 + 2 z^3 + z^4$ & $(7/10;1)$   & $-0.5 + 2.36187 i$ & $-448$  & $2.0575$ \\ \hline
25& $3 + 9 z + 5 z^2 - z^3 + z^4$ & $(13/17;1)$  & $1.15139 - 2.50596 i$ & $-507$  & $3.2480$ \\ \hline
26& $z^4+5 z^3+10 z^2+6 z+1$ & $(3/5;3)$   & $-2.07833 + 1.4203 i$ & $-331$  & $0.8948$ \\ \hline
27& $z^4-3 z^3+2 z^2+6 z+1$ & $(4/5;3)$   & $2.03623 + 1.43534 i$ & $-283$  & $0.6951$ \\ \hline
\end{tabular} 
\end{center}  
 
 For completeness we recall Takeuchi's result in this case.
\begin{theorem} Let $\Gamma=\langle f,g \rangle$ with $f$ and $g$ of finite order $4$ and $p$ be an arithmetic Fuchsian triangle group Then there are eleven such groups and $\{4,p,q\}$ form one of the following triples.
\begin{itemize}
\item noncompact : $(2,4,\infty)$, $(4,4,\infty)$, 
\item compact: $(2,4,12)$, $(2,4,18)$, $(3,4,4)$, $(3,4,6)$, $(3,4,12)$, $(4,4,4)$, $(4,4,5)$, $(4,4,6)$, $(4,4,9)$.  
\end{itemize}
\end{theorem}

 \section{Total degree bounds.}
 
 In this section we first outline the arithmetic criteria on the complex number 
 \[ \gamma=\gamma(f,g) = \tr[f,g]-2 \]
  that are necessary in order for the group $\Gamma=\langle f,g:f^4=g^p=
 \cdots=1\rangle$ to be a discrete subgroup of an arithmetic group, following on from \cite{MM,MMpq}. This is encapsulated in the Identification Theorem \ref{idthm}  below. It places rather stringent conditions on $\gamma=\gamma(f,g)$.  
 
 Next,  in order for $\Gamma$ to be a lattice,  it is necessary that is it not free on generators.  We then use the ``disjoint isometric circles test'' and the Klein ``ping-pong'' lemma to obtain a coarse bound on $|\gamma|$.  The arithmetic criteria and coarse bound then allow us to either directly give degree bounds or,  in some cases, adapt the method of Flammang and Rhin \cite{FR} to obtain a total degree bound for the field $k\Gamma=\IQ(\gamma)$.  It will turn out that the degree is actually no more than $4$ as we have noted above,  however a total degree bound of $7$ or $8$  brings us down into the range of feasible search spaces for integral monic polynomials with the root bounds we obtain exploiting arithmeticity and bounds on $|\gamma(f,g)|$,  and so this is the first bound we will seek.   Note that this gives first bounds on $p$ since $\IQ(\cos\frac{2\pi}{p})\subset k\Gamma$ and 
 \[ [\IQ:k\Gamma]=[\IQ:\IQ(\cos\frac{2\pi}{p})][\IQ(\cos\frac{2\pi}{p}):k\Gamma]\leq 8.\]

Then refining this degree bound down to $4$ or $5$ relies on us obtaining a coarse description of the moduli space ${\cal M}_{4,p}$ --- the deformation space of Riemann surfaces $\IS_{4,4,p,p}$ topologically homeomorphic to the sphere with two cone points of order $4$ and two cone points of order $p$,  an analogue of the Riley slice for groups generated by two parabolics where $\IS_{4,4,p,p}$ is replaced by the four times punctured sphere.     Actually only the cases $(4,3)$,$(4,3)$ and $(4,4)$ have very  large search spaces for (as we shall see below)  there is not necessarily an intermediate field real for us to work with.  Indeed the $(4,2)$ case is actually covered by the methods of Flammang and Rhin \cite{FR} after they achieved a degree $9$ bound (for the $p=q=3$ case) and their methods largely carry over but are a bit more complicated since these moduli spaces when $q=4$ and $p\neq 2,4$ are not symmetric and do have  significantly larger complements (where the groups we are looking for will lie).  Next, the arithmetic criteria show that $\IQ(\gamma)$ has one complex conjugate pair of embeddings and good bounds on the real embeddings of the number $\gamma$ due to a ramification requirement.  Thus the coefficients of the minimal polynomial for $\gamma$ are bounded with a bound on $|\gamma|$ (though very crudely).  Given the degree bound we can now look through the space of polynomials with integer coefficients within these ranges to produce a list of possibilities.  Without good bounds this can take quite a long time (many days).  Then a factorisation criteria significantly reduces the possibilities further.  This leaves us with a moderate number of values for $\gamma$ (perhaps several thousand).  At this point there is no issue about whether the group is discrete.  That is guaranteed.  

The only question is if it is a lattice and if so what is it.  This is where we need the refined descriptions of moduli spaces. Our earlier work  implemented a version of the Poincar\'e polyhedral theorem to construct a fundamental domain for the group acting on hyperbolic $3$-space and thereby hope to determine if this action is cocompact or not \cite{Cooper}. The case of cofinite volume but not cocompact is rather easier in the arithmetic case as then the invariant trace field $k\Gamma$ is quadratic \cite{MR}.  This implementation worked well when there were no points quite close to the boundary (where the geometrically infinite groups are dense). However it fell over for points close to the boundary, and  was a case by case analysis which was fine when there were only a few tens of groups to look at. Here we replace the last few steps by a finer description - and so tighter bounds -  on the set where $\gamma$ must lie,  giving smaller search spaces (a matter of minutes). Then we ``guess'' it is a Heckoid group and then set about finding what it might be by enumerating the possible words and associated trace polynomials.  This just happens to work --- and of course it is natural to conjecture that is always does.

 \subsection{The standard representation.}  Define the two matrices 
\begin{equation}\label{AB} A = \begin{pmatrix}  \cos \pi/p & i \sin \pi/p \\ i \sin \pi/p & 
\cos \pi/p \end{pmatrix}, \quad B = \begin{pmatrix}
 \cos \pi/q & iw \sin\ \pi/q \\ i w^{-1} \sin \pi/q & \cos \pi/q \end{pmatrix}. 
 \end{equation}
Then if $\G = \langle  f,g\rangle $ is a non-elementary Kleinian group with $o(f)=p$ (where $o(f)$ denotes the order of $f$) and 
$o(g) = q$, where $p \geq q \geq 3$ then $\G$ can be normalised so that 
$f,g$ are represented by the matrices $A,B$ respectively. The parameter
$\gamma$ is related to $w$ by
\begin{equation}\label{2}
\gamma = \sin^2 \frac{\pi}{p} \, \sin^2 \frac{\pi}{q} \; (w - \frac{1}{w})^2.
\end{equation}
Given $\gamma$, we can further normalise and choose $w$ such that $| w | 
\leq 1$ and ${\rm Re}(w) \geq 0$.

\subsection{Isometric circles.}

The isometric circles of the M\"obius transformations determined by (\ref{AB}) are 
\[ \{z:|i \sin(\frac{\pi}{p}) z \pm  \cos(\frac{\pi}{p}) |^2=1\}, \;\;\; \{ |i w^{-1} \sin(\frac{\pi}{q}) z \pm \cos (\frac{\pi}{q}) |^2=1\} \]
These circles are paired by the respective M\"obius transformations.
With our normalisation on $w$ as above, these two pairs of circles are disjoint precisely when
\begin{equation}\label{3}
| i w \cot \pi/q + i \cot \pi/p | + \frac{| w |}{\sin \pi/q} \leq \frac{1}{\sin \pi/p}.
\end{equation}
If $\gamma$ is real,   the isometric circles are disjoint if and only if $\gamma<-4$ or $\gamma > 4\big(\sqrt{2} \cos(\frac{\pi }{p})+1\big)+2\cos(\frac{2 \pi }{p})$.  In the case of equality here the isometric circles are tangent.

\bigskip

With our normalisations, the well known Klein ``ping-pong" argument quickly implies that, should the isometric circles be pairwise disjoint  or at worse tangent,  then the group generated by $A$ and $B$ is free on these two generators.  It therefore cannot be an arithmetic lattice. In fact in these circumstances (with $p=4$ and $q\geq 3$) it is easy to see
\[ \left[ \ID(i,\sqrt{2})\cap\ID(-i,\sqrt{2}) \right] \setminus \left[\ID\Big(i w\cot (\frac{\pi}{q}),\frac{|w|}{\sin(\frac{\pi}{q})}\Big) \cup \ID\Big(-i w\cot (\frac{\pi}{q}),\frac{|w|}{\sin(\frac{\pi}{q})}\Big)\right] \]
is a nonempty  fundamental domain for the action of $\Gamma=\langle f,g\rangle$ on $\oC\setminus \Lambda(\Gamma)$.

\medskip
 We now need to turn the inequality (\ref{3}) into a condition on $\gamma$.  From (\ref{3}) we see with $q=4$ and $w=r e^{i\theta}$,  $0\leq \theta\leq \pi/2$, that 
\begin{eqnarray*}  |  r \cos(\theta) +  \cot \pi/p +i r\sin(\theta)| + \sqrt{2}r & \leq &\frac{1}{\sin \pi/p} \\
\Big( r \cos(\theta) +  \cot (\frac{\pi}{p})\Big)^2 +  r^2\sin^2(\theta) & \leq &(\frac{1}{\sin\frac{\pi}{p}}- \sqrt{2}r )^2 \\
2 r \cos(\theta) \cot (\frac{\pi}{p}) & \leq &1+r^2 - \frac{2\sqrt{2}r}{\sin\frac{\pi}{p}}  \\
\end{eqnarray*}
On the curve where equality holds,  we obtain a parametric equation for a region $\Omega_p$ bounded by this curve. 
\[\Omega_p\; :\;0=1+r^2 -2 r \; \frac{\sqrt{2}- \cos(\theta) \cos (\frac{\pi}{p})}{\sin\frac{\pi}{p}}\]

 \scalebox{0.5}{\includegraphics[viewport=-50 250 600 750]{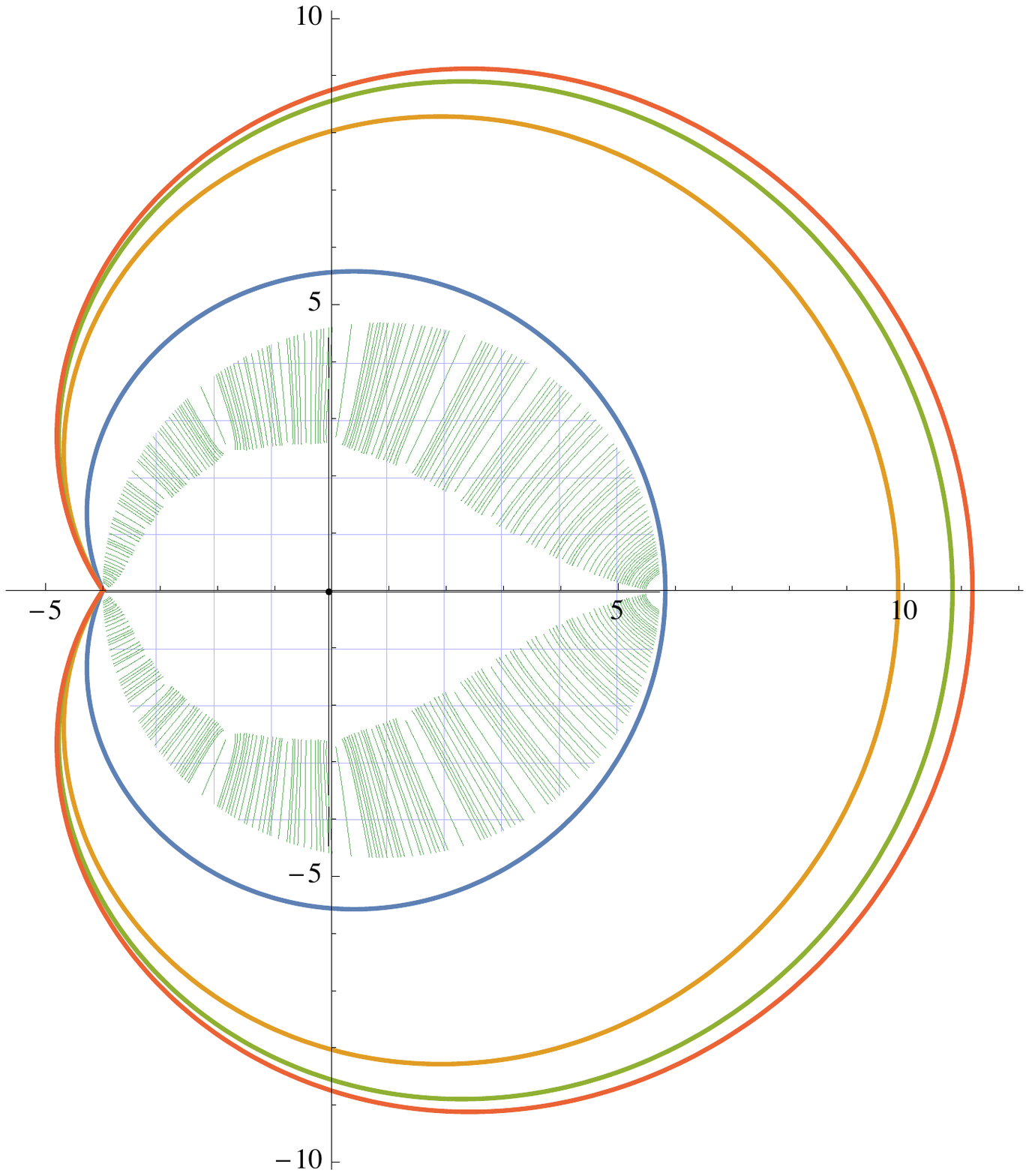}}\\
 {\bf Figure 2.} {\em The smooth parametric regions $\Omega_p$ from smallest to largest,  $p=3,6,9$ and $12$.  If $\gamma \not\in \Omega_p$,  then $\langle f,g\rangle$ is freely generated. Inset is the actual region (in the case $p=3$) where free generation occurs,  illustrated by interior pleating rays landing on its boundary.}
 
 \medskip
 
The $\gamma$ values for groups we seek are inside the curves $\Omega_p$.  However the  regions $\Omega_p$ where the simplest criterion giving a region where disjoint isometric circles holds is apparently much larger than that required to identify groups that are freely generated.   However these regions are sufficient to provide initial bounds especially for $p$ large enough that $[\IQ(\cos 2\pi/p):\IQ]>17$,  that is $31$ possible values no more than $60$.

 \bigskip
\subsection{Real points.} At this point we would like to remove the cases $\gamma\in \IR$ from further consideration.  This result is primarily based on \cite{MM3}

\subsubsection{If $\gamma>0$,}  then the group $\Gamma$ is a Fuchsian triangle group,  or is freely generated,  and is therefore not an arithmetic lattice.  However, the $(4,p,q)$- triangle groups $\langle a,b|a^4=b^p=(ab)^q \rangle$ have parameters 
\begin{eqnarray*}(\gamma,\tr^2(f)-4,\tr^2(g)-4)& = & (\gamma,-2,-4\sin^2(\frac{\pi}{p})), \\  \gamma &=& 4 \cos (\frac{\pi }{q})\big(\sqrt{2} \cos(\frac{\pi }{p})+\cos(\frac{\pi }{q})\big)+2\cos(\frac{2 \pi }{p})\end{eqnarray*} 
\subsubsection{If $\gamma<0$,}   and if  $\langle f,g\rangle$ is not free on generators,  then $-4<\gamma <0$.  It follows that $-2< \tr[f,g] < 2$ and $[f,g]$ is elliptic.  All groups generated by two elements of finite order and whose commutator has finite order were identified in \cite{MM3}.  As examples, the $(4,p,q)$-generalised triangle groups $\langle a,b|a^4=b^p=[a,b]^q \rangle$ have parameters
\[ (\gamma,-2,-4\sin^2(\frac{\pi}{p})),  \gamma =-2-2\cos(\frac{\pi}{p}) \]  

\bigskip
From that classification we have the following result.
\begin{theorem} Let $\Gamma=\langle f,g\rangle$ be a non-elementary Kleinian group where $f,g$ are
elliptic elements with $o(f) =4$, $o(g) = q$ and $o([f,g]) = n$.   Assume that $\gamma$ is nonelementary and does not have an invariant
hyperbolic plane. Then the generators give rise to the sets of parameters described are
\begin{enumerate}
\item {\em Generalised Triangle groups}  $\langle x,y|x^4=y^p=[x,y]^n=1\rangle$. 
\[ (\gamma,-2,\beta)= (-2-2 \cos(\pi/n),-2,-4 \sin^2(\pi/p)). \]
\item {\em Tetrahedral groups with $n$ odd}. \begin{eqnarray*} 
\langle x,y|x^4 & = & y^p=[x,y]^n=1, (x[y, x^{-1}]^{(n-1)/2})^2 = 1, \\
&& (x[y, x^{-1}]^{(n-1)/2}y)^2 = ([y, x^{-1}]^{(n-1)/2}y)^2 =1\rangle
\end{eqnarray*}
\[ (\gamma,-2,\beta)= (-2-2 \cos(2\pi/n),-2,-4 \sin^2(\pi/p)). \]
\item   {\em Tetrahedral groups with $n$ odd and $n\geq 7$}. In this case we have a group with the presentation as determined above,  but additional possibilities for the value of $\gamma$.
\[ (\gamma,-2,\beta)= (-2-2 \cos(4\pi/n),-2,-4 \sin^2(\pi/p)). \]
\end{enumerate} \
The only cocompact group is the group  $(4,p;n)=(4, 5; 3)$  which is arithmetic. The only non-cocompact groups of finite
covolume are $(4,p;n)=GT(4,3; 2), Tet(4, 6; 3), GT(4,4,2)$,
which are all arithmetic. 
\end{theorem}
Futher data on these groups can be found in \cite{MM3}.

\medskip

Notice here that if $n=2$,  then $tr[f,g]=0$ and $\gamma(f,g)=-2=\beta(f)$.  Then $\gamma(f,gfg^{-1})=\gamma(f,g)(\gamma(f,g)-\beta(f))=0$ and $f$ and $gfg^{-1}$ share a fixed point on $\oC$.  In this last case here $(4,p;n)=(4, 4; 2)$ is the $\IZ_2$ extension with parameters $(1+i,-2,-4)$.  In summary we have the possible values for $\gamma\in \IR$ for arithmetic Kleinian groups.
\[
p=5; \gamma=-3, \quad
p=3 ; \gamma=-2,\quad
p=4 ; \gamma=-2\quad
p=6 ; \gamma=-3
\]

\subsection{Arithmetic restrictions on the commutator $\gamma(f,g)$.}

Having dispensed with the case that $\gamma=\gamma(f,g)\in \IR$ in the rest of this paper we shall assume that $\gamma$ is complex.
We require the following preliminaries. Let
$\G$ be any non-elementary finitely-generated subgroup  of $\PSL(2,\IC)$. 
Let $\G^{(2)} = \langle  g^2 \mid g \in \G \rangle $ so that $\G^{(2)}$ is a subgroup of finite
index in $\G$. Define
\begin{equation} \left.
\begin{array}{lll}
k\G  & = & \IQ(\{ \tr(h) \mid h \in \G^{(2)} \}) \\ 
A\G  & = & \{ \sum a_i h_i \mid a_i \in k\G, h_i \in \G^{(2)} \}
\end{array}\;\;\;\;  \right\}
\end{equation}
where, with the usual abuse of notation, we regard elements of $\G$ as matrices,
so that $A\Gamma \subset M_2(\IC)$. 

Then $A\G$ is a quaternion algebra over $k\G$ and the pair $(k\G, A\G)$ is 
an invariant of the commensurability class of $\G$. If, in addition, $\G$ 
is a Kleinian group of finite co-volume,  then $k\G$ is a number field.

We state  the identification theorem as follows:
\begin{theorem} \label{idthm} Let $\G$ be a  subgroup of 
$\PSL(2,\IC)$ which is finitely-generated and non-elementary. Then $\G$
is an arithmetic Kleinian group if and only if the following conditions all hold:
\begin{enumerate}
\item  $k\G$ is a number field with exactly one complex place,
\item for every $g \in \G$, $\tr(g)$ is an algebraic integer,
\item $A\G$ is ramified at all real places of $k\G$.
\item $\G$ has finite co-volume.
\end{enumerate}
\end{theorem}
It should be noted  that the first three conditions together imply that $\G$ is 
Kleinian, and without the fourth condition, are  sufficient to imply that $\G$ is a subgroup of an arithmetic Kleinian group.

The first two conditions clearly depend on the traces of the elements of $\G$.
In addition, we may also find a Hilbert symbol for $A\G$ in terms of the 
traces of elements of $\G$ so that the third condition also depends on the traces 
(for all this, see \cite{MR1},\cite[Chap. 8]{MR}).

 \subsection{The factorisation condition.}\label{factorisation}
 
 For an arithmetic group generated by elliptic elements of order $4$ and $p$  \cite{MMpq} observed that as a consequence of the Fricke identity the number $\lambda= \tr f \tr g \tr fg$ is an algebraic integer and satisfies the quadratic equation
\begin{equation}\label{eqn5}
x^2 - 8\cos^2\frac{\pi}{p}\, x 
  + 8\cos^2\frac{\pi}{p}(-2 + 4\cos^2\frac{\pi}{p} -\gamma) = 0,
\end{equation}
Thus 
\[ \lambda = 4 \cos ^2 \frac{\pi }{p} \pm 2 \sqrt{2} \cos \frac{\pi }{p}  \sqrt{\gamma + 2\sin^2\frac{ \pi }{p}} \]
and 
\begin{equation}\label{lp}
 \lambda_p=  2 \sqrt{2} \cos \frac{\pi }{p}  \sqrt{\gamma + 2\sin^2\frac{ \pi }{p}} 
\end{equation}
is an algebraic integer.
Further, if $p$ is odd,  then $2 \cos \frac{\pi }{p} $ is a unit and so 
\[ \alpha_p =   \sqrt{2}   \sqrt{\gamma +2\sin^2\frac{ \pi }{p}} =\sqrt{2\gamma -\beta} \]
is an algebraic integer. Since the field $\IQ(\gamma)$ has one complex place we must have $\IQ(\gamma)=\IQ(\lambda_p)$ and when $p$ is odd $\IQ(\gamma)=\IQ(\alpha_p)$.

\medskip

This factorisation criterion yielding $\IQ(\gamma)=\IQ(\lambda_p)$ is a powerful obstruction for an algebraic  integer $\gamma$ to satisfy when $p\neq 2$.   In particular we will apply it in the following form which is easily seem to be equivalent in the case at hand in fields with one complex place.
\begin{theorem} \label{samedegree} Let $\Gamma=\langle f,g\rangle$ be an arithmetic Kleinian group generated by elliptic elements of order $4$ and $p$,  with $p\neq 2$.  Then the minimal polynomial for $\lambda_p$ as at (\ref{lp}) has the same degree as the minimal polynomial for $\gamma=\gamma(f,g)$.
\end{theorem} 

For $p\not\in\{2,3,4,6\}$ there is an intermediate real field as $\cos(2\pi/p)\not \in\IQ$ lies in the invariant trace field.  We set
\[ L = \IQ(\cos(2\pi/p)),\quad\quad [\IQ(\gamma):\IQ]=[L:\IQ] \times [\IQ(\gamma):L] \]

Now if $\G$ is arithmetic,  $k\G$ will be  a number field with one complex place if and only if $L(\gamma)$
has one complex place and the quadratic at (\ref{eqn5}) splits into linear factors over
$L(\gamma)$. This implies that, if $\tau$ is any real embedding of 
$L(\gamma)$, then the image of the discriminant of  (\ref{eqn5}), which is 
$32 \cos^2\frac{\pi}{p}(-2\sin^2\frac{\pi}{p} +   \gamma)$, 
under $\tau$ must be positive. Clearly this is equivalent to requiring that
\begin{equation}\label{eqn6}
\tau( 2\cos^2\frac{\pi}{p} + \gamma) > 0.
\end{equation}
Thus $k\G$ has one complex place if and only if (i) $\IQ(\gamma)$ 
has one complex 
place, (ii) $L \subset \IQ(\gamma)$, (iii) for all real embeddings $\tau$ of $\IQ(\gamma)$,
(\ref{eqn6}) holds and (iv) the quadratic at (\ref{eqn5}) factorises over $\IQ(\gamma)$.

 Now, still in the cases where $p> 2$,(\cite[\S 3.6]{MR})  
\begin{equation}\label{eqn7}
A\G = \HS{-1}{ 2\cos^2(\frac{\pi}{p}) \,\gamma}{k\G}.
\end{equation}
Under all the real embeddings of $k\G$,  $A\G$ is ramified at
all real places of $k\G$ if and only if, under any real embedding $\tau$
of $k\G$,
\begin{equation}\label{eqn8}
\tau(\gamma) < 0.
\end{equation}
Thus, summarising, we have the following theorem which we will use to determine the possible $\gamma$ values for the groups we seek.
\begin{theorem}[The Identification Theorem] \label{2genthm}
Let $\G = \langle  f , g \rangle $ be a non-elementary
subgroup of $\PSL(2,\IC)$ with $f$ of order $4$ and $g$ of order $p$, $p \geq 3$.  Let $\gamma(f,g) = \gamma \in \IC \setminus \IR$. Then $\G$ is an arithmetic Kleinian group if and only
if 
\begin{enumerate}
\item $\gamma$ is an algebraic integer,
\item $\IQ(\gamma) \supset L = \IQ(\cos 2 \pi/p)$ and $\IQ(\gamma)$ is a number field with exactly one complex place,
\item if $\tau : \IQ(\gamma) \rightarrow \IR$ such that $\tau |_L = \sigma$, then
\begin{equation}\label{eqn10}
 - \sigma( 2\sin^2 \frac{\pi}{p}) < \tau(\gamma) < 0,
\end{equation}
\item the algebraic integer $\lambda_p$ defined at (\ref{lp}) has the same degree as $\gamma$,
\item $\G$ has finite co-volume.
\end{enumerate}
\end{theorem}

\subsection{The possible values of $p$.}

If the generator $f$ has order $4$,  then Theorem \ref{2genthm}, implies first, that $\gamma$ is an algebraic integer,
secondly, that $\IQ(\gamma)$ has exactly one complex place and thirdly, 
that $\IQ(\gamma)$  must contain $L = \IQ(\cos 2 \pi/p)$. Let 
\begin{equation} [\IQ(\gamma) : L ] = r
\end{equation}
Since the minimal polynomial for $\gamma$ factors over $L$ we know $n=r\mu=[\IQ(\gamma):\IQ]$. Further, all real embeddings of $\gamma$ lie in the interval $[-2,0]$.

Next,  from the disjoint isometric circles criteria,  we obtain the following information if the group $\Gamma=\langle f,g\rangle$ is not freely generated --- a necessary condition if $\G$ is to be arithmetic.

\begin{lemma} \label{modlem} Let $\langle f,g\rangle$ be a discrete group with $o(f)=4$ and $o(g)=p$ and which is not free on these two generators. Let $\gamma=\gamma(f,g)$,  Then 
\begin{enumerate}
\item $\Re e(\gamma)>-4$ and this is sharp as for each $p\geq 2$, $GT(4,p;\infty)$ has $[f,g]$ parabolic and $\gamma=-4$.  This group is free on generators.  Further, $GT(4,p;n)$ has $\gamma=-2-2\cos(\frac{\pi}{n})$ which tends to $-4$ with $n$.  None of these groups are free on generators, \cite{MM3}.
\item For each $p\geq 2$ we have
\begin{equation}\label{trivial bound} |\gamma|\leq  4 \big(\sqrt{2} \cos(\frac{\pi }{p})+1\big)+2\cos(\frac{2 \pi }{p}).\end{equation}
This estimate is sharp and achieved by the $(4,p,\infty)$ triangle group.  The $(4,p,n)$ triangle groups have 
\[\gamma=2 \left(\cos(\frac{2 \pi }{n})+2\cos(\frac{\pi }{p}) \left(\sqrt{2} \cos(\frac{\pi }{n})+\cos(\frac{\pi }{p})\right)\right)\]
which tends to $ \big(\sqrt{2} \cos(\frac{\pi }{p})+1\big)+2\cos(\frac{2 \pi }{p})$ as $n\to\infty$.  The $(4,p,n)$-triangle groups are not free on generators.
\end{enumerate}
\end{lemma}
\noindent{\bf Proof.} To obtain the bound $|\gamma|\leq  4 \big(\sqrt{2} \cos(\frac{\pi }{p})+1\big)+2\cos(\frac{2 \pi }{p})$ for each $p$ we simply find the point of maximum modulus on the region $\Omega_p$ using calculus.  Next,  there are points within $\Omega_p$ whose real part is less than $-4$.  We have to show these give rise to groups freely generated.  We do this by a cut and paste argument to produce a fundamental domain from the isometric circles,  which are no longer disjoint.  We can normalise so that the isometric circles of $f$ and of $g$ are
\[ I(f): |z\pm i|=\sqrt{2},\quad I(g): |z\pm i\omega\cot \frac{\pi}{p}|=\frac{1}{\sin\frac{\pi}{p}} \]
and that $\Im m(i\omega)\geq 0$.  Let $D_1=\{z:|z+i\omega\cot \frac{\pi}{p}|=\frac{1}{\sin\frac{\pi}{p}}\}$ and $D_2=\{z:|z-i\omega\cot \frac{\pi}{p}|=\frac{1}{\sin\frac{\pi}{p}}\}$ be the disks bounded by the isometric circles of $g$.
\begin{lemma}\label{fog}  Suppose that $f^{-1}(D_1)\cap (D_1\cup D_2)=f(D_2)\cap (D_1\cup D_2)=\emptyset$  Then $\langle f,g \rangle$ is free on generators.
\end{lemma}
\noindent{\bf Proof.} We replace the fundamental domain $D(i,\sqrt{2})\cap D(-i,\sqrt{2})$ for $f$ with the domain
\[ [D(i,\sqrt{2})\cap D(-i,\sqrt{2})]\cup D_1\cup D_2)\setminus (f^{-1}(D_1) \cup f(D_2)) \]
By construction the exterior of the isometric circles of $g$ lie in this region and this region is a fundamental domain for $f$.  This the ``ping-pong'' hypotheses are satisfied.
\hfill $\Box$

\medskip
 \scalebox{0.5}{\includegraphics[viewport=-100 300 500 730]{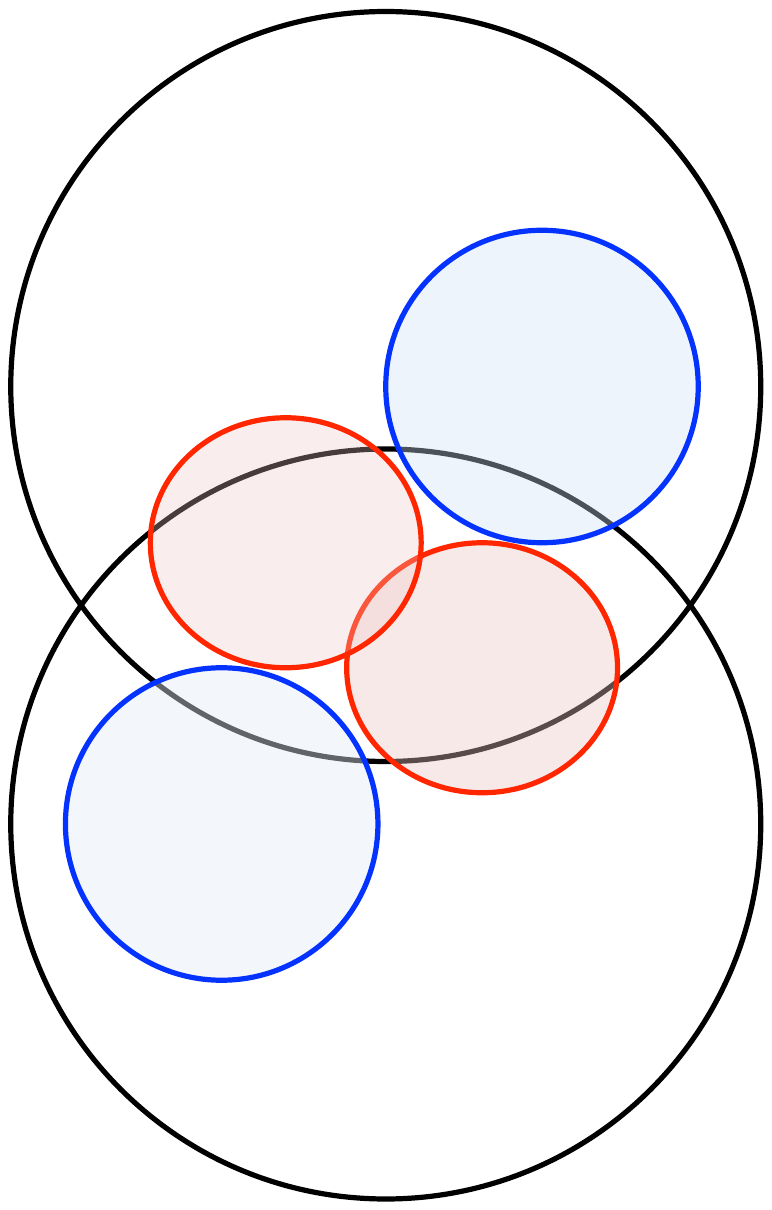}}
     
     \medskip
     \noindent{\bf Figure 3.}{\em  Isometric circles of $f$ and $g$ (in red).  The blues disks are the images of the leftmost isometric circle of $f$ under $f^{-1}$,  and the rightmost under $f$.The modified fundamental domain for $f$ consists of the intersection of the fundamental disks for $f$,  together with the fundamental disks for $g$ (red) and then deleting their images (two blue disks).}
     
     \medskip
     
This cut and paste technique for moving the fundamental domains around to a good configuration can be extended,  but quickly becomes to complex to obtain useful information. It apparent that is is only possible when the angle between the axes of $f$ and $g$ is small compared with the radius of the disk. We want to turn the hypotheses of Lemma \ref{fog} into a restriction on $\gamma$.  Things have been set up so 
\[ f(z) = \frac{z+i}{iz+1}. \]
There is also an evident rotational symmetry.  We therefore need to assert $f(D)$ is disjoint from the isometric circles of $g$ where $D$ is the left-most isometric circle of $g$. We write
\[\zeta=i\omega\cot \frac{\pi}{p}+\frac{e^{i\theta}}{\sin\frac{\pi}{p}}, \quad \mbox{that is $\zeta\in \partial D$} \]
and we need to show
\[ |f(\zeta)\pm i\omega\cot \frac{\pi}{p}|\geq \frac{1}{\sin\frac{\pi}{p}}  \]
These yield the two inequalities
\[ \left|  -i+i w \cot \frac{\pi }{p} +\frac{2}{-i+i w \cot \frac{\pi }{p} +e^{i \theta } \csc \frac{\pi }{p} }\right|\sin \frac{\pi }{p}  \geq 1\]
and 
\[\left|  -i -iw \cot \frac{\pi }{p} +\frac{2}{-i+i w \cot \frac{\pi }{p} +e^{i \theta } \csc \frac{\pi }{p} }\right| \sin  \frac{\pi }{p}  \geq 1\]

\subsection{First bounds.}

 We now make use of these facts,
together with  the inequalities that $\gamma$ and its conjugates must satisfy 
values for the triple $(4,p,r)$ for which there may exist a $\gamma$-parameter 
corresponding to an arithmetic Kleinian group which is not free using the criteria above. Our first estimate is straightforward and is simply used to cut down the set of triples we search through.

\begin{lemma}\label{easybound}  Let $\langle f,g\rangle$ be a discrete group with $o(f)=4$ and $o(g)=p$ and which is not free on these two generators. Let $\gamma=\gamma(f,g)$.  Then the degree $[\IQ(\gamma):\IQ]\leq 16$.
\end{lemma}
 \noindent{\bf Proof.} Let $P(z)$ be the minimial polynomial for $1+\gamma$ of degree $n$. Then $P$ has complex roots $1+\gamma$ and $1+\bar\gamma$ and,  as $\IQ(\gamma)$ has one complex place, $n-2$ real roots $r_i$ in the interval $[-1,1]$. As $P$ is irreducible and integral
 \begin{eqnarray*}
1 &\leq & |P(0)|^2|P(-1)||P(1)| = |1+\gamma|^4|\gamma|^2|\gamma+2|^2 \prod_{i=1}^{n-2} r_i^2(1-r_i^2) \\
&\leq &  \left(7+4 \sqrt{2}\right)^4\left(6+4 \sqrt{2}\right)^2\left(8+4 \sqrt{2}\right)^24^{2-n} \end{eqnarray*}
where we have used the trivial bounds on $|\gamma|$ from  Lemma \ref{modlem}. This is not possible if $n\geq 17$.  \hfill $\Box$

\medskip

There is far more information to be had here when we identify a specific $p$.  For instance if $p=4$  and $P$ is the minimal polynomial for $\gamma$,  following the above argument we see $r_i\in [-1,0]$ and 
 \begin{eqnarray*}
1 &\leq & |P(0)||P(-1)| = |\gamma|^2|1+\gamma|^2 \prod_{i=1}^{n-2} |r_i |(1-|r_i|)  
\leq   \left(10884 + 7696\sqrt{2}\right)^2 4^{2-n} \end{eqnarray*}
This gives us a total degree bound of at most $9$.   In fact this (where the orders of $f$ and $g$ are the same) is typically the worst case, and degree $9$ yields a feasible search space (as per \cite{FR}).

\medskip

The field $L$ is totally
real and the embeddings $\sigma : L \rightarrow \IR$ are defined by 
\[ \sigma( \cos \frac{2 \pi}{p})  =   \cos \frac{2 \pi j}{p}, \hskip15pt (j,p)=1,  \;\; j< [p/2].\] 
 Let us denote these embeddings by $\sigma_1, \sigma_2, \ldots , \sigma_{\mu}$, 
with $\sigma_1 = {\rm Id}$.  Here $\mu=[L:\IQ]$ and from the online list of integer sequences A055034 we find that
\[  [L:\IQ] = \Phi[p]/2, \]
where $\Phi$ is Euler's function. Since $\gamma$ is complex,  our bound from Lemma \ref{easybound} implies that $\Phi[p] \leq 16$.  Hence $p\leq 60$,  and in fact only $30$ possibilities for $p$ remain.

\subsection{Discriminant bounds.}

We now follow \cite[\S 4.2]{MM} to obtain an estimate on the relative discriminant. The discriminant of the minimal polynomial for $\gamma$  is
 \[ disc(\gamma) = |\gamma-\bar\gamma|^2 \prod_{ i=1}^{n-2} |\gamma - r_i|^4   \prod_{1 \leq i < j \leq n-2} (r_i - r_j)^2, \]
where the $n-2$ real roots $r_i\in [-2,0]$. Let
\begin{equation}
\mu = [L:\IQ]
\end{equation}
so that the total degree of the field $\IQ(\gamma)$ is $n=r\mu$.

For $n \geq 2$, let $D_n$ denote the minimum absolute value of the discriminant of 
any field of degree $n$ over $\IQ$ with exactly one complex place. For small
values of $n$ the number $D_n$  has been widely investigated (\cite{CDO,Di,DO})
and lower bounds for $D_n$ for all $n$ can be
computed (\cite{Mull,Od,Rodgers,Stark}). In \cite{Od}, the bound 
is given in the form $D_n > A^{n-2} B^2 \exp(-E)$ for varying values of 
$A,B$ and $E$. Choosing, by experimentation, suitable values from this table 
we obtain the bounds shown in Table 1.  We will use more precise data later.

\begin{table}[h]
\begin{center}
\begin{tabular}{clcl}
Degree $n$ & Bound & Degree $n$ & Bound  \\
2 & 3 &
3 & 27 \\
4 & 275 &
5 & 4511 \\
6 & 92779 &
7 & 2306599 \\
8 & 68856875* &
9 & $0.11063894 \times 10^{10} $ \\
10 & $0.31503776 \times 10^{11}$ &
11 & $0.90315026 \times 10^{12}$ \\
12 & $0.25891511 \times 10^{14}$ &
13 & $0.74225785 \times 10^{15}$ \\
14 & $0.21279048 \times 10^{17}$&
15 & $0.61002775 \times 10^{18}$ \\
16 & $0.17488275 \times 10^{20}$ &
17 & $0.50135388 \times 10^{21}$ \\
18 & $0.14372813 \times 10^{23} $ &
19 & $0.41203981 \times 10^{24}$ \\
20 & $0.11812357 \times 10^{26}$ 
\end{tabular}
\caption{Discriminant Bounds}
\end{center} 
\end{table}
\subsection{Schur's bound.}
We will need to use Schur's bound \cite{S} which gives that, 
if $-1 \leq x_1 < x_2 < \cdots < x_r \leq 1$ with $r \geq 3$ then
\begin{equation}\label{eqn31}
 \prod_{1 \leq i < j \leq r} (x_i - x_j)^2 \leq M_r = \frac{2^2\,3^3\, \ldots r^r\, 2^2 \, 3^3 \, \ldots (r-2)^{r-2}}{3^3\, 5^5 \, \dots (2r-3)^{2r-3}}.
\end{equation}
 
The bounds we have at hand give 
\begin{equation}\label{disc} D_n \leq disc(\gamma) \leq  |\gamma-\bar\gamma|^2  (2+|\gamma|)^{2(n-2)}   M_{n-2}, \quad\quad n \geq \Phi[p] \end{equation}
 where $\Phi$ is the Euler phi function.  However,  this can be improved to capture the fact that the real roots of the minimal polynomial for $\gamma$ are not so well distributed as those of the sharp estimate for Schur's formula (roots of Chebychev polynomials).  Following \cite{MM} we 
let $\Delta_1 = \delta_{\IQ(\gamma) \mid L}$, the relative discriminant 
of the field extension $\IQ(\gamma) \mid L$, and let $\Delta$ denote the 
discriminant of the basis $1, \gamma, \gamma^2, \ldots , \gamma^{r-1}$ over
$L$. Then 
$$| N_{L \mid \IQ}(\Delta) | \geq | N_{L \mid \IQ}(\Delta_1) |.$$
 
 \begin{equation}\label{eqn33}
|N_{L \mid \IQ}(\delta_{\IQ(\gamma) \mid L})| = |\Delta_{\IQ(\gamma)}|/\Delta_L^r.
\end{equation}
Next set
 \begin{equation}\label{18}
 K(p,r) =   M_{r-2}[2(\sqrt{2}+ \cos\frac{\pi}{p})^2]^{4(r-2)}\left(\sin\frac{\pi}{p} \right)^{2(r-2)(r-3)}|\gamma-\bar\gamma|^2
 \end{equation}
 
The discriminant $\Delta_p$ of the field  $\IQ(\cos \frac{2\pi}{p})$ is given in \cite[(30)]{MM}.    Then (see \cite[(31)]{MM}) we have the inequality
 \begin{equation}\label{19}
 K(p,r) \left(\sin\frac{\pi}{p}\right)^{-2r(r-1)}\left( \frac{\delta_{p}}{8}\right)^{\mu r(r-1)} M_r^{\mu-1} \geq \max\{1, D_n/\Delta_L^r\}
 \end{equation}
 where here
 \[ \delta_p = \left\{ \begin{array}{ll}
                1 & {\rm if~}p \neq \pi^{\alpha}, \;\; \pi~{\rm a~prime} \\
               \pi & {\rm if~}p = \pi^{\alpha}, \;\; \pi ~{\rm a~prime},
                 \end{array}
        \right.  \]

If (\ref{disc}) fails to eliminate a case,  then we put together inequalities (\ref{18}) and (\ref{19}) to try to eliminate it.  We obtain the following remaining possibilities.  
\begin{center}
\begin{tabular}{|c|c|c|l|}
\hline
$p$ & $[\IQ(\gamma):L]$ & $\Delta_p$ & total degree $n$ \\ \hline 
3 & 1 & $1$& $n\leq 11$ \\ \hline 
4 & 1 & $1$&$n\leq 8$\\ \hline
5 & 2 & $5$&$n=4,6,8$\\ \hline
6 & 1& $1$&$n=2,3,4$ \\ \hline
7 & 3 & $49$&$n=6,9$ \\ \hline
8 & 2 & $8$&$n=4,6,8$  \\ \hline
9 & 3 & $81$&$n=6,9$\\ \hline
10 & 2 & $5$&$n=4,6$\\ \hline
11 & 5 & $14641$&$n=10$  \\ \hline
12 & 2& $12$&$n=4,6,8$ \\ \hline
15 & 4 & $1125$&$n=8$ \\ \hline
18 & 3 & $81$&$n=6,9$ \\ \hline
20 & 4 & $2000$&$n=8$ \\ \hline
24 & 4 & $2304$&$n=8$ \\ \hline
30 & 4  & $1125$& $n=8$ \\ \hline
\end{tabular} 
\end{center}
 
For all but small $p$ we find that the minimal polynomial for $\gamma$ is quadratic or cubic over the base field $\IQ(\cos 2\pi/p)$.  We will work through some explicit examples to find all these algebraic integers in a moment,  but all these search spaces are feasible with a little work except $p=3$ for which we will find alternative arguments, and  $p=4$ which will follow from \cite{FR} once we have better bounds on the shape of the moduli space.  The case $p=5$ also looks troublesome.  Typically we search for $1+\gamma$ as then all real roots are in $[-1,1]$ and it is easier to get bounds on the symmetric functions of the roots, and therefore the coefficients of the polynomials.

\subsection{$p=7$, $p=9$, $p=11$ \& $p=18$.}  There are four cases not yet eliminated where the total degree exceeds $8$, these are $p=11$ and total degree $n=10$ and $p=7, 9,18$ with total degree $6$ or $9$.  Let us work through these case, first with $p=11$. The minimal polynomial $P$ is quadratic over $\IQ(\cos(2\pi/11))$. It has the complex conjugate pair of roots $\gamma$ and $\bar\gamma$ and $8$ real roots which come in pairs,  two in $[-2\sin(2\pi/11),0]$,  two in $[-2\sin(3\pi/11),0]$, two in $[-2\sin(4\pi/11),0]$ and two in $[-2\sin(5\pi/11),0]$.  Therefore the norm of the relative discriminant is bounded by 
\begin{eqnarray*}
\lefteqn{|\gamma-\bar \gamma|^2 4^4 \sin^2(2\pi/11)  \sin^2(3\pi/11)  \sin^2(4\pi/11) \sin^2(5 \pi/11) \Delta_{11}^2}\quad\quad\\
& = & |\gamma-\bar \gamma|^2 \; \frac{11}{2\sin^2(\pi/11)} \; (14641)^2    \geq D_{10}= 0.31503776 \times 10^{11}.
\end{eqnarray*}
This implies 
\[ |\gamma-\bar \gamma|=2\Im m(\gamma) \geq 25.6581, \]
and this of course is a contradiction.

\medskip
Next with $p=18$.  The minimal polynomial $P$ is cubic over $\IQ(\cos(2\pi/18))$. It has the complex conjugate pair of roots $\gamma$ and $\bar\gamma$ and a real root in $[-2\sin^2\pi/18,0]$ and  $6$ real roots which come in triples,  three in $[-2\sin(5\pi/18),0]$,  and three in $[-2\sin(7\pi/18),0]$.  Therefore the discriminant is bounded by 
\begin{eqnarray*}
\lefteqn{|\gamma-\bar \gamma|^2 (|\gamma|+ 2\sin^2\pi/18)^4 16^{-2} (2\sin(5\pi/18))^6(2\sin(7\pi/18))^6 \Delta_{18}^3}\quad\quad\\
& = & |\gamma-\bar \gamma|^2 (|\gamma|+ 2\sin^2\pi/18)^4 \times164609 \geq D_{9}=  0.11063894 \times 10^{10} .
\end{eqnarray*}
This implies 
\[ |\gamma-\bar \gamma|^2 (|\gamma|+ 2\sin^2\pi/18)^4  \geq 6721.32, \]
which gives $|\gamma|\geq 14$ and this is a contradiction.

We next consider the degree $6$ case to show in a simple case how we can use more refined information about the discriminant to eliminate a value.  Following the above argument we would achieve the estimate
\begin{equation}\label{deg6}
 {\cal N}(\IQ(\gamma)|L) \times 81^2  \geq D_{6}= 92779.
\end{equation}
The norm ${\cal N}$ is a rational integer but $92779$ is prime.  This shows the inequality can never be sharp.  It gives the estimate $|\gamma-\bar \gamma|\geq 1.8$.
However,  from the useful resource $https://hobbes.la.asu.edu/NFDB/$, \cite{JR} we can identify all the polynomials giving number fields with one complex place, degree $6$ and discriminant less than $3\times 10^6$ - there are $2352$ such fields.   This list is proven complete.  The one we seek for $\IQ(\gamma)$ has $\IQ(\cos \frac{2\pi}{18})$ as a subfield and the formula $(\ref{deg6})$ shows the discriminant has $3^8$ as a factor and Galois group $A_4\times C_2$.   This leaves $19$ candidates. These are identified in the table below with $\delta$ as the distance between the smallest and largest roots.  The polynomials presented are unique up to integer translation.  

\medskip

\begin{tabular}{|c|c|c|l|c|}
\hline
 & $-\Delta(\IQ(\gamma))$ &  $p(z)$ & $\delta$ \\
\hline
1 &$ 3^8 \times 19 $ & $x^6 - 3x^4 - 2x^3 + 3x^2 + 3x - 1.$ & $2.76$\\
\hline
2 &$ 3^9 \times 17$ & $   x^6 - 3x^4 - 5x^3 + 3x + 1.$ & $2.91$ \\
\hline
3 &$ 3^8 \times 107 $ & $ x^6 - 3x^5 + 3x^4 - 9x^2 + 6x - 1.$ & $3.53$ \\
\hline
4 &$ 3^8 \times 2^6   $ & $ x^6 - 3x^2 + 1.$ & $2.47$ \\
\hline
5 &$ 3^9 \times 2^6$   & $x^6 - 3x^4 + 3.$& $ 3.18$\\
\hline
6 &$ 3^9 \times 37  $ & $  x^6 - 3x^5 + 3x^4 - x^3 - 3x^2 + 3x + 1.$ & $2.66$\\
\hline
7 &$ 3^9\times 53 $ & $ x^6 - 3x^4 - 4x^3 + 6x + 1.$ & $3.47$ \\
\hline
8 &$ 38 \times 1271  $ & $ x^6 - 3x^4 - 2x^3 - 6x^2 - 6x - 1.$ & $4.00$ \\
\hline
9 &$  3^8\times 163  $ & $ x^6 - 3x^5 + 3x^4 + 6x^3 - 15x^2 + 6x + 1.$ & $3.10$ \\
\hline
10 &$  3^8 \times 179 $ & $ x^6 - 3x^5 + 5x^3 - 3x^2 + 3.$ & $3.33$ \\
\hline
11 &$  3^9\times 73  $ & $x^6 - 8x^3 + 9x + 1$ & $2.59$\\
\hline
12 &$  3^8 \times 251 $ & $ x^6 - 3x^5 + 6x^2 + 6x - 1.$ &$ 3.24$ \\
\hline
13 &$  3^ 9\times 89 $ & $ x^6 - 3x^4 - 4x^3 - 9x^2 - 3x + 1.$ & $4.22$\\
\hline
14 &$  3^ 8\times 271 $ & $ x^6 - 3x^5 - 3x^3 + 6x^2 + 6x + 1.$ & $3.54$ \\
\hline
15 &$  3^ 8\times 17\times 19 $ & $x^6 - 3x^5 + 5x^3 - 12x^2 + 9x + 3 .$ & $4.47$ \\
\hline
16 &$   3^ 8\times 17\times 19 $ & $x^6 - 7x^3 - 6x^2 + 18x - 3.$ & $3.32$\\
\hline
 17 &$  3^9\times 109 $ & $x^6 - 3x^5 - 3x^4 + 3x^3 + 9x^2 + 9x + 3.$ &$4.06$\\
\hline
 18 &$  3^ 8\times 359 $ & $ x^6 - 6x^4 - 4x^3 - 3x^2 + 3.$ & $4.94$ \\
\hline
 19 &$  3^ 8\times 431 $ & $x^6 - 3x^5 + x^3 + 18x - 9.$ & $3.93$\\
\hline
\end{tabular}

\medskip

None of these polynomials satisfy the root restrictions we require.  In particular $\delta > 2$ implies no integral translate of any of these polynomials has all its real roots in an interval of length $2$. Therefore no integral translate can have the root distributions we have established as necessary.  In fact in this case $\delta > 2 \sin^2(7\pi/18)\approx 1.766$ will suffice.   We therefore have 
\begin{equation} 
 |\gamma-\bar\gamma|^2\big(2\sin^2 \frac{5\pi}{18}\big)^2 \big(2\sin^2 \frac{7\pi}{18}\big)^2  81^2 \geq  {\cal N}(\IQ(\gamma)|L) \times \Delta_p^2  \geq 3\times 10^6. 
\end{equation}
Hence
\[ \Im m(\gamma) \geq   5.15829  \]
and this easily implies the group is free on its generators.  

\medskip

The remaining cases in this subsection are $p=7,9$. As with $p=18$ there are two cases when the extension over $L$ is quadratic or a cubic.  In the cubic case Galois group in question is the wreath product $S(3)wr3$ (of order $2^3 3^4$ - T28 in the notation of \cite{JR})  There are $285$ fields to consider with discriminant less than $3.84\times 10^{10}$. In the case $p=7$, we must have $7^6$ as a factor of the discriminant (the smallest discriminant here is $7^6 \times 22679$ ) and in the case $p=9$ we must have $3^{12}$ as a factor (the smallest discriminant here is $3^{12}\times53\times 163$).
\begin{eqnarray*}
|\gamma-\bar \gamma|^2 (|\gamma|+ 2\sin^2\pi/9)^4 16^{-2} (2\sin(2\pi/9))^6(2\sin(4\pi/9))^6 \Delta_{9}^3\geq 3.84\times 10^{10}\\
|\gamma-\bar \gamma|^2 (|\gamma|+ 2\sin^2\pi/7)^4 16^{-2} (2\sin(2\pi/7))^6(2\sin(4\pi/7))^6 \Delta_{7}^3 \geq   3.84\times 10^{10} . 
\end{eqnarray*}
The first case $p=7$ quickly gives either 
\begin{itemize}
\item If $|\gamma-\bar\gamma|\leq 8$,  then $|\gamma|\geq 9.16471$, or 
\item if $|\gamma-\bar\gamma|\leq 5$,  then  $|\gamma|\geq 11.6923$.
\end{itemize} 
and for $p=9$,  simply $|\gamma|\geq 11.1919$.  Together  these inequalities  imply that the associated group is free.

In the quadratic cases we have have a field of degree $6$.  We enumerate the degree $6$ fields (actually we are enumerating the polynomials) with one complex place and discriminant less than $1.7 \times 10^6$.  In $p=7$ there is a factor of $7^4$ and when $p=9$ a factor of $3^8$.  The Galois group is $A_4\times C_2$ and there are $40$ such fields (polynomials).  Then
\begin{eqnarray*}
|\gamma-\bar \gamma|^2  (2\sin(2\pi/9))^2(2\sin(4\pi/9))^2 \Delta_{9}^2\geq  2.2\times 10^6.\\
|\gamma-\bar \gamma|^2  (2\sin(2\pi/7))^2(2\sin(4\pi/7))^2 \Delta_{7}^2 \geq   2.2\times 10^6 . 
\end{eqnarray*}
Together these inequalities yield $|\gamma-\bar \gamma|>10 $ and again this implies the groups is discrete and free on generators.

\medskip

Apart from the case $p=3,4$ and $5$ we have now established the total degree is at most  $8$.  We will actually run through searches to verify these results independently and we will discuss that later.

\subsection{$p=15,20,24,30$,  degree $8$.}

These four cases arise from a quadratic polynomial over the base field $L$.  In each case we can compute the relative discriminant as
\begin{equation}
|\gamma-\gamma|^2 \; \delta_2^2\; \delta_3^2\; \delta_4^2 \; \Delta_p^2 \geq D_8 = 68856875
\end{equation}
Here $\delta_i=\sigma_i(2\sin^2\frac{\pi}{p})$ as $i$ runs through the three non-identity Galois automorphisms $\sigma_i$ of $L$,  since a pair of roots lies in the interval $[-2\sigma_i( \sin^2\frac{\pi}{p}),0]$.  For $p=15$ we find 
\begin{eqnarray*}
p=15,& & |\gamma-\gamma|^2 \; 2^6 \;  \sin^2 \frac{2\pi}{15}  \sin^2 \frac{4\pi}{15}  \sin^2 \frac{7\pi}{15}  \times  1125^2 \geq  68856875.\\
\end{eqnarray*}
Hence $\Im m(\gamma) \;  \geq  5.10151.$ and this implies the group is free.  For the other cases the discriminant bound is not enough to eliminate them.  Thus we seek a stronger bound.  
 From the resource $https://hobbes.la.asu.edu/NFDB/$, \cite{JR} once again we can identify all the polynomials giving number fields with one complex place, degree $8$ and discriminant less than $2\times 10^9$.   When $p=20$ $2^8\times 5^6$ is a factor of the discriminant (smallest such is $2^8\times 5^6\times 79$), $p=24$ finds $2^16 \times 3^4$ as a factor (smallest such is $2^{16}\times 3^4 \times 47$) and $p=30$ finds $3^4 \times 5^6$ as a factor (smallest such is $3^4 \times 5^6 \times 59$). The Galois groups are either $[2^4]^4$ ($T27$ in \cite{JR}) when $p=20,30$ or $[2^4]E(4)$ ($T31$ in \cite{JR}) when $p=24$.
These lists are proven complete.  When $p=20$ we have the following $12$ possibilities.

\medskip

\begin{tabular}{|c|c|l|c|}
\hline
 &  $-\Delta$ & $p(z)$ &$\delta$ \\
\hline
1 &$2^85^679$ & $ x^8 - 2x^7 + x^6 - 6x^5 - x^4 + 12x^3 - x^2 - 4x + 1$ & $3.29$\\
\hline
2 &$2^85^719$ & $  x^8 - 4x^7 + 2x^6 + 8x^5 - 10x^4 + 2x^3 + 7x^2 - 6x + 1.$ &$3.55$\\
\hline
3 &$ 2^85^6199 $ & $x^8 - 4x^7 + 3x^6 + 6x^4 - 8x^2 + 2x + 1.$ & $3.30$\\
\hline
4 &$ 2^85^6239 $& $x^8 - 4x^7 + 10x^6 - 16x^5 - x^4 + 24x^3 - 10x^2 - 4x + 1.$ & $2.97$ \\
\hline
5 &$ 2^85^759$   & $x^8 - 2x^7 - 7x^6 + 16x^5 + 5x^4 - 24x^3 + 8x^2 + 8x - 4.$ & $4.75$ \\
\hline
6 &$2^{12}5^619$ & $ x^8 - 2x^7 - 8x^6 + 10x^5 + 11x^4 - 10x^3 - 2x^2 + 4x + 1.$ &$5.53$\\
\hline
7 &$ 2^{12}5^619  $ & $ x^8 - 8x^6 - 4x^5 + 8x^4 - 8x^3 - 12x^2 + 4x + 1.$ & $4.82$\\
 \hline
8 &$ 2^85^6359 $ & $x^8 - 9x^6 - 10x^5 + 11x^4 + 30x^3 + 11x^2 - 10x + 1.$ & $5.16$\\
 \hline
9 &$ 2^83^45^7$ & $x^8 - 10x^4 + 15x^2 - 5.$ & $2.64$ \\
 \hline
10 &$ 2^83^45^7$ & $x^8 - 5x^6 + 5x^4 + 5x^2 - 5.$ & $3.36$ \\
 \hline
11&$2^85^6439$ & $x^8 - 2x^7 - 3x^6 + x^4 + 10x^3 + 3x^2 - 6x + 1.$ & $4.08$ \\
 \hline
 12&$2^85^6479$ & $x^8 - 3x^6 - 11x^4 - 30x^3 - 17x^2 + 1.$ & $4.13$\\
 \hline
\end{tabular} 

When  $p=30$ there are $26$ possibilities (which we do not enumerate)

Finally  when $p=24$,  there are $7$ possibilities.

\begin{tabular}{|c|c|l|c|}
\hline
 &  $-\Delta$ & $p(z)$ & $\delta$ \\
\hline
1 &$2^{16}3^4 47 $ & $ x^8 - 8x^5 + x^4 + 12x^3 - 2x^2 - 4x + 1$ & $2.29$\\
\hline
2 &$2^{18}3^4 23 $ & $   x^8 - 6x^6 - 8x^5 + 4x^4 + 12x^3 - 2x^2 - 4x + 1.$ & $4.29$ \\
\hline
3 &$ 2^{16}3^4191 $ & $x^8 - 4x^6 - 11x^4 - 12x^3 + 30x^2 + 36x + 9.$ & $4.42$\\
\hline
4 &$  2^{16}3^6 23 $& $ x^8 - 4x^7 + 2x^6 + 8x^5 - 14x^4 + 4x^3 + 14x^2 - 8x + 1.$ & $4.06$\\
\hline
5 &$ 2^{16}3^4 239$   & $x^8 - 6x^6 + 5x^4 - 12x^3 - 12x^2 + 1.$ & $4.49$ \\
\hline
6 &$ 2^{18}3^4 71 $ & $ x^8 - 4x^7 + 10x^6 - 12x^5 - 20x^4 + 40x^3 - 18x^2 + 1.$ & $3.84$\\
\hline
7 &$ 2^{20}3^4 23  $ & $ x^8 - 8x^6 - 4x^5 + 8x^4 - 8x^3 - 12x^2 + 4x + 1.$ & $4.82$\\
 \hline
\end{tabular}

 Hence we may use the bound $2\times 10^9$ to see
\begin{eqnarray*}
p=20,& & |\gamma-\gamma|^2 \; 2^6 \;  \sin^2 \frac{3\pi}{20}  \sin^2 \frac{7\pi}{20}  \sin^2 \frac{9\pi}{20}  \times  2000^2 \geq  2\times 10^9.\\
p=24,& & |\gamma-\gamma|^2 \; 2^6  \; \sin^2 \frac{5\pi}{24}  \sin^2 \frac{7\pi}{24}  \sin^2 \frac{11\pi}{24}  \times  2304^2 \geq  2\times 10^9.\\
p=30,& & |\gamma-\gamma|^2 \; 2^6  \; \sin^2 \frac{7\pi}{30}  \sin^2 \frac{11\pi}{30}  \sin^2 \frac{13\pi}{30}  \times  1125^2 \geq  2\times 10^9.
\end{eqnarray*} 
Then 
\begin{eqnarray*} 
p=20,& &  \Im m(\gamma)  \;  \geq  8.75528\\
p=24,& &  \Im m(\gamma)  \;   \geq  5.29112.\\
p=30,& &  \Im m(\gamma)  \;   \geq  6.94947.
\end{eqnarray*}
Thus all these groups are free as well.

 \begin{center}
\scalebox{0.7}{\includegraphics[viewport=80 340 800 800]{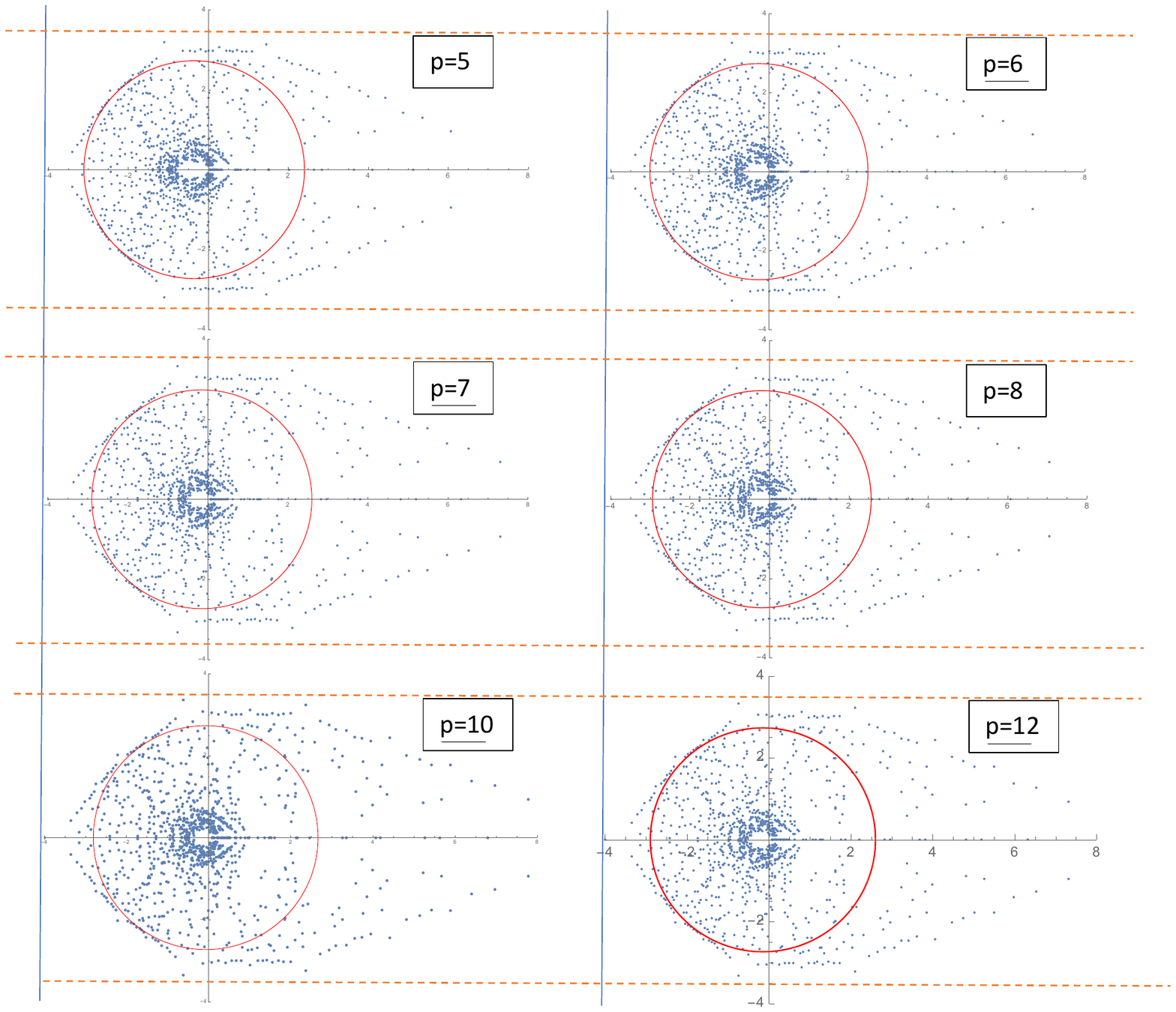}}
 \end{center}
         
\noindent{\bf Figure 2.} {\em Rough descriptions of the moduli space of  $\;\IZ_4*\IZ_p$ for $p=5,6,7,8,10$ and $12$. These moduli spaces contain the convex hull of the points (which are roots of the Farey polynomials of degree up to $120$, and contain about $513$ boundary cusp groups.)}

\medskip

From these rough descriptions of the moduli spaces one directly obtains concrete information such as the following.  We will give a careful (and sharper) proof for this later in the $(4,5)$ case and then again when we examine the $(3,4)$ and $(4,4)$ case in considerable detail.  The details for other $p$ are exactly the same except that since we only give rough bounds it is only necessary to examine pleating rays of small slope.

\begin{lemma} Let $\Gamma=\langle f,g \rangle$ be an arithmetic lattice generated by $f$ of order $4$ and $g$ of finite order $p\in \{5,6,7,8,10,12\}$. Then
\begin{equation}
|\Im m(\gamma(f,g))|\leq 4.
\end{equation}
This bound is also true if $\Gamma$ is not a lattice and is not free.
\end{lemma}
The sharp bound here is nearly 3.5. 

\medskip

Next we consider the following remaining possibilities.  
\begin{center}
\begin{tabular}{|c|c|c|l|}
\hline
$p$ & $[\IQ(\gamma):L]$ & $\Delta_p$ & total degree $n$ \\ \hline 
3 & 1 & $1$& $n\leq 11$ \\ \hline 
4 & 1 & $1$&$n\leq 8$\\ \hline
5 & 2 & $5$&$n=4,6,8$\\ \hline
6 & 1& $1$&$n=2,3,4$ \\ \hline 
8 & 2 & $8$&$n=4,6,8$  \\ \hline 
10 & 2 & $5$&$n=4,6$\\ \hline 
12 & 2& $12$&$n=4,6,8$ \\ \hline 
\end{tabular} 
\end{center}

\subsection{$p=5,8,10$ and $p=12$.}

These are quadratic, cubic and quartic extensions over the base field which is a degree two extension of $\IQ$.  We next examine the limits of the approach we have used so far.

\subsubsection{$p=12$ : quartic.} 

As always,  we estimate the discriminant,  using the relative discriminant.  There are $4$ real roots in $[-2\sin^2\frac{5\pi}{12},0]$ and two in $[-2\sin^2\frac{\pi}{12},0]$.  We use Schur's bound to deal with the real roots in $[-2\sin^2\frac{5\pi}{12},0]$
\begin{equation}
|\gamma-\gamma|^2 \big(|\gamma|+2\sin^2\frac{\pi}{12}\big)^8 \; \big(\sin^2\frac{5\pi}{12}\big)^{12} \; \frac{2^23^34^42^2}{3^3 5^5}\; 12^4 \geq D_8 \geq 68856875
\end{equation}
We quickly see that in order to get useful information on $\gamma$ we will need to bound the discriminant by something of the order $10^{12}$.  The Galois group in question is $[2^4]E(4)$ of order $64$.  However there are (provably) $15,648$ polynomials which meet our criteria.   
\subsubsection{$p=12$ : cubic.} We have 
\begin{equation}
|\gamma-\gamma|^2 \big(|\gamma|+2\sin^2\frac{\pi}{12}\big)^4 \; 16^{-2} \big(2\sin^2\frac{5\pi}{12}\big)^{6} \;  12^3 \geq D_6 
\end{equation}
To get useful information we need a discriminant bound of about $10^8$ with a factor of $2^63^3$.  The Galois group is $A_4C_2$ ($T6$ in \cite{JR}).  There are (provably) $2319$ such polynomials.  
\subsubsection{$p=12$ : quadratic.} We have 
\begin{equation}
|\gamma-\gamma|^2   \big(2\sin^2\frac{5\pi}{12}\big)^{2} \;  12^2 \geq D_4 
\end{equation}
To get useful information we need a discriminant bound of about $4\times 10^5$ with a factor of $2^43^2$.  The Galois group is $D4$ ($T3$ in \cite{JR}).  There are (provably) $3916$ such polynomials (the smallest discriminant is $2^43^223$ with polynomial $x^4 - 2x^3 - x^2 + 2^x - 2$).  With a bit of work this method will work here but for quadratics it is far easier to search. 

\medskip

The point to discussing these three cases is that the methods applied above produce too many polynomials to be really useful.  It is oftentimes simpler to search for the polynomials directly.  The several thousand polynomials identified above will not have their real roots in the intervals we require,  nor their complex roots bounds bounded appropriately.   Finally the factorisation condition discussed in \S \ref{factorisation}, and which we will come to rely on heavily to eliminate cases, doesn't seem to be directly applicable in any way.

\section{Searches.}  In this section we give a few pertinent examples of the searches we performed to eliminate cases.  First off is a case we have already dealt with,  but which admits some useful features which are easier to follow.

\subsection{The case $p=4$ and $q=7$, degree $2$ over $\IQ(2\cos\frac{\pi}{7})$, total degree $6$. }
 
The intermediate field is
\[ L=\IQ(2\cos\frac{\pi}{7}) \]
of discriminant $49$.  We suppose $\gamma$ is quadratic over $L$.  The minimal polynomial for $\gamma$ is irreducible over $\IZ$ and then factors as 
\[ (z^2+a_1z+a_0)(z^2+b_1z+b_0)(z^2+c_1z+c_0) \]
where we may suppose that $\gamma$ and $\bar \gamma$ are roots of the first factor.  In order for the arithmeticity conditions to be satisfied we must have
\begin{itemize}
\item $z^2+b_1z+b_0$ has two real roots, say $r_\pm= \frac{1}{2}\big(-b_1\pm \sqrt{b_1^2-4b_0}\big)$, in $-[2\sin^2\frac{2\pi}{7},0]$,  thus $b_1> 0$
\begin{equation}\label{T1}
-4\sin^2\frac{2\pi}{7} \leq 2r_- = -b_1-\sqrt{b_1^2-4b_0},  
\end{equation}
and note the implication $b_1^2\geq 4b_0>0$,  strict inequality because of irreducibility.  We will always take square roots with positive real part. Also
\item $z^2+c_1z+c_0$ has two real roots (say $s_\pm$) in $-[2\sin^2\frac{3\pi}{7},0]$,  thus $c_1>0$ and 
\begin{equation}\label{T2}
-4\sin^2\frac{3\pi}{7} \leq 2s_-=-c_1-\sqrt{c_1^2-4c_0},  
\end{equation}
with $c_1^2\geq 4c_0>0$.
\end{itemize}
We also know from our  criterion regarding free groups used above that 
\begin{eqnarray}\label{T3}
0 < |\gamma|^2 & = &  a_0 <\left[ 2 \left(2+2 \sqrt{2} \cos  \frac{\pi }{7} +\sin \frac{3 \pi }{14} \right)\right]^2 \approx 106.991 \end{eqnarray}
It is efficient to improve this bound on $a_0$ and we can do this using the relative discriminant.  The total degree of $\gamma$ is $6$ and $|\Delta(\IQ(2\cos\frac{\pi}{7}))|=49$. As above there are $13$ polynomials yielding fields with one complex place an discriminant divisible by $7^4$ and less than $10^6$.  None of these polynomials have the properties we seek.  Hence
\[ |N(\gamma-\bar\gamma)| \times 49^2 \geq 10^6 \]
where $N$ is the relative norm, 
\[ |N(\gamma-\bar\gamma)| =| \gamma-\gamma|^2(r_+-r_-)^2(s_+-s_-)^2 \leq |\gamma-\gamma|^2 16 \sin^4\frac{2\pi}{7} \sin^4\frac{3\pi}{7}  \]
and hence
\begin{equation} \label{T4}
\Im m[\gamma]> 4.39079,  
\end{equation} 
and although we know this is enough to remove this case,  we keep calculating as this bound further implies  
\begin{equation} \label{T5}
19.279 \leq   |\gamma|^2=a_0 < 106.7
\end{equation} 
An improvement on  (\ref{T3}) alone.

\medskip
An integral basis for $\IQ(2\cos\frac{\pi}{7})$ can be found from $1$, $2\cos\frac{\pi}{7}$ and $2\cos\frac{3\pi}{7}$.  We also us the remaining Galois conjugate $2\cos\frac{5\pi}{7}$.  We can therefore write,  for rational integers $p_i,q_i,r_i$,  $i=1,2$,
\begin{eqnarray*}
a_0= p_0+ 2q_0\cos\frac{\pi}{7} +2r_0 \cos\frac{3\pi}{7},&& a_1= p_1+ 2q_1\cos\frac{\pi}{7} +2r_1 \cos\frac{3\pi}{7},\\
b_0=p_0+ 2q_0\cos\frac{3\pi}{7} +2r_0 \cos\frac{5\pi}{7}, && b_1=p_1+ 2q_1\cos\frac{3\pi}{7} +2r_1 \cos\frac{5\pi}{7} ,\\
c_0=p_0+ 2q_0\cos\frac{5\pi}{7} +2r_0 \cos\frac{\pi}{7}, && c_1=p_1+ 2q_1\cos\frac{5\pi}{7} +2r_1 \cos\frac{\pi}{7}. 
\end{eqnarray*}
From (\ref{T1}) and (\ref{T2}) we also have the bounds
\begin{eqnarray*}
0 < b_0< \left(2\sin^2\frac{2\pi }{7}\right)^2, && 0< b_1<4  \sin^2\frac{2\pi }{7}  ,\\
0 < c_0 < \left(2\sin^2\frac{3\pi }{7}\right)^2,  &&0< c_1<4  \sin^2\frac{3\pi }{7}.
\end{eqnarray*}
We can solve $p_0,q_0,r_0$ in terms of $a_0,b_0,c_0$.

\begin{eqnarray*}
p_0&=& \frac{1}{14} (4 (b_0+c_0+(b_0-c_0) \cos\big(\frac{\pi }{7}\big))+c_0 \text{csc}\big(\frac{\pi }{14}\big)-a_0 (-4+\text{csc}\big(\frac{3 \pi }{14}\big)+4 \sin\big(\frac{\pi }{14}\big))),\\
q_0 & = &  \frac{2}{7} (a_0 \cos\big(\frac{\pi }{7}\big)-b_0 \cos\big(\frac{\pi }{7}\big)+b_0 \sin\big(\frac{\pi }{14}\big)-c_0 \sin\big(\frac{\pi }{14}\big)+(a_0-c_0) \sin\big(\frac{3 \pi }{14}\big)),\\
r_0 &=&\frac{2}{7} (-b_0 \cos\big(\frac{\pi }{7}\big)+c_0 \cos\big(\frac{\pi }{7}\big)+a_0 \sin\big(\frac{\pi }{14}\big)-c_0 \sin\big(\frac{\pi }{14}\big)+(a_0-b_0) \sin\big(\frac{3 \pi }{14}\big))
\end{eqnarray*}
From this we deduce that
\begin{eqnarray*}
0\leq p_0 \leq 7,& 0\leq q_0\leq 21, &0\leq r_0\leq 12
\end{eqnarray*}
In exactly the same way we find the bounds
\begin{eqnarray*}
-1\leq p_1 \leq 3,& -7\leq q_1\leq 3, &-4 \leq r_1\leq 2.
\end{eqnarray*}
This gives us now $880880$ cases to consider.  Of these there is only one case which passes all the above tests.  That is $p_0=3$, $q_0= 11$, $r_0= 6$, $p_1= 1$, $q_1= -2$ and $r_1= -1$. The first pair of real roots are $-0.89457$ and $-0.462326$ both greater than $-2\sin^2\frac{2\pi}{7} = -1.22252$,  while the second pair  is $-1.63397$ and  $-0.0580492$,  both greater than $-2\sin^2\frac{3\pi}{7}  =-1.90097$.  However,  the complex roots are $\gamma = 1.52446 \pm 4.81327 i$ are easily seen to be in the space of groups freely generated as the imaginary part of $\gamma$ is too large.
\subsection{$p=12$}
 The intermediate field is
\[ L=\IQ(2\cos\frac{2\pi}{4},2\cos \frac{2\pi}{12})=\IQ(\sqrt{3}) \]
 \begin{center}
        \scalebox{0.8}{\includegraphics[viewport=80 520 500 800]{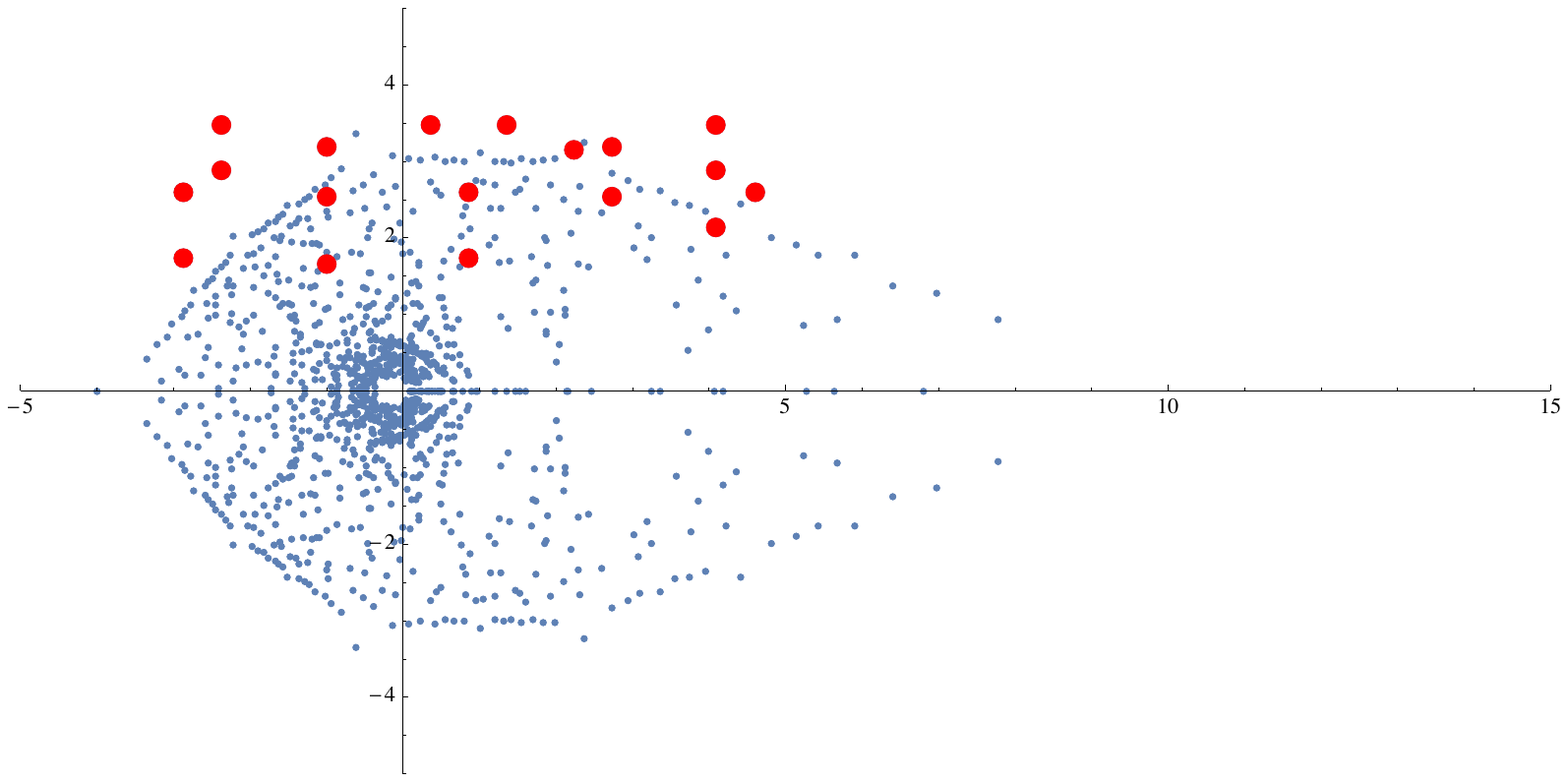}}
         \end{center}
\noindent{\bf Figure 3.}
{\em The moduli space of discrete groups generated by elliptics of order $p=4$ and $q=12$ obtained from identifying the roots of Farey polynomials.  The red points are those $\gamma$ satisfying the criteria:
\begin{enumerate}
\item the associated group is not obviously free,  that is 
\begin{itemize}
\item $0< Im(\gamma(f,g))< 4$,  
\item $\Re e(\gamma)\geq -4$,  (achieved in the $GT(4,12,\infty)$-generalised triangle group).
\item $\Re e(\gamma)\leq 3(2+\sqrt{3})$,  (achieved in the $(4,12,\infty)$-triangle group).
\end{itemize}
\item $\gamma$ is the complex root of a quadratic over $L$.
\item the two real embeddings of $\gamma$ lie in the interval $[-2\sin^2\frac{5\pi}{12},0]$.
\end{enumerate}}

\subsubsection{Degree 2 over $L=\IQ(\sqrt{3})$}
The minimal polynomial for $\gamma$ has the form
\[ p(z)=(z^2+a_1 z+ a_0)(z^2+b_1 z + b_0) \]
We always assume the first polynomial has a complex conjugate pair of roots.  Here $a_0,a_1,b_0,b_1\in \IQ(\sqrt{3})$ and $b_i$ is the Galois conjugate of $a_i$.  We write
\begin{eqnarray*}
a_0=p_0+q_0\sqrt{3}, &&b_0=p_0-q_0\sqrt{3}\\
a_1=p_1+q_1\sqrt{3}, &&b_1=p_1-q_1\sqrt{3}
\end{eqnarray*}
The first polynomial factor has roots $\gamma,\bar\gamma$ and the second factor has roots $r,s \in [-2\sin^2\frac{5\pi}{12},0]$.  Thus
\begin{eqnarray*}
a_0=|\gamma|^2, &&b_0=rs\\
a_1= - 2 \Re e[\gamma], &&b_1=-(r+s)
\end{eqnarray*}
From this we see that 
\begin{eqnarray*}
0< a_0 < 9  (7+4 \sqrt{3} ), && 0< b_0 < 4 \sin^4\frac{5\pi}{12} =\frac{7}{4}+\sqrt{3} \\
-3(2+\sqrt{3})  < a_1< 8, && 0< b_1<2+\sqrt{3}.
\end{eqnarray*}
Adding these inequalities provides bounds on $p_i$ and $q_i$.
\begin{eqnarray*}
0 < a_0 + b_0 = 2 p_0 < \frac{37}{4} \left(7+4 \sqrt{3}\right) \\
-3(2+\sqrt{3}) < a_1 + b_1=2p_1 < 10+ \sqrt{3}
\end{eqnarray*}
Hence,  as $q_0=(p_0-b_0)/\sqrt{3}$,  we find the following bounds.
\begin{eqnarray*}
1\leq  p_0 \leq 64,  && -5\leq p_1 \leq 5 \\
\big[\frac{p_0-\frac{7}{4}}{\sqrt{3}}\big]  \leq   q_0 \leq \big[\frac{p_0}{\sqrt{3}}\big], && \big[\frac{p_1-2 }{\sqrt{3}}\big] \leq   q_1 \leq \big[\frac{p_1}{\sqrt{3}}\big].
\end{eqnarray*}
This gives us $2,688$ cases to consider.  In higher degree searches we will want to break this up,  but this presents no computational problems. We have a few conditions to satisfy.  First $\gamma$ is complex so $a_1^2<4 a_0$.  Next $b_1^2>4 b_0$ for two real roots and additionally,  for the two real roots to be in the correct interval,
\[  -1-\frac{\sqrt{3}}{2} < \frac{-b_1}{2} \pm \frac{1}{2}\sqrt{ b_1^2-4 b_0 } <0 \]
This gives us the condition
\[   2+\sqrt{3}  > b_1+  \sqrt{ b_1^2-4 b_0 }  \]
These simple tests reduce the search space to $202$ candidates. We next realise that $|\Im m(\gamma)|<3.5$,  giving the simple additional test $-49< a_1^2-4a_0 $.  There are only now only $18$ values for $\gamma$ satisfying our criteria and these are illustrated above in Figure X.   Of these all but $7$ are well outside our region and can be shown to be free on generators or alternatively one can check the factorisation as we shall now do.  The seven points and their minimal polynomial as as follows.
\begin{eqnarray*}
\gamma_1 =& \frac{1}{2} \left(\sqrt{3}+i \sqrt{5+4 \sqrt{3}}\right) ,& 1 + 6 z + z^2 + z^4.\\
\gamma_2 =& -1+i \sqrt{1+\sqrt{3}}  , & 1 + 8 z + 8 z^2 + 4 z^3 + z^4.\\
\gamma_3 =& \frac{1}{2} \left(\sqrt{3}+i \sqrt{13+8 \sqrt{3}}\right) ,&4 + 12 z + 5 z^2 + z^4. \\
\gamma_4 =& -1+i \sqrt{3+2 \sqrt{3}}, & 4 + 16 z + 12 z^2 + 4 z^3 + z^4. \\
\gamma_5 =& 1+\sqrt{3}+i \sqrt{3+2 \sqrt{3}} ,&1 + 20 z + 6 z^2 - 4 z^3 + z^4 .\\
\gamma_6 =&  \frac{1}{2} +\sqrt{3}+i \frac{1}{2}\sqrt{19+12\sqrt{3}}, & 16 + 32 z + 5 z^2 - 2 z^3 + z^4.\\
\gamma_7 =&  \frac{1}{2}\left(3 +\sqrt{3}+i \sqrt{8+6\sqrt{3}}\right), & 13 + 42 z + 4 z^2 - 6 z^3 + z^4.
\end{eqnarray*}
Then we compute the minimal polynomial for $\lambda_p$
\begin{eqnarray*}
\lambda_1: & & 33 - 240 z^2 + 238 z^4 - 16 z^6 + z^8.\\
\lambda_2 : & & 1 - 180 z^2 + 182 z^4 + 12 z^6 + z^8.\\
\lambda_3 : & & 97 - 464 z^2 + 446 z^4 - 16 z^6 + z^8. \\
\lambda_4 : & & 33 - 372 z^2 + 390 z^4 + 12 z^6 + z^8.\\
\lambda_5 : & &  193 - 1004 z^2 + 870 z^4 - 44 z^6 + z^8.\\
\lambda_6 : & & -11 - 588 z^2 + 890 z^4 - 36 z^6 + z^8.\\
\lambda_7 : & &-3 - 1032 z^2 + 1306 z^4 - 64 z^6 + z^8.
\end{eqnarray*}
Since these are all degree eight we conclude these groups are not arithmetic.

\medskip

\subsection{p=10} A similar situation occurs when $p=10$.  
\[ \left(\sqrt{3}+1\right) \sqrt{\gamma -\frac{\sqrt{3}}{2}+1} \]
After a similar search there there are three points we must consider,  the complex roots of the three equations
\begin{center}
\begin{tabular}{|c|c|} \hline
$\gamma$ polynomial &  $\lambda$ polynomial \\ \hline
$1+4 z+2 z^2+z^3+z^4$ &  $-25 - 400 z^2 + 165 z^4 - 10 z^6 + z^8$ \\ \hline
  $1+12 z-2 z^2-3 z^3+z^4$  &$ -1775 - 700 z^2 + 485 z^4 - 40 z^6 + z^8 $\\ \hline
 $1+11 z+6 z^2+z^3+z^4$  & $-725 - 1000 z^2 + 385 z^4 - 10 z^6 + z^8$ \\ \hline
\end{tabular} \\
\end{center} 
We again find there are no arithmetic lattices.
 
 \subsection{p=8}
 \[ \lambda_8=\sqrt{2 \left(\sqrt{2}+2\right) \gamma +2}\]
\begin{center}
\begin{tabular}{|c|c|} \hline
$\gamma$ polynomial &  $\lambda$ polynomial  \\ \hline
 $2+8 z+8 z^2+4 z^3+z^4$  &   $16 - 224 z^2 + 120 z^4 + 8 z^6 + z^8$ \\ \hline
$1+8 z+4 z^2+z^4$  & $ 272 - 704 z^2 + 328 z^4 - 16 z^6 + z^8$ \\ \hline
$1+6 z+7 z^2+2 z^3+z^4$ &  $ 496 - 736 z^2 + 256 z^4 + z^8 $ \\ \hline
$1+10 z+13 z^2+6 z^3+z^4$ &  $ 112 - 448 z^2 + 160 z^4 + 24 z^6 + z^8$ \\ \hline
$7+14 z+3 z^2-2 z^3+z^4$ &$ 368 - 928 z^2 + 544 z^4 - 32 z^6 + z^8$   \\ \hline
$7+20 z+14 z^2+4 z^3+z^4$ & $ 144 - 672 z^2 + 392 z^4 + 8 z^6 + z^8 $  \\ \hline
$4+20 z+5 z^2-2 z^3+z^4$ & $ 112 - 1120 z^2 + 656 z^4 - 32 z^6 + z^8 $  \\ \hline
$14+28 z+2 z^2-4 z^3+z^4$ &  $ 16 - 1088 z^2 + 856 z^4 - 48 z^6 + z^8 $ \\ \hline
$7+40 z+10 z^2-8 z^3+z^4$ &  $ 656 - 2784 z^2 + 1480 z^4 - 72 z^6 + z^8 $  \\ \hline
\end{tabular}
\end{center}
 \subsection{p=6}
 
 As noted earlier,  this case is covered by the results of \cite{MM6}.  The methods used here are similar and lead to the same results.  The candidates we find are
\begin{center}
 \begin{tabular}{|c|c|} \hline
$\gamma$ polynomial &  $\lambda$ polynomial  \\ \hline
$z^2 +6$ &  $15 - 6 z + z^2$ \\  \hline
$z^3 +3z^2 +6z+2$ &   $9 + 9 z - 3 z^2 + z^3$ \\ \hline
$ z^3 +5z^2 +9z+3 $ &   $-9 + 15 z - 3 z^2 + z^3$   \\ \hline
 $ z^4+8z^2+6z+1$ &    $-9 + 18 z + 12 z^2 - 6 z^3 + z^4$ \\ \hline
\end{tabular}  
\end{center} 
 
\subsection{The case $p=4$ and $q=5$, degree $2,3$ and $4$ over $\IQ(2\cos\frac{\pi}{5})$, total degree $4,6$ and $8$. }

In this subsection we meet the first really challenging search where it is likely that there may be some groups to find --- we have already found one with $\gamma$ real.  We shall find that there are no more. The intermediate field is
\[ L=\IQ(2\cos \frac{2\pi}{5})=\IQ(\sqrt{5}) \]
 We have the following absolute bounds on the modulus of $\gamma$ for groups which are not free;
\begin{equation}
|\gamma|<4  \big(\sqrt{2} \cos\frac{\pi}{5}+1 \big) + 2\cos \frac{2\pi}{5}  = 9.19453.
\end{equation}
This is achieved in the $(4,5,\infty)$-triangle group.  We also use the bounds
\begin{equation}
|\Im m(\gamma)|\leq 4, \;\;\; \Re e(\gamma) \geq  -4.
\end{equation}
 
We have already established that the total degree is no more than $8$ using the discriminant method.  The minimal polynomial for $\gamma$ now factors into two polynomials.  Both of these polynomials are either degree $2$, $3$ or $4$ and we consider each case separately.
\subsection{Degree 2 over $L=\IQ(\sqrt{5})$}
The minimal polynomial for $\gamma$ has the form
\[ p(z)=(z^2+a_1 z+ a_0)(z^2+b_1 z + b_0) \]
We always assume the first polynomial has a complex conjugate pair of roots.  Here $a_0,a_1,b_0,b_1\in \IQ(\sqrt{5})$ are algebraic integers and $b_i$ is the Galois conjugate of $a_i$.   
\begin{eqnarray*}
 a_0=\frac{p_0+q_0\sqrt{5}}{2}, && b_0=\frac{p_0-q_0\sqrt{5}}{2}\\
 a_1=\frac{p_1+q_1\sqrt{5}}{2}, && b_1=\frac{p_1-q_1\sqrt{5}}{2}
\end{eqnarray*}
where $p_i,q_i\in\IZ$ are integers with the same parity.
We observe that
\begin{eqnarray*}
p_i = a_i+b_i,\quad 
q_i = (p_i-2b_i)/\sqrt{5}
\end{eqnarray*}
Our number theoretic restrictions on $\gamma$ imply the first polynomial factor has roots $\gamma,\bar\gamma$ and the second factor has roots $r,s \in [-2\sin^2\frac{2\pi}{5},0]$.  Thus
\begin{eqnarray*}
a_0=|\gamma|^2, &&b_0=rs\\
a_1= - 2 \Re e[\gamma], &&b_1=-(r+s)
\end{eqnarray*}
From this we see that 
\begin{eqnarray*}
0< a_0 < 84.5393 , && 0< b_0 < 4 \sin^4\frac{2\pi}{5} =3.27254 \\
-18.3891  < a_1 < 8, && 0< b_1< 4 \sin^2\frac{2\pi}{5} =3.61803 .
\end{eqnarray*}
It is quite apparent that all these bounds can be improved with work. This is really only helpful when the searches are very large as we will see.
Hence we find the following bounds.
\begin{eqnarray*}
1\leq  p_0=a_0+b_0 \leq  87,  && 
\big[ \big(p_0-\frac{5}{4}(\sqrt{5}+3)\big)/\sqrt{5} \big] +1 \leq   q_0 \leq \big[p_0/\sqrt{5} \big] \\ -18 \leq p_1 =a_1+b_1  \leq 11 ,&&\big[(p_1-\sqrt{5}-1)/\sqrt{5}  \big] +1 \leq   q_1 \leq \big[p_1/\sqrt{5} \big].
\end{eqnarray*}
We are only interested in complex values for $\gamma$,  and $|\Im m(\gamma)|<4$, so $a_1^2-4a_0< -16$.  Also the condition that the second polynomial have two real roots in $[-2\sin^2\frac{2\pi}{5},0]$ gives the inequality 
\[    0<  b_1^2-4b_0  < (4\sin^2\frac{2\pi}{5} - b_1)^2 \]
It is this root condition that will be the most challenging test for our polynomials. After additionally checking parity conditions this gives $20$ polynomials whose complex roots are illustrated below.  They well indicate the problems we face. Figure 1. shows a rough picture or the exterior of the closed space of faithful discrete and free representations of the group $\IZ_4*\IZ_5$ in $PSL(2,\IC)$ and the two points we have found.

 \begin{center}
 \scalebox{0.55}{\includegraphics[angle=-90,viewport=100 300 500 500]{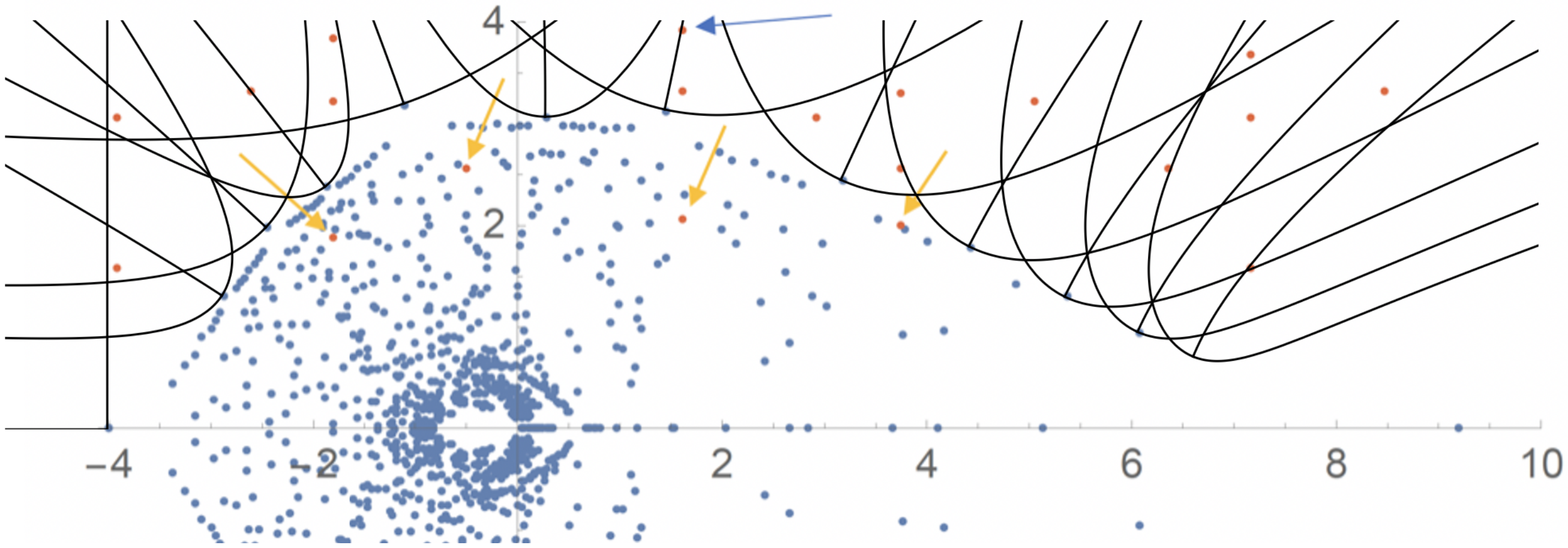}}\\
         \end{center}
         
\noindent{\bf Figure 4.}
{\em The $20$ potential $\gamma$ values. Four identified as not free,  or not discrete. One as a subgroup of an arithmetic group.}

\medskip

Of these $20$ points,  all but $4$ are well outside our region where non-free groups must lie and these $16$ are captured by neighbourhoods of pleating rays as per \cite{EMS2}.  One of these $16$ is in fact a discrete free subgroup of an arithmetic group and indicated by an arrow.  The remaining $4$ points are indicated with an arrow.

Notice that we are faced with the apparently real possibility of a point lying very close to the boundary and trying to decide if the associated group is free or not.  However at this point we cannot confirm any of these groups is discrete (other than those 16 we identified as free) as we have not used,  or satisfied, all the arithmetic information needed for the identification theorem. The remaining necessary and sufficient  condition we need to check is the factorisation condition.  In each case we set (where we have removed the factor $2\cos \frac{\pi }{5}$ from (\ref{lp}) as it is a unit).
 
 \begin{equation}\label{lp}
 \lambda_i=     \sqrt{2\gamma_i + 4\sin^2\frac{ \pi }{5}} 
\end{equation}
We know that it must be the case that $\IQ(\gamma_i)=\IQ(\lambda_i)$ and in particular the minimal polynomials for $\gamma_i$ and $\lambda_i$ must have the same degree, in this case degree $4$.  This is quite straightforward to check and we do it for {\em all} $20$ points.  The four points that might be at question are listed below with their minimal polynomials.

\medskip

{\small
\noindent\begin{tabular}{|c|c|l|l|}
\hline
& $\gamma $ & min. polynomial   $\gamma$ &   min. polynomial $\lambda$ \\
\hline
1 & $1.618 + 2.058 i $  &  $1 + 8 x + 3 x^2 - 2 x^3 + x^4$ & $181 - 226 x^2 + 87 x^4 - 14 x^6 + x^8$ \\
\hline
2 & $-0.5 + 2.569 i $ & $1 + 7 x + 8 x^2 + 2 x^3 + x^4$ & $171 - 144 x^2 + 37 x^4 - 6 x^6 + x^8$ \\
\hline
3 & $-1.809 + 1.892i$ &$1 + 10 x + 12 x^2 + 5 x^3 + x^4$ & $71 - 70 x^2 + 3 x^4 + x^8$ \\
\hline
4 & $3.736 + 1.996i $ & $1 + 26 x + 7 x^2 - 6 x^3 + x^4$ & $251 - 452 x^2 + 173 x^4 - 22 x^6 + x^8$\\
\hline
\end{tabular}
}

\medskip

Surprisingly there is the point $\gamma=1.61803 + 3.91487i$ (also indicated in Figure 4.) with minimal polynomial $x^4-2 x^3+14 x^2+22 x+1$ and with $\lambda= \sqrt{2\gamma + 4\sin^2\frac{ \pi }{5}}$ having minimal polynomial $x^4-6 x^3+11 x^2+4 x-19$. Both of degree $4$.  It is clear that $\gamma\in \IQ(\lambda)$ and hence the group generated by elliptics of order $4$ and $5$ with this value $\gamma$  for the commutator parameter is indeed a subgroup of an arithmetic Kleinian group.  However, as noted, it is captured by a pleating ray neighbourhood. We have established the following.

\begin{corollary} There are no arithmetic Kleinian groups generated by elements of order $4$ and $5$ with invariant trace field of degree $4$.
\end{corollary}

We now turn to the case of degree $6$, but offer fewer details.  We do note that the modulis space desciption and the capturing pleating ray neighbourhoods are independent of the degree of $\IQ(\gamma)$,  so we may use the same picture to remove possible values of $\gamma$.

\subsection{Degree 3 over $L=\IQ(\sqrt{5})$}
The minimal polynomial for $\gamma$ has the form
\begin{equation}\label{poly} p(z)=(z^3+a_2 z^2 +a_1 z+ a_0)(z^3+b_2 z^2+b_1 z + b_0) 
\end{equation} 
We continue to assume the first polynomial has a complex conjugate pair of roots.   We write,  for $p_i$ and $q_i$ of the same parity, 
\begin{eqnarray*}
2a_0=p_0+q_0\sqrt{5}, &&2b_0=p_0-q_0\sqrt{5}\\
2a_1=p_1+q_1\sqrt{5}, &&2b_1=p_1-q_1\sqrt{5} \\
2a_2=p_2+q_2\sqrt{5}, &&2b_2=p_2-q_2\sqrt{5} 
\end{eqnarray*}
The first polynomial factor has roots $\gamma,\bar\gamma,r$ with $r\in (-2\sin^2\frac{\pi}{5},0)$ and the second factor has real roots $r_1,r_2,r_3 \in [-2\sin^2\frac{2\pi}{5},0]$.  Thus
\begin{eqnarray*}
a_0=-|\gamma|^2 r, && b_0=-r_1r_2r_3.\\
a_1= |\gamma|^2+2 r   \Re e[\gamma], &&b_1=r_1r_2+r_2r_3+r_1r_3.\\
a_2= - 2 \Re e[\gamma]-r, &&b_2=-(r_1+r_2+r_3).
\end{eqnarray*}
From this,  with some easy estimates, we see that 
\begin{eqnarray*}
0< a_0 < 58.416, && 0< b_0 < 5.921\\
-1  < a_1< 97.246 , && 0< b_1<9.82 \\
-18.39 < a_2<  8.691, && 0< b_2<5.428
\end{eqnarray*}
As before, adding these inequalities provides bounds on $p_i$ and $q_i$.
\begin{eqnarray*}
1\leq p_0 \leq 65, \quad 
0\leq   p_1 \leq 107, \quad
-18 \leq   p_2 \leq 14 .
\end{eqnarray*}
\begin{eqnarray*}
\left[ \frac{ p_0-12}{\sqrt{5}} \right]+1 \leq \frac{ p_0-2b_0}{\sqrt{5}}= q_0 \leq \left[ \frac{ p_0}{\sqrt{5}} \right] \\
\left[ \frac{ p_1-20}{\sqrt{5}} \right]+1 \leq \frac{p_1-2b_1}{\sqrt{5}}= q_1\leq \left[ \frac{ p_1}{\sqrt{5}} \right] \\
\left[ \frac{ p_2-11}{\sqrt{5}} \right]+1 \leq \frac{p_2-2b_2}{\sqrt{5}}= q_2 \leq \left[ \frac{ p_2}{\sqrt{5}} \right] 
\end{eqnarray*}
Checking  parity alone gives us $6793308$ cases to consider.  Again we have a few conditions to satisfy. 
\begin{enumerate}
\item The first factor must be negative at $-2\sin^2\frac{\pi}{5}$.
\[ a_0+\frac{1}{4} (\sqrt{5}-5) a_1-\frac{5}{8}\big((\sqrt{5}-3) a_2 -2 \sqrt{5}+5\big) <0 \]
\item The second factor must be negative at $-2\sin^2\frac{2\pi}{5}$ and have  derivative with two roots in the interval $[-2\sin^2\frac{2\pi}{5},0]$
\begin{eqnarray*}
b_0-\frac{1}{4} (\sqrt{5}+5)b_1+\frac{5}{8} \left((\sqrt{5}+3)b_2-2 \sqrt{5}-5\right)<0\\
\frac{3}{4}  (\sqrt{5}+5 )> b_2+ \sqrt{b_2^2-3b_1}  
\end{eqnarray*}
\item The discriminant of the first factor is negative and the second factor positive.
\begin{eqnarray*}
0 > -27 a_0^2+a_1^2 \left(-4 a_1+a_2^2\right)+2a_0 \left(9 a_1a_2-2 a_2^3\right)  \\
0 <  -27 b_0^2+b_1^2 \left(-4 b_1+b_2^2\right)+2b_0 \left(9 b_1b_2-2 b_2^3\right)  
\end{eqnarray*}
\end{enumerate} 

If these tests are all passed, as they were by $114$ candidates,  then we numerically computed the roots to identify $\gamma$ and check that all requirements were satisfied.  This time there were only three points clearly lying outside the discrete free space.   However we found the minimal polynomial in $\lambda$ for all these cases and found it to be degree $12$.  This proves

\begin{corollary} There are no arithmetic Kleinian groups generated by elements of order $4$ and $5$ with invariant trace field of degree $6$.
\end{corollary}

\subsection{Degree 4 over $L=\IQ(\sqrt{5})$}
The minimal polynomial for $\gamma$ has the form
\begin{equation}\label{poly} p(z)=(z^4+a_3z^3+a_2 z^2 +a_1 z+ a_0)(z^4+b_3z^3+b_2 z^2+b_1 z + b_0) 
\end{equation} 
We continue to assume the first polynomial has a complex conjugate pair of roots.   We write,  for $p_i$ and $q_i$ of the same parity, 
\begin{eqnarray*}
2a_0=p_0+q_0\sqrt{5}, &&2b_0=p_0-q_0\sqrt{5}\\
2a_1=p_1+q_1\sqrt{5}, &&2b_1=p_1-q_1\sqrt{5} \\
2a_2=p_2+q_2\sqrt{5}, &&2b_2=p_2-q_2\sqrt{5} \\
2a_3=p_3+q_3\sqrt{5}, &&2b_3=p_3-q_3\sqrt{5}  
\end{eqnarray*}
 Following our earlier notation, direct estimates yield
 \begin{eqnarray*}
0<  b_0  < 10.71 \quad  
0< b_1<23.69, \quad 0< b_2 <19.64 \quad 0< b_3 < 7.24 \\
0< a_0< 40.364 , \quad -9.382<a_3<18.39 .
\end{eqnarray*}
 It is worth a little time to improve the obvious bounds on $a_1$ and $a_2$.  The polynomial 
 \[ q(x)=x^4+a_3x^3+a_2x^2+a_1x+a_0 \]
 has a complex conjugate pair of roots and two roots in the interval $[-2\sin^2\frac{\pi}{5},0]$.  Thus $q'(0)=a_1>0$ and $q'(-2\sin^2\frac{\pi}{5})<0$. These conditions are not sufficient.  Using elementary calculus,  which we leave the reader to consider,  one also obtains
 \[ 0 \leq a_1 <112.69, \quad\quad  -1.76 < a_2 < 88.03 \]
Hence 
\begin{eqnarray*}
1\leq p_0 \leq 40, \quad 
0\leq   p_1 \leq 114, \quad
-1 \leq   p_2 \leq 90, \quad -9 \leq   p_3 \leq 21. 
\end{eqnarray*}
\begin{eqnarray*}
\left[ \frac{ p_0-21.5}{\sqrt{5}} \right]+1 \leq \frac{ p_0-2b_0}{\sqrt{5}}= q_0 \leq \left[ \frac{ p_0}{\sqrt{5}} \right] \\
\left[ \frac{ p_1-48}{\sqrt{5}} \right]+1 \leq \frac{p_1-2b_1}{\sqrt{5}}= q_1\leq \left[ \frac{ p_1}{\sqrt{5}} \right] \\
\left[ \frac{ p_2-40}{\sqrt{5}} \right]+1 \leq \frac{p_2-2b_2}{\sqrt{5}}= q_2 \leq \left[ \frac{ p_2}{\sqrt{5}} \right] \\
\left[ \frac{ p_3-15}{\sqrt{5}} \right]+1 \leq \frac{p_3-2b_3}{\sqrt{5}}= q_3 \leq \left[ \frac{ p_3}{\sqrt{5}} \right] 
\end{eqnarray*}
This search is now about $10^{10}$ cases after checking parity.  To make this search more efficient we need to find tests on coefficients as soon as possible in the process.  The most obvious target is the fact that the polynomial
\[ p(z)=z^4+b_3z^3+b_2z^2+b_1z+b_0 \]
has four real roots in the interval $I_5=[-2\sin^2\frac{2\pi}{5},0]$. Then 
\begin{eqnarray*}
p'(z) & = & 4z^3+3 b_3z^2+2 b_2z +b_1, \quad \mbox{has three real root in $I_5$} \\
p''(z) & = & 12z^2+6 b_3z+2 b_2, \quad \mbox{has two real root in $I_5$} \\
p'''(z) & = & 24z+6 b_3, \quad \mbox{has a real root in $I_5$}  
\end{eqnarray*}
This last condition gives $-2\sin^2\frac{2\pi}{5}<-b_3/4<0$, which implies  $0<b_3< 8\sin^2\frac{2\pi}{5}$ which we have by construction.  The condition on $p''(z)$ implies
\[   0< \frac{b_3}{4} \pm \frac{\sqrt{9 b_3^2-24 b_2}}{12} < 2 \sin^2\frac{2\pi}{5} \]
The content here is that 
\begin{eqnarray*}
0  < & 3 b_3^2-8  b_2 & < 3\big(8 \sin^2\frac{2\pi}{5}-b_3\big)^2    
\end{eqnarray*} 
Now the cubic $p'(z)$ has three real roots in $I_5$.  It's discriminant must satisfy 
\[ 108 b_1^2+27b_1b_3^3+32 b_2^3<108 b_1b_2b_3+9b_2^2 b_3^2 \]
and $p'(-2\sin^2\frac{2\pi}{5})<0$, 
\[ b_1-\frac{1}{2}  (\sqrt{5}+5 )b_2+\frac{5}{8}  (3 \sqrt{5} b_3+9b_3-8 \sqrt{5}-20 ) <0 \]
Finally,  for the quartic itself we must have the discriminant positive and at least one real root.   

\medskip

After this search we still do not yet have candidate polynomials,  just a much shorter list of possibilities.  It is hard to assert that the two real roots of the first polynomial   lie in the correct interval without actually computing them.  Now that we have a shorter list of possibilities we do just that.  We compute the roots of  $q(x)=x^4+a_3x^3+a_2x^2+a_1x+a_0$ and check the imaginary part of the complex roots are smaller than $5$ is absolute value and that the real roots lie in the shorter interval.  This left us with $9$ possibilities.
{\small \begin{eqnarray*}
 \frac{1}{2} \left(7+3 \sqrt{5}\right)+\frac{1}{2} \left(45+19 \sqrt{5}\right) z+\frac{1}{2} \left(56+22 \sqrt{5}\right) z^2+\frac{1}{2} \left(-9-7 \sqrt{5}\right) z^3+z^4\\ \frac{1}{2} \left(18+8 \sqrt{5}\right)+\frac{1}{2} \left(67+29 \sqrt{5}\right) z+\frac{1}{2} \left(56+22 \sqrt{5}\right) z^2+\frac{1}{2} \left(-9-7 \sqrt{5}\right) z^3+z^4\\ \frac{1}{2} \left(7+3 \sqrt{5}\right)+\frac{1}{2} \left(37+15 \sqrt{5}\right) z+\frac{1}{2} \left(40+14 \sqrt{5}\right) z^2+\frac{1}{2} \left(-6-6 \sqrt{5}\right) z^3+z^4\\ \frac{1}{2} \left(7+3 \sqrt{5}\right)+\frac{1}{2} \left(44+18 \sqrt{5}\right) z+\frac{1}{2} \left(47+17 \sqrt{5}\right) z^2+\frac{1}{2} \left(-6-6 \sqrt{5}\right) z^3+z^4\\ \frac{1}{2} \left(3+\sqrt{5}\right)+\frac{1}{2} \left(28+10 \sqrt{5}\right) z+\frac{1}{2} \left(38+12 \sqrt{5}\right) z^2+\frac{1}{2} \left(-3-5 \sqrt{5}\right) z^3+z^4\\ \frac{1}{2} \left(18+8 \sqrt{5}\right)+\frac{1}{2} \left(67+29 \sqrt{5}\right) z+\frac{1}{2} \left(70+28 \sqrt{5}\right) z^2+\frac{1}{2} \left(16+4 \sqrt{5}\right) z^3+z^4\\ \frac{1}{2} \left(3+\sqrt{5}\right)+\frac{1}{2} \left(43+17 \sqrt{5}\right) z+\frac{1}{2} \left(75+29 \sqrt{5}\right) z^2+\frac{1}{2} \left(19+5 \sqrt{5}\right) z^3+z^4\\ \frac{1}{2} \left(7+3 \sqrt{5}\right)+\frac{1}{2} \left(58+24 \sqrt{5}\right) z+\frac{1}{2} \left(86+34 \sqrt{5}\right) z^2+\frac{1}{2} \left(19+5 \sqrt{5}\right) z^3+z^4\\ \frac{1}{2} \left(7+3 \sqrt{5}\right)+\frac{1}{2} \left(55+23 \sqrt{5}\right) z+\frac{1}{2} \left(90+36 \sqrt{5}\right) z^2+\frac{1}{2} \left(19+5 \sqrt{5}\right) z^3+z^4
\end{eqnarray*}
}
All of these cases are well within the space of free representations,  but in fact none of them satisfied the factorisation criterion.

\begin{theorem} There are no arithmetic Kleinian groups generated by elements $f$ of order $4$ and $g$ of order $5$ with $\gamma(f,g)\not\in\IR$.
\end{theorem}

In fact we have actually proven quite a bit more as we have shown every subgroup of an arithmetic Kleinian group generated by  $f$ of order $4$ and $g$ of order $5$,  if discrete and if $\gamma$ is complex, is free.  We therefore have the following theorem.

\begin{theorem} Let $\gamma$ be an arithmetic Kleinian group.  Suppose  $f$ of order $4$ and $g$ of order $5$ lie in $\Gamma$.  Then either $\langle f,g\rangle$ is free, or $\Gamma$ contains $Tet(4,5;3)$ with finite index. 
\end{theorem}

\section{$p=2,3$ and $4$.}

Apart from the case $p=6$ already covered,  it is here that we will find many examples.  Unfortunately there is not necessarily any intermediate field to help us in our searches.  The fact that earlier we had a polynomial factor with all real roots was enormously helpful. Therefore the cases here are of a quite different nature.  As we will see, the cases $p=2,4$ can be dealt with together as they will be commensurable,  and the polynomials we seek are largely covered by the results of Flammang and Rhin \cite{FR}.  The complexity of that case is in getting a very accurate description of the moduli space. That leaves the case $p=3$ which we now deal with.

\subsection{$p=3$.}

We first identify a useful total degree bound.  This is obtained as a modification of the method of auxilliary functions introduced by Smyth and used in this context by Flammang and Rhin, \cite{FR}.  We will have more to say about this a bit later.

Let $\alpha = 1+\gamma$ and note that then 
\begin{equation}\label{alphabounds} |\alpha-1|<3+2 \sqrt{2}, \quad |\alpha|<4+2 \sqrt{2} \end{equation}
\bigskip
Let $t_0$ be the unique solution to the equation 
\[  \left|\frac{t}{1+t}\right|^t\frac{1}{1+t}=\frac{3}{2^{t+1}} \]
Then $t_0\approx 4.28291$.  A little calculus reveals that the function $|x|^t|1-x|\leq 0.0770561$ on the interval $[-\frac{1}{2},1]$.  The graph is illustrated below.
  \scalebox{0.5}{\includegraphics[viewport=-40 520 570 800]{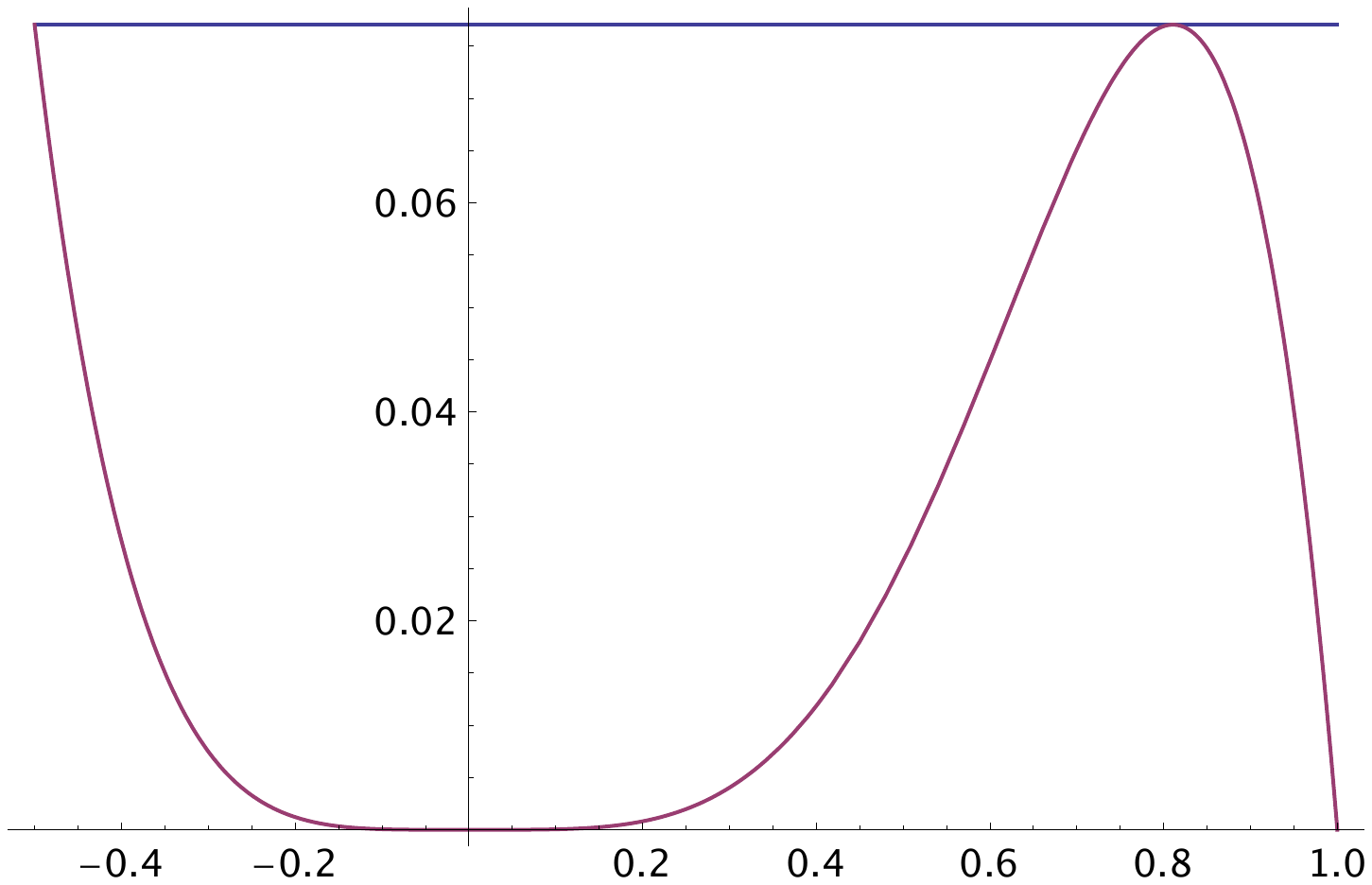}}\\
  
  \medskip
\noindent  {\bf Figure 5.} {\em The graph of $|x|^t|1-x|$ with $t=4.28291$.}
  
  \medskip

Let 
\[ q_\alpha(z)=z^n+a_{n-1}z^{n-1}+\cdots+a_1z+a_0 = (z-\alpha)(z-\bar\alpha)(z-r_1) \cdots(z-r_{n-2}) \]
be the minimal polynomial for $\alpha$.  Then,  as $r_i\in [-\frac{1}{2},1]$,
\begin{eqnarray} |q_\alpha(0)|^{t_0}|q_\alpha(1)| &=& |\alpha|^{2t_0}|\alpha-1|^2 \prod_{i=1}^{n-2}  |r_i|^{t_0} |1-r_{i}|\nonumber  \\
& \leq & (3+2 \sqrt{2})^2\Big(4+2 \sqrt{2}\Big)^{2t_0} (0.0770561)^{n-2}. \label{a0bound}
\end{eqnarray}
As the left hand side here is greater than one,  we deduce that $n\leq 9$.  Indeed if $p(0)\neq \pm 1$,  then $n\leq 8$.  In fact we have the following lemma which will be useful in our high degree searches.

\begin{lemma}\label{*lemma}  Let $\IQ(\alpha)$ be a field with one complex place, $\sigma(\alpha)\in [-\frac{1}{2},1] $ for all real embeddings,  and satisfying (\ref{alphabounds}). The the degree of this field is no more than $9$, and if it is $9$,  then $\alpha$ is a unit, and $1\leq q(1) \leq 7$,  where $q$ is the minimal polynomial for $\alpha$.  If the degree is $8$,  then $1\leq |q(0)|\leq 2$, and $|q(0)|=2$ implies $|q(1)|\leq 5$. Further,  in the degree $9$ case we also have $|\alpha|\geq 5.28769$,  and in the degree $8$ case $|\alpha|\geq 4.1139$.
\end{lemma}
For the last two bounds,  we also use the fact that $|\alpha-1|=|\gamma|$,  and our bounds on the values $|\gamma|$ when $\gamma=\gamma(f,g)$ and $\langle f,g\rangle$ is not free.

\bigskip

Having these degree bounds now allows us to search to find all the possibilities for the values $\gamma$ satisfying the criteria determined by the Identification Theorem \ref{2genthm}.  Some of these searches are at first sight barely feasible,  for instance the degree $9$ case.  However,  the factorisation condition described earlier significantly reduces the possibilities by imposing strong parity conditions on the coefficients of the minimal polynomial of 
\begin{equation} \lambda=\sqrt{2\alpha+1}=\sqrt{2\gamma+3}, \quad \IQ(\lambda)=\IQ(\alpha)=\IQ(\gamma)
\end{equation} as we now describe.  We discuss briefly the low degree cases of degree $2$ and $3$ as there are straightforward.  We then go into rather more detail in degrees $4$ and $5$.

\section{The case of degree 2.}
 
We have $\gamma^2+b\gamma+c=0$,  and $a_3=\sqrt{2\gamma+3}$ must also be a quadratic integer by the factorisation condition.  Also, using easier bound on the size of the moduli space, $\Im m(\gamma)<3$ and $-4<\Re e(\gamma) < 6$, and as $\gamma$ must be complex,  $4c>b^2$.  This gives us $\gamma=\frac{1}{2}(-b+\sqrt{b^2-4c})$ and hence
\[  8> b > -12, \quad  0>  b^2-4c >- 36 \] 
This gives $162$ cases to consider,  many of which are not in the moduli space. However once we implement the factorisation test that the minimal polynomial of $a_3$ is degree 2 this quickly reduces to the three complex values
$-2+\sqrt{2}i$,  $2i$ and $-3+2i$.  The first, $-2+\sqrt{2}i$,  is found from $(3,0)$, $(4,0)$ orbifold Dehn surgery on the Whitehead link. The second, $2i$, is the associated Heckoid group with $(W_{3/4})^2=1$ and the last group is free on its generators.  One can establish this later fact directly or it is a consequence of the results described in Figure 6 below. 

\section{The case of degree 3.}
We consider the polynomial $p(z)=z^3+az^2+bz+c$ which must have exactly one real root $r\in [-3/2,0]$.  Then
\[ p(z)= -\bar \gamma \gamma r + (\bar \gamma \gamma   + \bar\gamma r   + \gamma r )z  - 
 (\bar\gamma + \gamma + r)+ z^3 \]
With our bounds on $r$ and $\gamma$ we find
\[  -9 \leq a \leq 10, \quad 1\leq b \leq 25, \quad 1\leq c \leq  36,\]
As $p(0)=c>0, p'(0)=b>0$ and $p(-3/2)<0$ and $p'(-3/2)>0$ and the discriminant must be negative a simple loop with these tests finds about $5,500$ candidates. One the factorisation condition is applied to these potential candidates only $17$ remain.  This is illustrated below with the $4$ points outside the moduli space indicated. Proving that these are the only such points in the complement of the moduli space of groups free on generators follows from the arguments we provide in the next case of degree $4$ which will give a provable computational description of these spaces.  \\
  
  \scalebox{0.5}{\includegraphics[angle=-90,viewport=10 80 550 550]{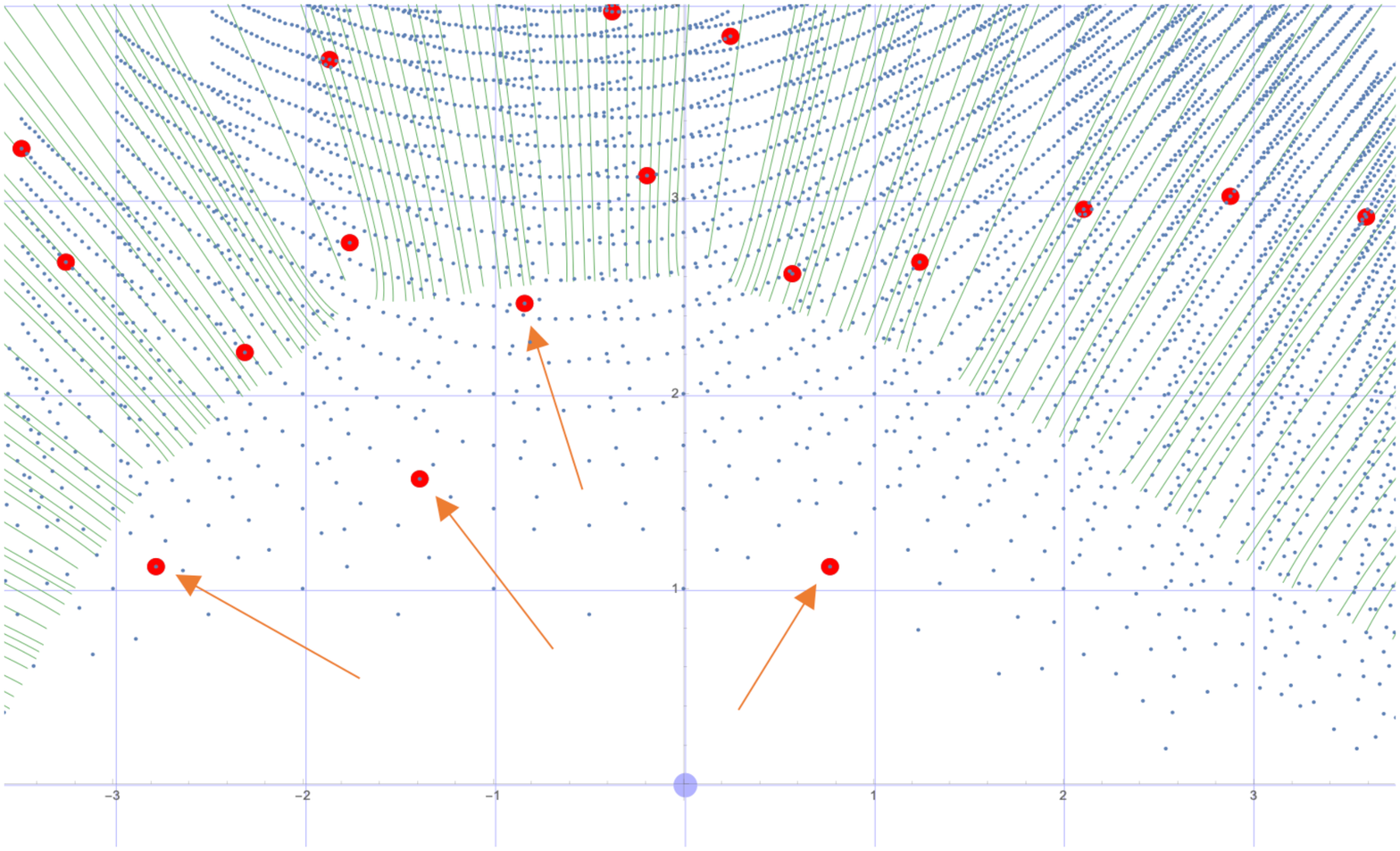}}\\

\noindent{\bf Figure 6.} {\em The possible $5,500$ candidate $\gamma$ values in the case of the degree $3$ search. Of these only  $17$ satisfy the factorisation condition and only $4$ are in the complement of the moduli space of groups free on generators of order $3$ and $4$.}

\bigskip

We now describe the degree $4$ and $5$ cases in some detail for two reasons.  First in degree $4$ we find several cases which do turn out to generate arithmetic lattices, and then in degree $5$ there are none,  but we need to develop more ideas to limit the search spaces.  Finally we carefully eliminate the highest degree search,  degree $9$.  We have searched the cases degree $6,7$ and $8$.  They all go the same way and there is nothing to report.
 
\section{The case of degree 4.}

The minimal polynomial for $\gamma$ and hence $\lambda$ has degree $4$,
 \[\lambda^4+a_3\lambda^3+a_2\lambda^2+a_1\lambda+a_0 = 0 \]
Thus 
\[(2\alpha+1)^2+a_2(2\alpha+1)+a_0 = -\sqrt{2\alpha+1}(a_1 +a_3(2\alpha+1) \]
Squaring both sides and rearranging gives us a polynomial for $\alpha$.

\begin{eqnarray*}
0& = & 16 \alpha ^4+\left(-8 a_3^2+16 a_2+32\right) \alpha ^3+\left(4 a_2^2+24 a_2-12 a_3^2+8 a_0-8 a_1 a_3+24\right) \alpha ^2 \\ && +\left(-2 a_1^2-8 a_3 a_1+4 a_2^2-6 a_3^2+8 a_0+4 a_0 a_2+12 a_2+8\right) \alpha \\&& +\left(a_0^2+2 a_2 a_0+2 a_0-a_1^2+a_2^2-a_3^2+2 a_2-2 a_1 a_3+1\right) 
\end{eqnarray*}
There is a unique monic polynomial for the algebraic integer $\alpha$ over $\IQ$ and so all these coefficients are divisible by $16$.
\begin{enumerate}
\item degree 3 coefficient: $a_3^2$ is even,  so \underline{$a_3$ is even}.
\item degree 2 coefficient: $  a_2^2+2a_2+2 a_0+2 $ is divisible by $4$. Thus \underline{$a_2$ is even}.  Then 
$2 a_0+2 $ is divisible by $4$,  so \underline{$a_0$ is odd}.
\item degree 1 coefficient:  $- a_1^2-3 a_3^2+4 a_0+2 a_0 a_2+6 a_2+4$ is divisible by $8$.   We reduce mod $2$ to see  \underline{$a_1$ is even}.
\item degree 0 coefficient  $\left(1+a_0-a_1+a_2-a_3\right) \left(1+a_0+a_1+a_2+a_3\right)$.  We reduce mod $2$ to see  \underline{$a_0$ is odd}
 \end{enumerate}
We now make the substitutions $a_i=2k_i$ or $2k_i+1$ as the case may be.  We expand out the polynomial for $\alpha$ and check parities again.  We find $k_2$ is even and $k_0$ is odd, and the sum $k_1+k_3$ is even.  Again write $k_2=2m_2$, $k_0=2m_0-1$ and $k_3=2m_3-k_1$.  Substitute back and expand out.  We have
\begin{align*}
&16 \left(\left(m_0+m_2\right){}^2-m_3^2\right)\\
&+32 \left(\left(1+m_2\right) \left(m_0+m_2\right)+k_1 m_3-3 m_3^2\right) \alpha \\
&+16 \left(1+2 m_0+m_2 \left(4+m_2\right)-\left(k_1-6 m_3\right) \left(k_1-2 m_3\right)\right) \alpha ^2\\
& -32 \left(-1-m_2+\left(k_1-2 m_3\right){}^2\right) \alpha ^3 +16 \alpha ^4
\end{align*}
as $16$ times the minimal polynomial.  Further,  we have proved the following.
\begin{lemma}  The minimal polynomial for $\lambda$ has the following coefficient structure.  There are integers $n_0,n_1,n_2$ and $n_3$ such that
\begin{itemize}
\item $a_0=-1 + 4 n_0$.
\item $a_1= 2n_1$.
\item $a_2=4n_2$.
\item $a_3=-2n_1 + 4 n_3$.
\end{itemize}
Moreover the coefficient structure of the minimal polynomial for $\alpha$ has (with the same $n_i$)
\begin{itemize}
\item $b_0=\left(n_0+n_2\right){}^2-n_3^2$.
\item $b_1= 2 \left(\left(1+n_2\right) \left(n_0+n_2\right)+n_1 n_3-3 n_3^2\right)$.
\item $b_2=2 \left(\left(1+n_2\right) \left(n_0+n_2\right)+n_1 n_3-3 n_3^2\right)$.
\item $b_3=-2 \left(-1-n_2+\left(n_1-2 n_3\right){}^2\right)$.
\end{itemize}
Note that also $b_1$ and $b_3$ are even,  $b_2$ has the same parity as $b_3/2$.
\end{lemma} 
With this information we can new set up a simple search.  To do so we use Vieta's formulas.  First for $\alpha$ whose real roots are in $[-\frac{1}{2},1]$.

First note that once we fix $\alpha$ the coefficients of its minimal polynomial  are monomials in the real roots $r_i$. Thus either increasing or decreasing.  Thus the extrema occur at $r_i=1$ or $r_i=-\frac{1}{2}$ so it is just a small finite problem to find these. We replace $\alpha+\bar\alpha$ by either $|\alpha|=4+2\sqrt{2}$ or $-4$. This is unchallenging in low degree,  but a bit of work in higher degree.

\medskip

 \begin{tabular}{|l|c|c|} \hline 
$b_0=|\alpha|^2 r_1r_2$ &  $-16 \leq b_0\leq 33$ \\ \hline
$b_1=-|\alpha|^2(r_1+r_2)-(\alpha+\bar\alpha)r_1r_2 =2x_1$&  $-53 \leq x_1\leq 24$ \\ \hline
$b_2=|\alpha|^2+(\alpha+\bar\alpha)(r_1+r_2)+r_1r_2 $&  $0 \leq b_2\leq 74$ \\ \hline
$b_3 =-( \alpha+\bar\alpha+r_1+r_2)=2x_3 $ &  $-7 \leq x_3 \leq 4$ \\ \hline
\end{tabular}

\medskip

With the one other parity restriction we have found this is a search of about $1.5$ million cases.  

For the $\lambda$ search we may assume additionally that $\lambda+\bar\lambda \geq 0$ with the choice of square root to find the search bounds (this time $r_1,r_2\in [-\sqrt{3},\sqrt{3}]$).  Recall $\lambda=\sqrt{2\alpha+1}$ and $|\lambda|<\sqrt{9+4\sqrt{2}}\approx 3.828$.

\medskip

{\small
 
\noindent \begin{tabular}{|l|c|c|c|} \hline 
$a_0=|\lambda|^2 r_1r_2$&  $-43  \leq a_0\leq  43 $ & $ -10 \leq n_0 \leq 11 $ \\ \hline
$a_1=-|\lambda|^2(r_1+r_2)-(\lambda+\bar\lambda)r_1r_2$&  $-73  \leq a_1\leq 27$& $ -36 \leq n_1 \leq 13$\\ \hline
$a_2=|\lambda|^2+(\lambda+\bar\lambda)(r_1+r_2)+r_1r_2 $&  $-2 \leq a_2\leq 30$& $ -1 \leq n_2 \leq 7$ \\ \hline
$a_3=-( \lambda+\bar\lambda+r_1+r_2)$&  $-11\leq a_3\leq 3 $& $\Big[\frac{2n_1-11}{4}\Big] \leq n_3\leq\Big[ \frac{3+2n_1}{4}\Big] $ \\   \hline
\end{tabular} }

\medskip

 This is a search of about $22$ thousand cases - nearly two orders of magnitude smaller.  So we of course choose to do the latter search.  However as the degree grows,  the symmetric functions start to get quite big quite quickly (for $b_i$ they are evaluated at $-\frac{1}{2}$ or $1$,  while for $a_i$ they are evaluated at $\pm \sqrt{3}$) and from about degree $7$ or $8$  it is best to search through the coefficients $a_i$.   For these larger degrees it is worthwhile obtaining much better coefficient bounds as well.
 
 \medskip
 
 Within these bounds on the coefficients we tested the real root condition,  that $|\lambda|<4$ and that $|\gamma|<5.85$ and $|\Im m(\gamma)|<3.75$ we found $55$ possibilities for $\gamma=\alpha-1$ up to complex conjugacy.  These are illustrated below.

  \scalebox{0.5}{\includegraphics[viewport=-40 520 570 800]{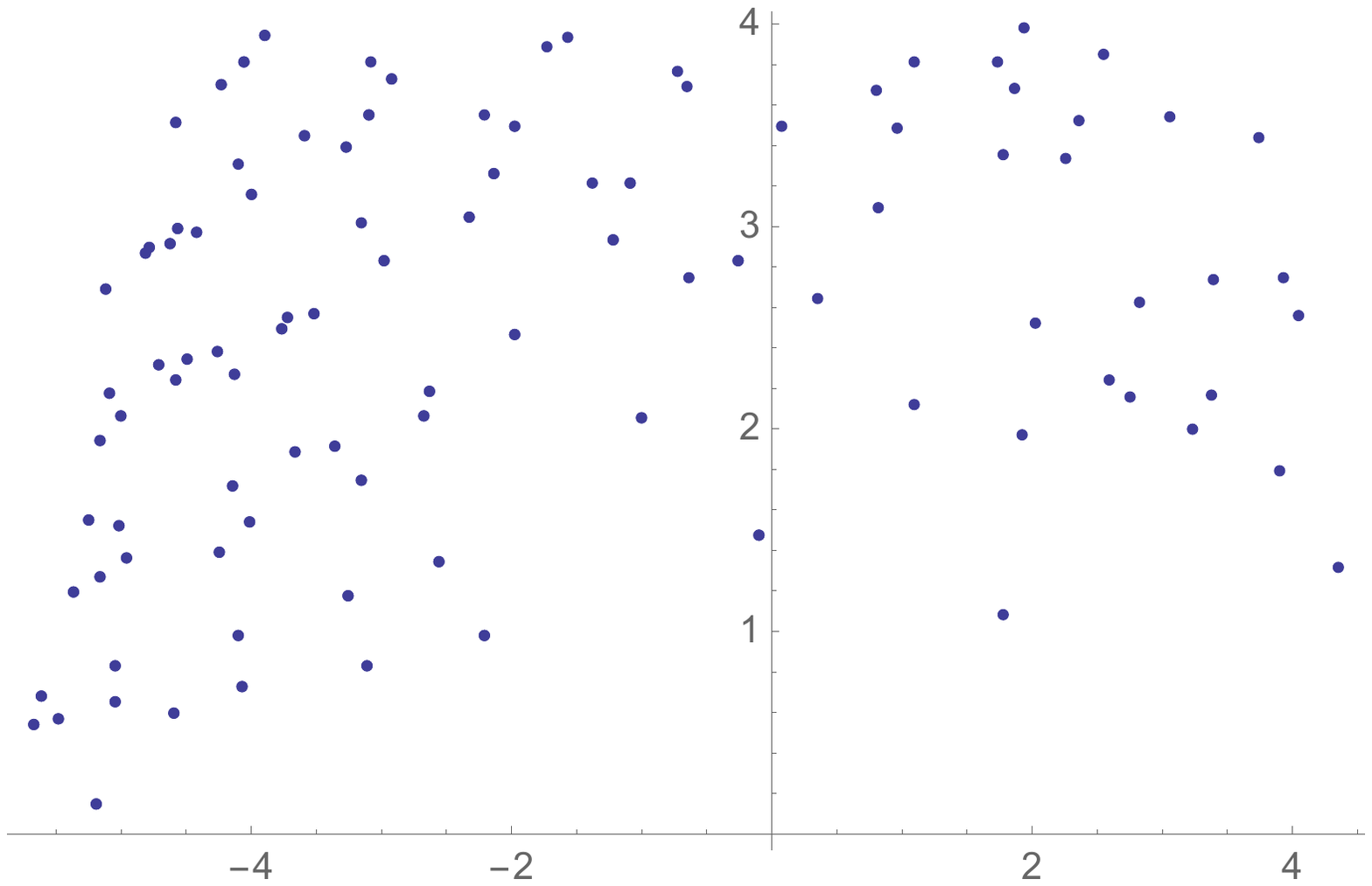}}\\
  
  \medskip
\noindent{\bf Figure 7.} {\em The degree $4$ candidate values for $\gamma$ satisfying all arithmetic criteria}

\medskip 

We now have to decide whether these groups are free.  Each is definitely a discrete subgroup of an arithmetic Kleinian group at this point.

  \scalebox{0.85}{\includegraphics[viewport=100 520 570 800]{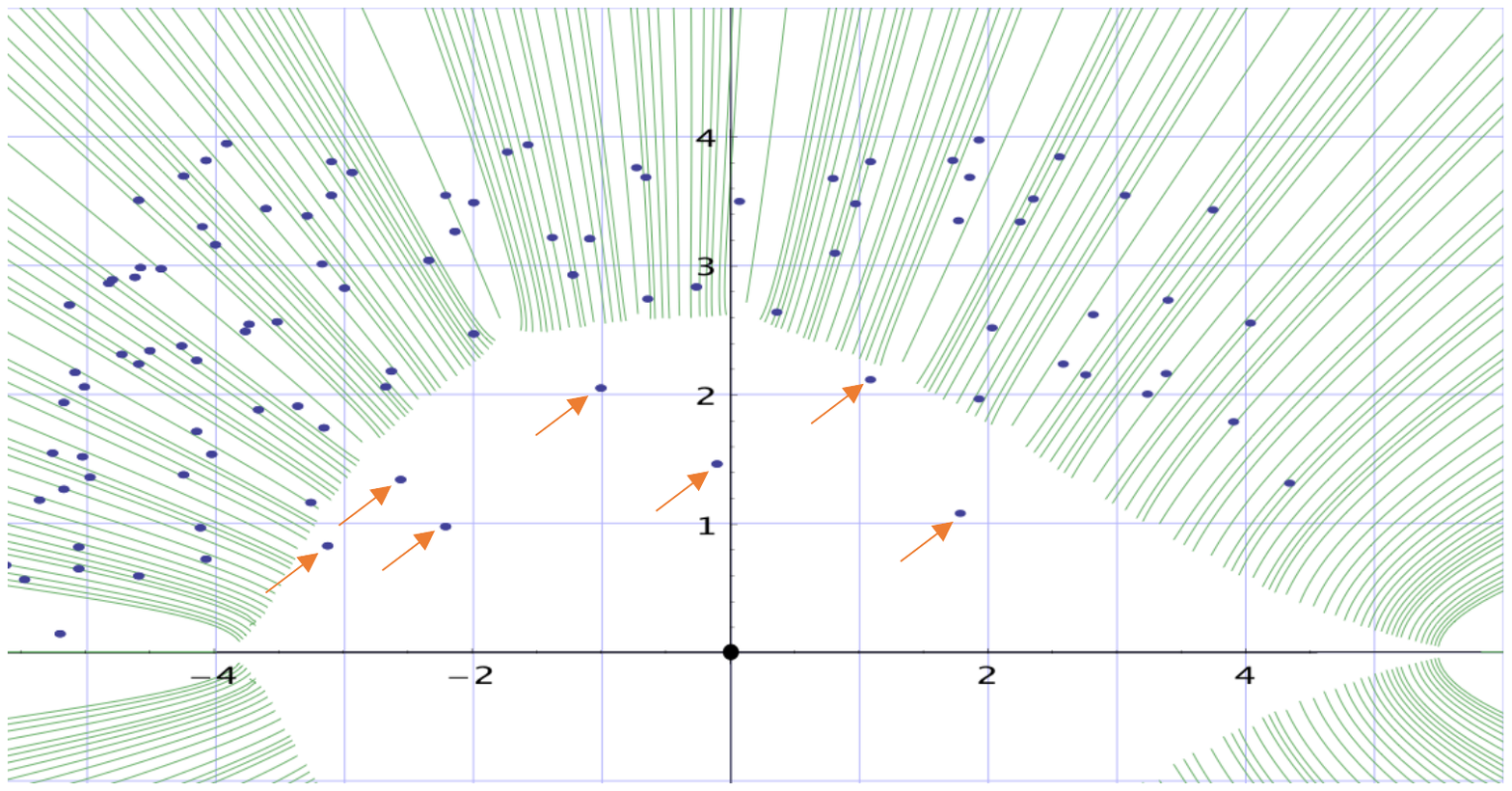}}\\
  \noindent{\bf Figure 8.} {\em The degree $4$ candidate values for $\gamma$ satisfying all arithmetic criteria mapped against the moduli space.}
  
  \medskip
  
  Now to eliminate these points we first examine the $17$ Farey polynomials and the pleating rays of low slopes where the denominator (which is the degree of the polynomial) is small.
  {\small \begin{equation}\label{slp}\{1/2, 3/5,4/7,6/7,5/8,5/9,7/11,8/13,9/14,11/17,13/20,16/21,15/26\}.\end{equation}}
The inverse image of the half-space $H=\{z=x+i y: x\leq -2\}$ under the branch of the Farey polynomial $P_{r/s}$ which yields the rational pleating ray $r/s$ is proved in \cite{EMS2} to consist only of geometrically finite groups free on the enerators of order $3$ and $4$.  For these low slopes these preimages capture all but $4$ points easily.   There are some computation issues in drawing the figure below for the last $4$ polynomial inverse images,  it is really only for visual confirmation.  In fact all we need to do is evaluate the associated polynomial with integer coefficients on the point in question, show that the image lies in $\{z=x+i y: x\leq -2\}$,  then just check that we have the right branch - but this is determined by the pleating ray and path lifting which are direct to check as we can numerically identify the other roots and critical points.
 
  \scalebox{0.5}{\includegraphics[angle=-90,viewport=-40 100 550 550]{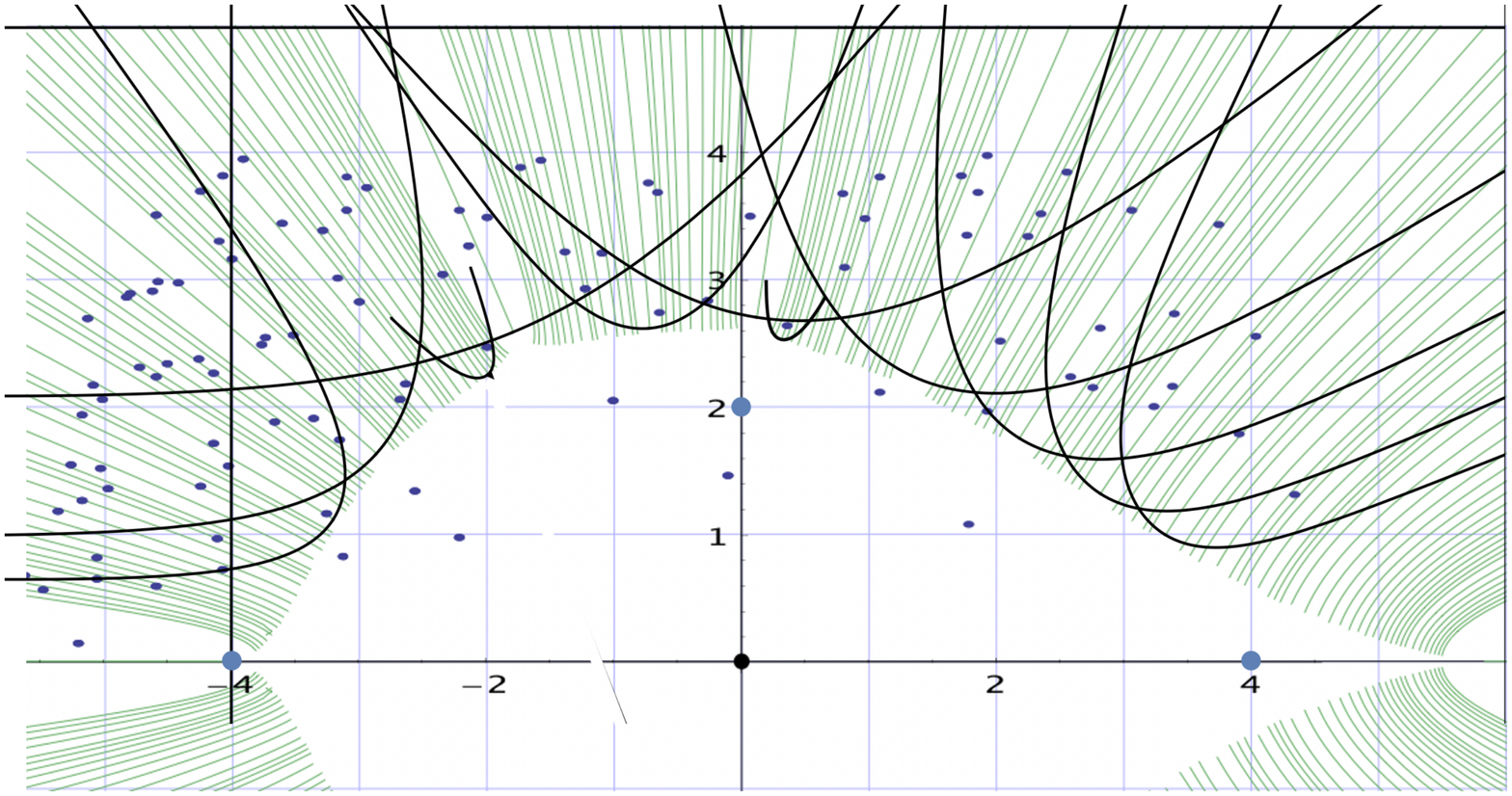}}\\

\noindent{\bf Figure 9.} {\em Neighbourhoods of the pleating rays with slopes give at (\ref{slp}) capturing all but $7$ points. }

\medskip

There are now $7$ points that remain and the associated groups are not going to be free on the two generators.  These appear in the table below.  We will explain how these groups are determined later.
 
\section{Degree $5$ search.}

As these searches get bigger it is efficient to look for further simple tests to eliminate possible candidates before we numerically calculate the roots.  We first look at the parity considerations as above.   Following the above strategy,  we assume minimal polynomial for $\lambda$ has degree $5$,
 \[p(\lambda)=\lambda^5+a_4+a_3\lambda^3+a_2\lambda^2+a_1\lambda+a_0 = 0 \]
 Thus 
\[a_4(2\alpha+1)^2+a_2(2\alpha+1)+a_0 = -\sqrt{2\alpha+1}(a_1 +a_3(2\alpha+1)+(2\alpha+1)^2) \]
Squaring both sides and rearranging gives us a polynomial for $\alpha$.
{\small
\begin{eqnarray}\label{mpa}
p(\alpha)& = & 32 \alpha ^5+16 \left(-a_4^2+2 a_3+5\right) \alpha ^4+8 \left(a_3^2+8 a_3-4 a_4^2+2 a_1+2 a_2 a_4+10\right) \alpha ^3 \\ && \nonumber +4 \left(-a_2^2+2 a_4 a_2+3 a_3^2-6 a_4^2+12 a_3+2 a_1 \left(a_3+3\right)-2 a_0 a_4+10\right) \alpha ^2\\ && +\left(2 \left(a_1+a_3+1\right){}^2+4 \left(a_3+2\right) \left(a_1+a_3+1\right)+4 \left(a_0+a_2+a_4\right) \left(a_2-2 a_4\right)\right) \alpha  \nonumber \\ && \nonumber+\left(\left(a_1+a_3+1\right){}^2-\left(a_0+a_2+a_4\right){}^2\right)
\end{eqnarray}}
There is a unique monic polynomial for the algebraic integer $\alpha$ over $\IQ$ and so all these coefficients are divisible by $32$.
We reduce these coefficients modulo $2$ to determine the parity necessary.
\begin{itemize}
\item \underline{$a_4$ is odd},  $a_4=2k_4+1$.
\item $a_3^2+8 a_3-4 a_4^2+2 a_1+2 a_2 a_4+10$ is divisible by $4$,  so \underline{$a_3$ is even}, $a_3=2k_3$.  
Then $1+ a_1- a_2 a_4$ is even.   Hence  $ a_1+ a_2+1$ is odd.
\item $-a_2^2+2 a_4 a_2+3 a_3^2-6 a_4^2+12 a_3+2 a_1 \left(a_3+3\right)-2 a_0 a_4+10$ is divisible by $8$. So \underline{$a_2$ is even}, $a_2=2k_2$ and \underline{$a_1$ is odd},  $a_1=2k_1+1$. Then
\[-4 a_0 k_4-2 a_0-4 k_2^2+4 k_2+4 k_3^2+4 k_1+4 k_3+2\]
is divisible by $8$,  so \underline{$a_0$ is odd},  $a_0=2k_0+1$.  Then 
\begin{equation}\label{*1} k_0 +  k_1+  k_4 \end{equation}
is even.
\item Next,  making the substitutions identified above,  we have
\begin{equation}\label{*2}  k_1^2+2 k_2^2+3 k_3^2-2 k_0+2 k_0 k_2+2 k_3-2 k_2 k_4-2 k_4+1 \end{equation}
divisible by $4$. Thus $k_1^2 +3 k_3^2$ is odd and so \underline{$k_1+k_3$ is odd}.  
\item Next,  expanding out the constant term gives 
\[ k_0+k_1+k_2+k_3+k_4 \]
is even. Thus $k_0+k_2+k_4$ is odd.  Hence (\ref{*1}) gives $k_1+k_2$ odd.
\end{itemize}
We now record this information.
\begin{align} 
& \label{*5} a_0=2k_0+1,  \quad a_1=2k_1+1,\quad  a_2=2k_2,   \quad a_3=2k_3,\quad a_4=2k_4+1.  \\
& \label{*6} k_1+k_2 \mbox{ is odd},   \quad   k_1+k_3 \mbox{ is odd},\quad  k_0+k_2+k_4 \mbox{ is odd}.
\end{align} 
Writing $k_2=2m_2+1-k_1, \; k_3=2m_3+1-k_1, \; k_4=2m_4-k_0-k_1$ now yields the polynomial
\begin{align*}
& 16  \Big((2-k_1+m_2+m_3+m_4) (k_1-m_2+m_3-m_4) \Big)\\
&+32 \Big(1-3 k_1^2+3 m_3 (2+m_3)-6 m_4+2 k_0 (1-k_1+m_2+m_4)\\ & \quad      +k_1 (4+5 m_2-m_3+7 m_4)-2 (m_2^2+2 m_4^2+m_2 (2+3 m_4)) \Big) \alpha \\
&+32  \Big(4-2 k_0^2-6 k_1^2-5 m_2-2 m_2^2+13 m_3+6 m_3^2-(13+12 m_2) m_4-12 m_4^2  \\ & \quad   +k_0 (6-8 k_1+6 m_2+10 m_4)+k_1 (5+8 m_2-4 m_3+18 m_4) \Big) \alpha ^2
\\&+32 \Big(-4 k_0^2-5 k_1^2+2 k_0 (3-5 k_1+2 m_2+8 m_4)+k_1 (2+4 m_2-4 m_3+20 m_4)   \\ & \quad   -2 (-3+m_2-2 m_3 (3+m_3)+6 m_4+4 m_2 m_4+8 m_4^2)\Big) \alpha ^3
\\ &-64 \Big(k_0^2+k_0 (-1+2 k_1-4 m_4)+(k_1-2 m_4){}^2-2 (1+m_3-m_4)\Big) \alpha ^4+32 \alpha ^5
\end{align*}
Only the constant term 
\[ (2-k_1+m_2+m_3+m_4) (k_1-m_2+m_3-m_4) =(1+m_3)^2 -(1+m_2+m_4-k_1)^2 \]
is now an issue. As it is a difference of squares it is not congruent to $2$ mod $4$.  It also must be even,  thus the constant term in $q_\alpha(z)$ is an even integer.

Now (\ref{*5}) and (\ref{*6})  give necessary and sufficient conditions for the factorisation condition to hold.  
\medskip 

With a little simplification, putting together what is above,  we now  have the following information.

\begin{eqnarray*}
a_0=2k_0+1,\quad a_1=2k_1+1,\quad a_2=4m_2+2-2k_1, \\a_3=4m_3+2-2k_1,\quad a_4=4m_4-2k_0-2k_1+1.   
\end{eqnarray*}

There are another couple of quick tests we can utilise to cut down the search space.  The polynomial $p(z)$ has odd degree,  irreducible and all its real roots in the interval $[-\sqrt{3},\sqrt{3}]$. Thus
\begin{eqnarray*}
0 &< & p(\sqrt{3}) = \sqrt{3} k_1+3 \sqrt{3} k_3+k_0+3 k_2+9 k_4+5 \sqrt{3}+5.\\
0 &> & p(-\sqrt{3}) =-\sqrt{3} k_1-3 \sqrt{3} k_3+k_0+3 k_2+9 k_4-5 \sqrt{3}+5.
\end{eqnarray*}
Subtract the second from the first and substituting,  we find
\begin{equation}
k_1+3k_3+5=k_1+3(2m_3+1-k_1)+5>0
\end{equation}
and so we can implement $m_3\geq (k_1-4)/3$ in the search we  set up.  
\subsection{The search.} At this point we can begin our search.  Let $\lambda,\bar\lambda, r,s,t$ be the roots of the minimal polynomial.  We can choose $\lambda$ to have non-negative real part.  Without any great attempt to get the best possible bounds we obtained 

\medskip

 {\small
\begin{tabular}{|l|c|c|} \hline 
$a_0=|\lambda|^2 rst=2k_0+1$&  $-33 \leq k_0\leq 32 $ \\ \hline
$a_1=|\lambda|^2(rs+rt+st)=2k_1+1$&  $-66 \leq k_1\leq 74$ \\ \hline
$a_2=|\lambda|^2(r+s+t)+(\lambda+\bar\lambda)(rs+rt+st)+rst =2k_2$&  $-32 \leq k_2\leq 67$ \\ \hline
$a_3=|\lambda|^2+rs+rt+st+(\lambda+\bar\lambda)(r+s+t) =2k_3$&  $-21\leq k_3\leq 29$ \\ \hline
$a_4= \lambda+\bar\lambda+r+s+t =2k_4+1 $&  $-3\leq k_4\leq 5$ \\ \hline
\end{tabular}
}

\medskip

With the parity checks and the simple tests suggested this loop has about  $84$ million candidates.  Many of these candidates will have $4$ complex roots.  

We remark that if  there is exactly one complex conjugate pair of roots for the irreducible  minimal polynomial for $\lambda$ of degree $5$, then the same can be said of $\alpha$ if it is complex, since $\alpha\in \IQ(\lambda)$, so  $ \IQ(\alpha)\subset \IQ(\lambda)$ and since the only proper subfields of a field with one complex place are totally real $ \IQ(\alpha)= \IQ(\lambda)$.  However $\alpha$ cannot be real as then
\[ 5=[\IQ:\IQ(\lambda)]=  [\IQ(\alpha):\IQ(\lambda)][\IQ:\IQ(\alpha)] \]
so $[\IQ:\IQ(\alpha)]=1$, $\alpha\in \IZ$ and $\lambda$ is an integer in a quadratic field.

\medskip

We could implement a further test of the discriminant in these searches,  but polynomial root finding algorithms are quick enough that simply computed the roots and checked that $\lambda$ was a viable candidate,  that is all the real roots are in $[-\sqrt{3},\sqrt{3}]$.  This gave us a list with about 5100 candidates.  The first we found was $x^5-x^4+6x^3+22x^2-23x-51$ with real roots $-1.674, -1.490,1.62$ and $\lambda=1.269+i3.308$.  However,  as we did not implement {\em all} the parity checks in the search to keep the structure simple,
\[ \gamma=\frac{1}{2}(\lambda^2-3) = -6.168+i4.20 \]
has a minimal polynomial of degree $5$ with real roots in $[-\frac{1}{2},1]$,  but it is not monic.  

We then searched through these polynomials to find those minimal polynomials for $\gamma=(\lambda^2-3)/2$ which were monic,  had $-4<\Re e(\gamma)$ and $|\Im(\gamma)|<4$.  The values of $\gamma$ we found are illustrated below.

  \scalebox{0.7}{\includegraphics[viewport=50 430 770 850]{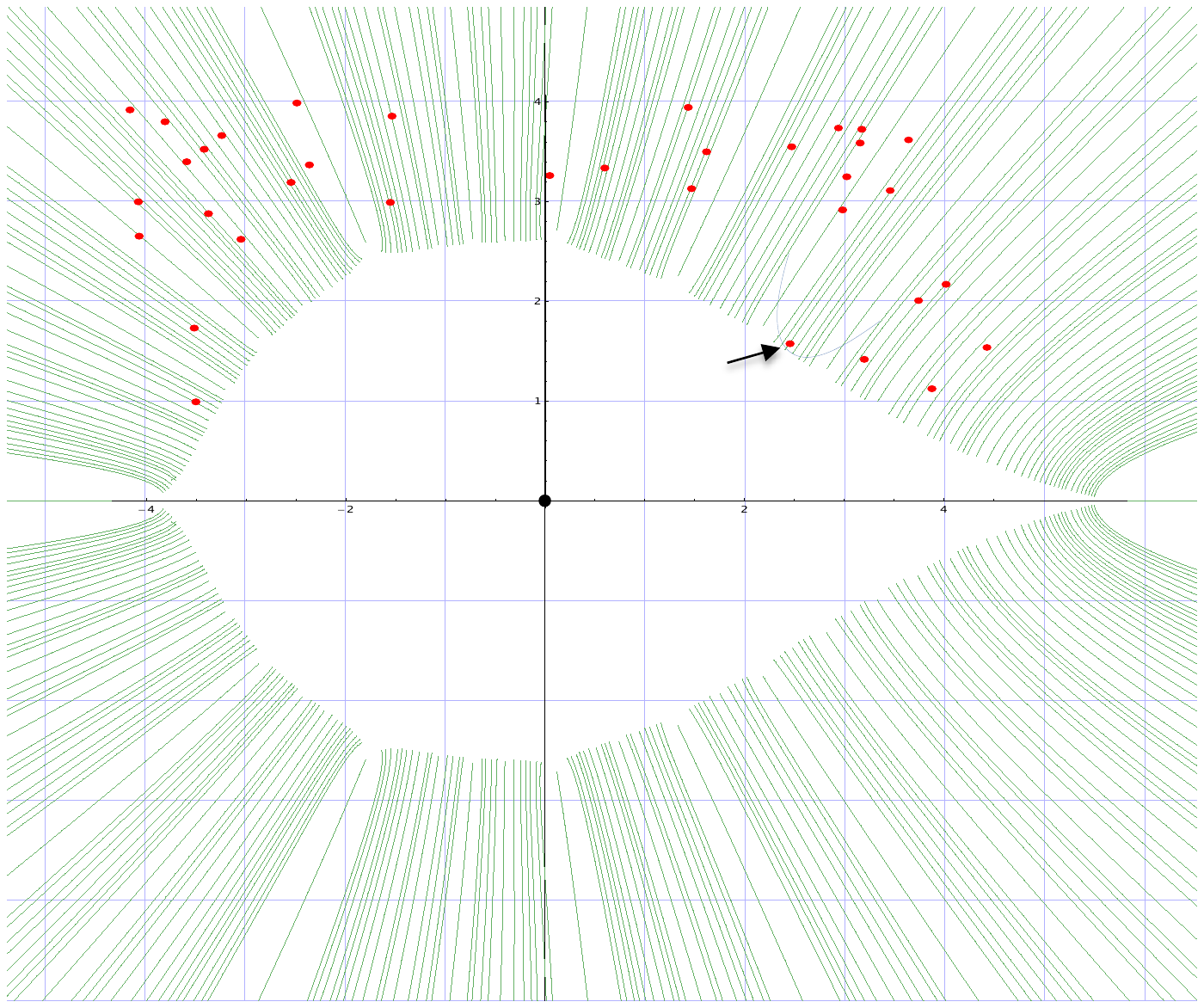}}\\
  
  \noindent{\bf Figure 10.} {\em The degree $5$ candidates with $11/13$ pleating ray neighbourhood illustrated.}

The most computationally challenging polynomial we found was
\[ 1+7 x+12 x^2-3 x^4+x^5  \]
with real roots $-1.29375, -0.378269, -0.240317$ and complex roots $2.45617 \pm i 1.57165$.  This was captured by the $11/13$ pleating ray neighbourhood illustrated above in Figure 10.

\subsection{Additional searches} As the degree grows the search space is growing rapidly and the fact that all the real roots of $\alpha$ are in $[-\frac{1}{2},1]$ means symmetric functions grow slowly and we can obtain better bounds for our searches (although this does not happen in degree $5$).  We now outline how this goes in this case and present data.

\medskip

Let $q_\alpha(z)=z^5+b_4z^4+b_3z^3+b_2z^2+b_1z^1+b_0$ be the minimal polynomial for $\alpha$.
Equations (\ref{*5}) and (\ref{*6}) together with the difference of squares condition we found yield the following information on the coefficients $b_i$.
\begin{itemize}
\item $b_4$ is even,  $b_4=2n_4$ and $n_4$ has the same parity as  $m_2+m_3+m_4$.
\item $b_3$ has the same parity as  $m_2+m_3+m_4$.
\item $b_2$ is even,  $b_2=2n_2$ 
\item $b_1$ has the same parity as  $1+m_3$
\item $b_0$ is even,  $b_0=2n_0$ and $n_0$ has the same parity as  $1+m_3$.
\item $b_3$ has the same parity as $n_4$
\item $b_1$ has the same parity as $n_0$
\end{itemize}
We now go about finding bounds on these quantities.
\subsection{Loop bounds.}
We recall $-8\leq \alpha+\bar\alpha \leq 6+4\sqrt{2}$,  $|\alpha|\leq 3+2\sqrt{2}$, and $r,s,t\in [-\frac{1}{2},1]$.
The following information which follows from calculus.  We set $x=|\alpha|^2$ and $2x=\alpha+\bar\alpha$,  $-4<x<(3+2\sqrt{2})$ and calculate gradients to find extrema are on the boundary and then minimise. So for instance, $-0.737 < rs+rt+st <3$. 
\begin{eqnarray*}
-33.98 &< b_0 =&- |\alpha|^2 rst  <   |\alpha|^2/2 < 16.99 \\
 -25.3 &< b_1=&|\alpha|^2(rs+rt+st)+(\alpha+\bar\alpha)(r+s+t) < 136.9 \\
  -108.2 &< b_2=&-|\alpha|^2(r+s+t)-(\alpha+\bar\alpha)(rs+rt+st) < 55.3 \\
-24 &< b_3=&(\alpha+\bar\alpha)(r+s+t) < 34.9 \\
-14.7 &< b_4=&-(\alpha+\bar\alpha + r+s+t) < 9.5 
\end{eqnarray*}
Thus
\begin{eqnarray*}
-16 \leq n_0 \leq 8 &
 -25\leq b_1 \leq 136 &
  -54 \leq n_2\leq  26 \\
-23  \leq b_3 \leq   34 &&
-7 \leq n_4 \leq 4 
\end{eqnarray*}
We actually ran this search (which took significantly longer) to confirm our results.

\section{The case of degree $6$ \& $7$}
\subsection{Degree $6$ parity considerations}

As before we obtain another expression for the minimal polynomial for $\alpha$ with the coefficients coming from the minimal polynomial for $\lambda$.

\begin{align}
&\left(-1-a_0+a_1-a_2+a_3-a_4+a_5\right) \left(1+a_0+a_1+a_2+a_3+a_4+a_5\right) \nonumber \\
&+2 \left(-2 a_0 \left(3+a_2+2 a_4\right)-2 \left(1+a_2+a_4\right) \left(3+a_2+2 a_4\right)+\left(a_1+a_3+a_5\right) \left(a_1+3 a_3+5 a_5\right)\right) \alpha \nonumber \\ &+\left(8 a_1 a_3+12 a_3^2-8 \left(3+a_4\right) \left(1+a_0+a_2+a_4\right)-4 \left(3+a_2+2 a_4\right){}^2+24 a_1 a_5+48 a_3 a_5+40 a_5^2\right) \alpha ^2\nonumber \\ &+8 \left(a_3^2-2 \left(a_0+a_2 \left(4+a_4\right)+2 \left(5+a_4 \left(5+a_4\right)\right)\right)+2 a_1 a_5+8 a_3 a_5+10 a_5^2\right) \alpha ^3\nonumber \\&-16 \left(15+2 a_2+a_4 \left(10+a_4\right)-2 a_3 a_5-5 a_5^2\right) \alpha ^4\nonumber\\ &+32 \left(-6-2 a_4+a_5^2\right) \alpha ^5-64 \alpha ^6 \label{afor}
\end{align}
The necessary parity conditions in degree 6 are
\begin{itemize}
\item $a_0=2m_0+1$
\item $a_1=2m_1$
\item $a_2=4m_2+3$
\item $a_3=4m_3$
\item $a_4=8m_4-4m_2-2m_0+3$
\item $a_5=8m_5-4m_3-2m_1 $
\end{itemize}
With the addition of another parity condition $m_4+m_5$ odd,  these are sufficient as well.
Flammang and Rhin give a method for bounding values of symmetric functions and coefficients using various identities.  We adopt those in degree $7$ and $8$.  However here we simply adopt the following strategy.  We replace $|\lambda|$ by $x\in [0,\sqrt{7+\sqrt{2}}]$ and $\lambda+\bar\lambda$ by $2x$.  Calculate the gradient of the function of $5$ variables to see there are no interior extrema in the region $x\in[0,\sqrt{7+4\sqrt{2}}]$, $r_i\in [-\sqrt{3},\sqrt{3}]$.  Then symmetry allows us to move the $r_i$ in an increasing or decreasing manner to the value $\pm \sqrt{3}$.  We could improve these estimates by realising that $r_i\approx \pm \sqrt{3}$ implies the corresponding real embedding gives $\sigma_i(\alpha)\approx 1$ and we have seen that the real roots of the minimal polynomial for $\alpha$ cannot pile up near $1$.
\[ -113 \leq a_0 = |\lambda|^2 r_1 r_2 r_3 r_4 \leq 113 \]
\[   a_1 =- |\lambda|^2 (r_1 r_2 r_3 +r_1r_2r_4+r_2r_3r_4)-(\lambda+\bar\lambda)(r_1r_2r_3r_4)   \]
\[ -261 \leq a_1\leq 129 \]
\[ a_2 = |\lambda|^2 (r_1 r_2 +r_1r_3 +r_1r_4+r_2r_3+r_2r_4+r_3r_4)+(\lambda+\bar\lambda)(r_1 r_2 r_3 +r_1r_2r_4+r_2r_3r_4) \leq   \]
\[ -188 \leq a_2   \leq  338 \]
\[ a_3 = -|\lambda|^2 (r_1 +r_2+r_3 +r_4)-(\lambda+\bar\lambda)(r_1 r_2 +r_1r_3 +r_1r_4+r_2r_3+r_2r_4+r_3r_4)    \]
\[ -215 \leq a_3   \leq 43  \]
\[  a_4=|\lambda|^2+(\lambda+\bar\lambda)(r_1+r_2+r_3+r_4)+(r_1r_2+r_1r_3+r_1r_4+ r_2r_3+ r_2r_4+r_3r_4) \]
\[-11 \leq a_4\leq 79 \]
\[ -14\leq  a_5=-\lambda-\bar\lambda-r_1-r_2-r_3-r_4\leq 6 \]
At first sight this loop seems to allow $20\times 91 \times 258\times 526\times 390\times 226 \approx 21\times 10^{12}$ possibilities.  However the parity considerations reduce this to only $5\times 10^9$ possibilities.  For instance,  $a_5$ lies in an interval of width $20$ and so after having determined $n_1,\ldots,n_4$ the value $m_5$ must lie in an interval of length $20/8$ and so admits at most three values. There is more than can be said here. If we review the formula (\ref{afor}) we see that the coefficient for $\alpha^5$ there is $32 \left(-6-2 a_4+a_5^2\right)$,  which must be divisible by $64$ and is since $a_5$ is even. In terms of the real roots $s_i$ of the minimal polynomial for $\alpha$ we find that there is an integer $m$ so that
\[ \alpha+\bar\alpha+s_1+s_2+s_3+s_4 = (-3- a_4+\frac{1}{2} a_5^2) m \]
In the loop described above we can calculate $m_1,m_3$ and hence $m_5$ before we calculate $m_4$, and hence $a_4$. Now
\[   -9\leq    \alpha+\bar\alpha+s_1+s_2+s_3+s_4 \leq 17 \]
Hence,  given $a_5$, 
\[   -9\leq  \big(-3- a_4+\frac{1}{2} a_5^2\big) m  \leq 17 \]
If $a_4$ lies in an interval of width $\frac{26}{|m|}$ and hence $m_4$ lies in an interval of width $\frac{26}{8|m|}$.  If $|m|=1$,  then that gives at most $4$ possibilities for $a_4$.  Decreasing the search space by a factor of at least $3$. If $|m|=2$, then only three possibilities and as soon as $m\geq 4$ there is at most one possibilty,  with 
\[   2 \leq  - 4m_4+2m_2+m_0+\frac{1}{4} a_5^2   \leq 10 \]
Other observations, coming from the additional parity condition, are that  the constant term coefficient in the minimal polynomial for $\alpha$ is  
\[ \frac{1}{64} \left(-1-a_0+a_1-a_2+a_3-a_4+a_5\right) \left(1+a_0+a_1+a_2+a_3+a_4+a_5\right) \]
\[ = \frac{1}{4} (3 + 4 m_4 - 4 m_5) (1 + m_4 + m_5)     \]
 so we must have $m_4+m_5$ odd,   and this halves the search space again.  These, and a few other simple observations, reduce the search space by another order of magnitude and additionally give simple tests to eliminate candidates before we use a root finding algorithm.

   \scalebox{0.65}{\includegraphics[viewport=-10 350 670 850]{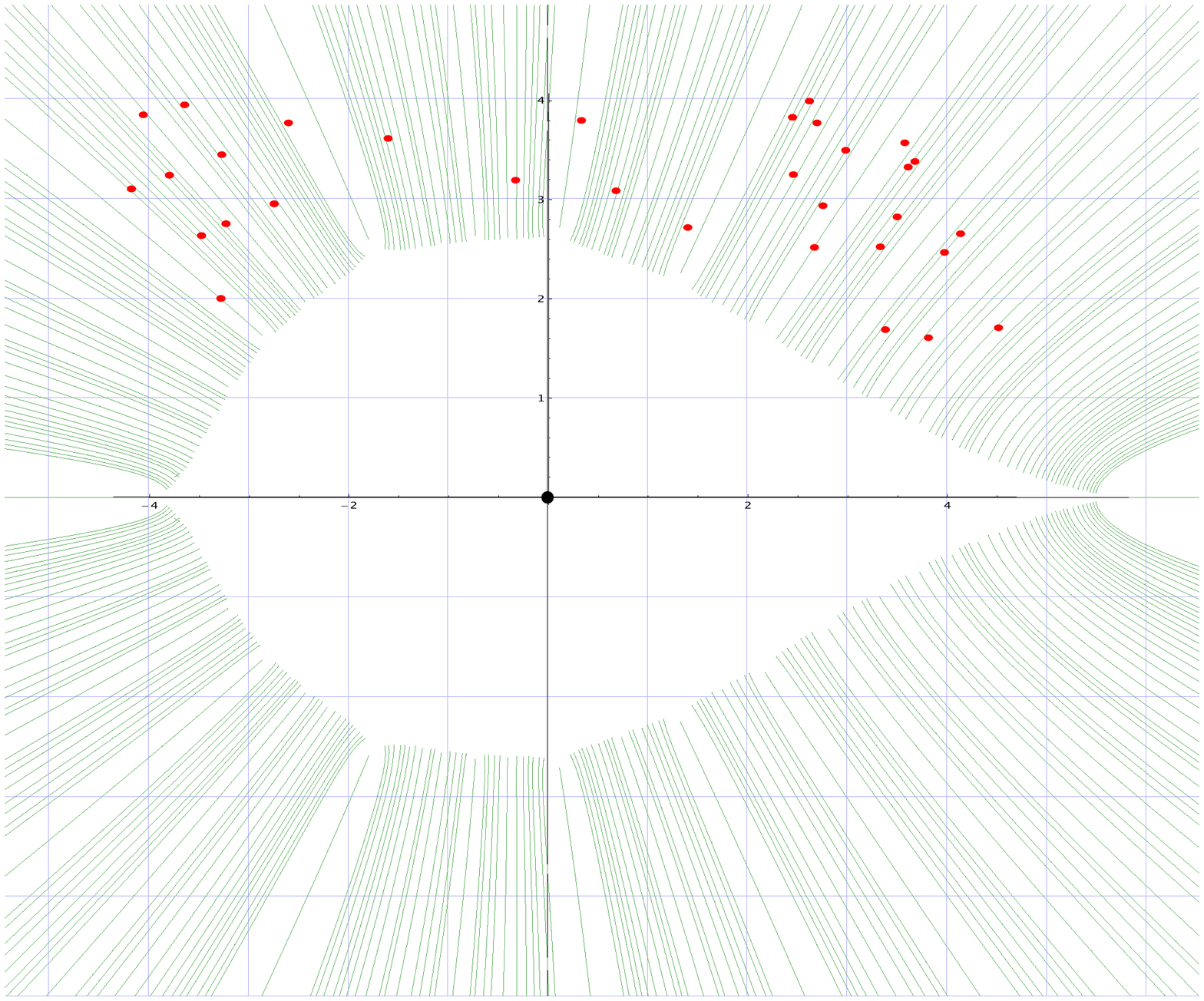}}\\
   
     \noindent{\bf Figure 11.} {\em The degree $6$ candidates satisfying all arithmetic criteria.}
     
     \medskip
     
We ran this search to find $32$ points,  illustrated here, and that we should consider further.  These were not in the required region.

\subsection{Degree $7$ parity considerations}
\begin{align*}
&128 \alpha ^7+64 \alpha ^6 \left(2 a_5-a_6 ^2+7\right)+32 \alpha ^5 (2 a_3-2 a_6  (a_4+3 a_6 )+a_5 (a_5+12)+21)
\\
&+16 \alpha ^4 \left(2 a_1-2 a_2 a_6 +2 a_3 (a_5+5)-a_4^2-10 a_4 a_6 +5 a_5 (a_5+6)-15 a_6 ^2+35\right)\\ 
&+8 \alpha ^3 \left(-2 a_0a_6 +2 a_1 (a_5+4)-2 a_2 (a_4+4 a_6 )+a_3^2+4 a_3 (2 a_5+5) \right. \\
& \left. \quad -4 \left(a_4^2+5 a_4 a_6 +5 a_6 ^2\right)+10 a_5 (a_5+4)+35\right)\\
&+4 \alpha ^2 \left(-2 (a_4+3 a_6 ) (a_0+a_2+a_4+a_6 )+2 a_1 (a_3+3 a_5+6)-(a_2+2 a_4+3 a_6 )^2\right.\\
& \left. \quad+3 a_3^2+4 a_3 (3 a_5+5)+10 a_5^2+30 a_5+21\right)\\
&+\alpha  \left(-4 (a_2+2 a_4+3 a_6 ) (a_0+a_2+a_4+a_6 )+2 (a_1+a_3+a_5+1)^2 \right. \\
& \left. \quad +4 (a_3+2 a_5+3) (a_1+a_3+a_5+1)\right) \\ & +\left((a_1+a_3+a_5+1)^2-(a_0+a_2+a_4+a_6 )^2\right)
\end{align*}
We continue to reduce these coefficients mod $2$ making substitutions as appropriate. Then \underline{$a_6$ is odd},  $a_6=2k_6+1$.  Then \underline{$a_5$} is odd, $a_5=2k_5+1$. Hence $a_3+a_4$ is even,  and subsequently \underline{$a_4$, hence $a_3$} are odd,  $a_4=2k_4+1$, $a_3=2k_3+1$.  It follows \underline{$a_0$ is odd}, $a_0=2k_0+1$. Subsequently \underline{$a_1$ and $a_2$ are odd}, $a_2=2k_2+1$,$a_1=2k_1+1$.  Expanding out and looking at the coefficients in terms of $k_i$ we find
\begin{itemize}
\item $k_1+k_2+k_4$ is even.
\item $k_0+k_1+k_2+k_4+k_5+k_6$ is odd.
\item $k_2+k_3+k_6$ is even.
\item $k_1+k_3+k_5$ is even.
\item $k_0+k_1+k_2+k_3+k_4+k_5+k_6$ is even.
\end{itemize}
Hence $k_3$ is odd.
\begin{itemize}
\item $k_0+k_1+k_2$ is odd. $k_2=2m_2+1-k_0-k_1$
\item $k_0+k_4$ is odd. $k_4=2m_4+1-k_0$
\item $k_2+k_6$ is odd. $k_6=2\tilde{m_6}+1-k_2=2m_6+k_0+k_1$
\item $k_1+k_5$ is odd. $k_5=2m_5+1-k_1$
\end{itemize}
At this point the coefficients of $\alpha^4,\alpha^5$ and $\alpha^6$ are divisible by $128$.  We continue.  Looking at the coeffficient of $\alpha^2$ gives $k_3=2m_3+1$ odd and deals with the coefficient of $\alpha^3$.  Next $m_3+m_5$ is even,  as is $m_2+m_4+m_6$.  Subsequently $m_4+m_5$,  is even, and $m_3+m_4$,  is even.  Then we have
\begin{itemize}
\item  $k_2=2m_2+1-k_0-k_1$
\item  $k_3=2m_3+1$
\item $k_4=2(2n_4-m_3)+1-k_0$
\item $k_6=2(2n_6-(2n_4-m_3)-m_2)+k_0+k_1$
\item $k_5=2(2n_5-m_3)+1-k_1$
\end{itemize}
Relabeling we now have an inductive structure to determine the coefficients.
\begin{itemize}
\item  $a_0=2n_0+1$
\item  $a_1=2n_1+1$
\item  $a_2=4n_2+3-2(n_0+n_1) $    
\item  $a_3=4n_3+3$
\item $a_4=8n_4-4n_3-2n_0+3$
\item $a_5=8n_5-4n_3-2n_1+3$
\item $a_6=8n_6- 8n_4-4n_3-4n_2+2n_0+2n_1+1  $ 
\end{itemize}

\subsection{The case of degree $8$}

Let
\[ \lambda ^8+a_7\lambda ^7+a_6\lambda ^6+a_5\lambda ^5+a_4\lambda ^4+a_3 \lambda ^3+a_2\lambda ^2+a_1\lambda +a_0=0 \]
be the minimal polynomial for $\lambda$.  Let
\[ \alpha ^8+b_7\alpha ^7+b_6\alpha ^6+b_5\alpha ^5+b_4\alpha ^4+b_3 \alpha ^3+b_2\alpha ^2+b_1\alpha +b_0 = 0 \]
be the minimal polynomial for $\alpha$. \\

\begin{lemma} The disk $\{|\alpha|<4.4960\}$ is excluded.  Hence 
\[ 4.49 \leq |\alpha| \leq 4+2\sqrt{2} \]
and 
\[ 4\leq \Re e(\alpha) \leq 4+2\sqrt{2}  \]
\end{lemma}

  \scalebox{0.75}{\includegraphics[viewport=40 450 470 800]{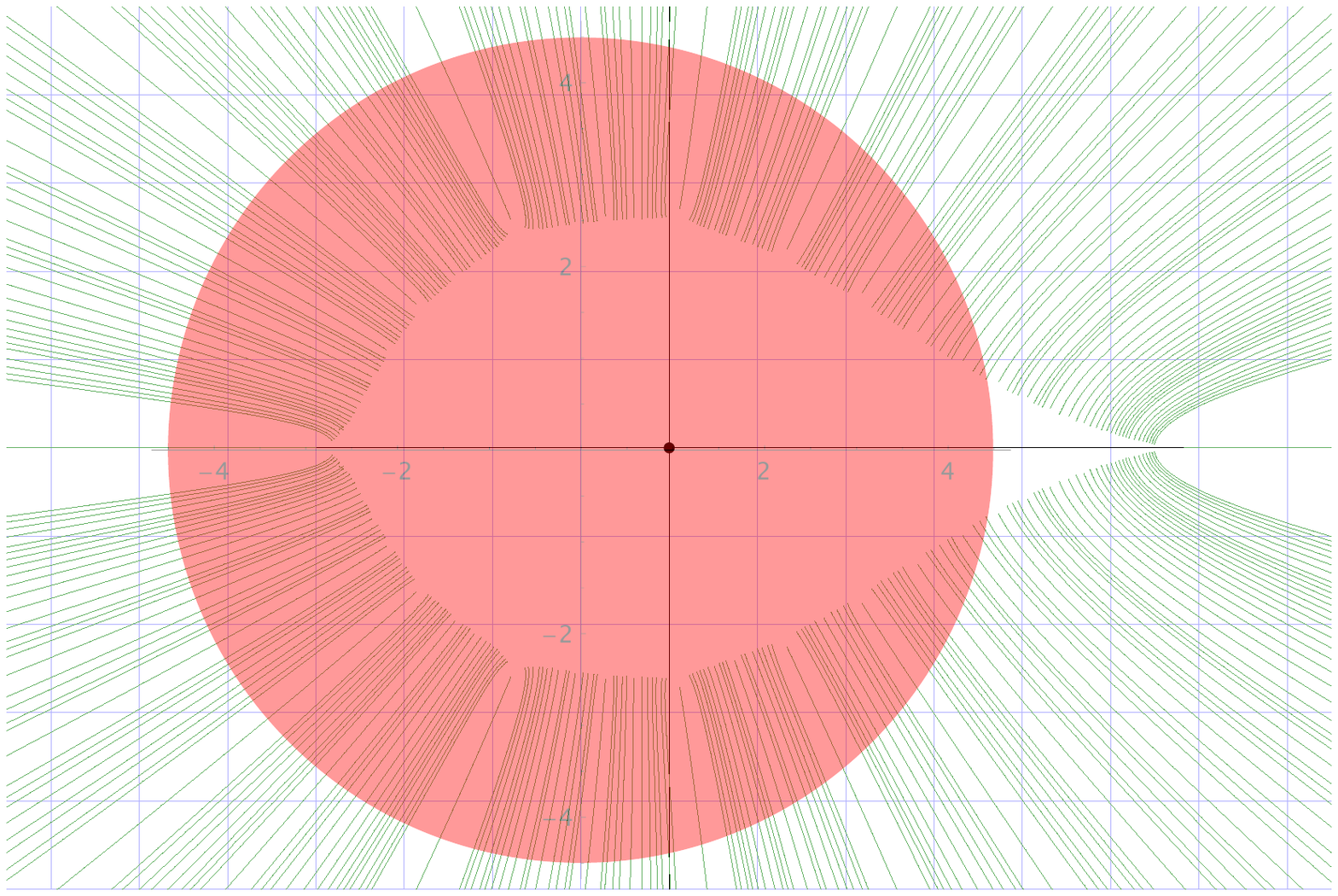}}\\

\noindent {\bf Figure 13.}{\em The excluded disk in the $\alpha$-plane.}

\medskip

\noindent{\bf Proof.}   Let $\alpha_0$ be the root of the equation
\[ \alpha_0^{2 t_0} (\alpha_0 - 1)^2 (0.0770561)^6 = 1 \]
Then (noting $\gamma=\alpha-1$) and $|\alpha|< \alpha_0$ Theorem \ref{a0bound} gives
\[ |p(0)|^{t_0}|p(1)|\leq |\alpha|^{2 t_0}|\alpha - 1|^2 (0.0770561)^6 <  1 \]
and this implies that $p(0)=0$ or $p(1)=0$ which is impossible. \hfill $\Box$

\medskip

Thus the possible values of $\alpha$ that we seek are confined to a rather small region. It is rather unsurprising then that there are no degree $8$ cases.  Against this however,  we do know there are infinitely many lattices in this region.  

We next recall Lemma \ref{*lemma} that implies $|p(0)|=|b_0|\leq 2$ in this case.  We calculate as before to find $2^8$ times the minimal polynomial for $\alpha$ to be

{\small \begin{align*}
&(-1 - a_0 + a_1 - a_2 + a_3 - a_4 + a_5 - a_6 + a_7) (1 + a_0 + a_1 + a_2 + a_3 + 
    a_4 + a_5 + a_6 + a_7) \\ &+ 
 2 (-2 a_0 (4 + a_2 + 2 a_4 + 3 a_6) - 
    2 (1 + a_2 + a_4 + a_6) (4 + a_2 + 2 a_4 + 3 a_6)  \\ &  + (a_1 + a_3 + a_5 + 
       a_7) (a_1 + 3 a_3 + 5 a_5 + 7 a_7)) \alpha \\ & + 
 4 (2 a_1 a_3 + 3 a_3^2 + 6 a_1 a_5 + 12 a_3 a_5 + 10 a_5^2 - 
    2 (1 + a_0 + a_2 + a_4 + a_6) (6 + a_4 + 3 a_6) \\& - (4 + a_2 + 2 a_4 + 
      3 a_6)^2 + 12 a_1 a_7 + 20 a_3 a_7 + 30 a_5 a_7 + 
    21 a_7^2) \alpha^2 \\ & + 
 8 (a_3^2 + 2 a_1 a_5 + 8 a_3 a_5 + 10 a_5^2 - 
    2 (4 + a_6) (1 + a_0 + a_2 + a_4 + a_6)\\& - 
    2 (6 + a_4 + 3 a_6) (4 + a_2 + 2 a_4 + 3 a_6) + 8 a_1 a_7 + 20 a_3 a_7 + 
    40 a_5 a_7 + 35 a_7^2) \alpha^3 \\ &+ (-32 a_0 - 
    16 (70 + a_4^2 - a_5 (2 a_3 + 5 a_5) + 10 a_4 (3 + a_6) + 
       2 a_2 (5 + a_6) \\ & + 5 a_6 (14 + 3 a_6)) + 
    32 (a_1 + 5 (a_3 + 3 a_5)) a_7 + 560 a_7^2) \alpha^4 \\& + 
 32 (a_5^2 - 2 (28 + a_2 + a_4 (6 + a_6) + 3 a_6 (7 + a_6)) + 2 a_3 a_7 + 
    12 a_5 a_7 + 21 a_7^2) \alpha^5 \\ & - 
 64 (28 + 2 a_4 + a_6 (14 + a_6) - 2 a_5 a_7 - 7 a_7^2) \alpha^6 \\ &+ 
 128 (-8 - 2 a_6 + a_7^2) \alpha^7  - 256 \alpha^8
\end{align*}}
 
We then inductive reduce modulo $2$ to gain the parity of the coefficients,  make the appropriate substitutions and repeat until all the remaining coefficients are a multiple of $2^8$ (we cary out this process more carefully in the next subsection for degree $9$).  After this we find that there are integers $k_0,k_1,\ldots,k_8$ so that
\begin{eqnarray*}
a_0 & = &  1 + 2 k_0 \\
a_1 & = &  -2+4 k_1 \\
a_2 & = &   -2 k_0+4 k_2 \\
a_3 & = &   2-4 k_1+4 k_3 \\
a_4 & = &  -2 \left(1+k_0+4 k_1-4 k_4\right) \\
a_5 & = & 2+4 k_1-8 k_4+8 k_5  \\
a_6 & = &  -2 \left(7 k_0-4 k_1+2 k_2+4k_4-8 k_6\right) \\
a_7 & = &   -2 \left(1+2 k_1+2 k_3+4 k_4+4 k_5-8 k_7\right)
\end{eqnarray*}

This gives us the minimal polynomial for $\alpha$ as 
\[ x^8+b_7x^7+b_6x^6+b_5x^5+b_4x^4+b_3x^3+b_2x^2+b_1 x+b_0 \] 
and where
\begin{eqnarray*} 
b_0 &= &\big(k_0-k_6\big){}^2-\big(k_4-k_7\big){}^2,\\
b_1 & = & 2\big(\big(6 k_0-k_1+k_2+k_4-6 k_6\big) \big(k_0-k_6\big)\\ && \quad -\big(k_1+k_3+5 k_4+k_5-6 k_7\big) \big(k_4-k_7\big)-\big(k_4-k_7\big){}^2\big),\\
b_2 & =&  2 \big(-1+11 k_0-4 k_1+3 k_2+4 k_4-12 k_6\big) \big(k_0-k_6\big)\\ && \quad -\big(k_1+k_3+5 k_4+k_5-6 k_7\big){}^2-2 \big(1+2 k_1+3 k_3+8 k_4+4 k_5-12 k_7\big) \big(k_4-k_7\big)\\ && \quad -4 \big(k_1+k_3+5 k_4+k_5-6 k_7\big) \big(k_4-k_7\big)\\ && \quad +\big(6 k_0-k_1+k_2+k_4-6 k_6\big){}^2,\\
b_3 & = & =2\big(73 k_0^2+k_1^2-k_2+3 k_2^2-k_3\\ && \quad-4 k_3^2-9 k_4+7 k_2 k_4-41 k_3 k_4-81 k_4^2-k_5\\ && \quad -9 k_3 k_5-50 k_4 k_5-5 k_5^2+k_0 \big(-8-39 k_1+31 k_2+39 k_4-153 k_6\big)+8 k_6\\ && \quad-32 k_2 k_6-40 k_4 k_6+80 k_6^2-k_1 \big(7 k_2+7 k_3+42 k_4+8 k_5-40 k_6-42 k_7\big)+9 k_7\\ && \quad+50 k_3 k_7+220 k_4 k_7+60 k_5 k_7-140 k_7^2\big),\\ 
b_4& = & -2 k_0+2 \big(-2+7 k_0-4 k_1+2 k_2+4 k_4-8 k_6\big) \big(6 k_0-k_1+k_2+k_4-6 k_6\big)+2 k_6\\ && \quad+\big(1-11 k_0+4 k_1-3 k_2-4 k_4+12 k_6\big){}^2\\ && \quad-4 \big(\big(1+2 k_1+3 k_3+8 k_4+4 k_5-12 k_7\big) \big(k_1+k_3+5 k_4+k_5-6 k_7\big)\\ && \quad+\big(1+2 k_1+2 k_3+4 k_4+4 k_5-8 k_7\big) \big(k_4-k_7\big)\big)\\ && \quad-\big(1+2 k_1+3 k_3+8 k_4+4 k_5-12 k_7\big){}^2\\ && \quad-2 \big(1+2 k_1+2 k_3+4 k_4+4 k_5-8 k_7\big) \big(k_1+k_3+5 k_4+k_5-6 k_7\big)
,\\ b_5 & = & 2 \big(\big(-1+11 k_0-4 k_1+3 k_2+4 k_4-12 k_6\big) \big(-2+7 k_0-4 k_1+2 k_2+4 k_4-8 k_6\big) \\ && \quad -6 k_0+k_1-k_2-k_4+6 k_6-\big(1+2 k_1+3 k_3+8 k_4+4 k_5-12 k_7\big){}^2\\ && \quad-\big(1+2 k_1+3 k_3+8 k_4+4 k_5-12 k_7\big) \big(1+2 k_1+2 k_3+4 k_4+4 k_5-8 k_7\big)\\ && \quad-2 \big(1+2 k_1+2 k_3+4 k_4+4 k_5-8 k_7\big) \big(k_1+k_3+5 k_4+k_5-6 k_7\big)\big),\\
b_6&=& 2-22 k_0+8 k_1-6 k_2-8 k_4+24 k_6+\big(2-7 k_0+4 k_1-2 k_2-4 k_4+8 k_6\big){}^2\\ && \quad-4 \big(1+2 k_1+3 k_3+8 k_4+4 k_5-12 k_7\big) \big(1+2 k_1+2 k_3+4 k_4+4 k_5-8 k_7\big)\\ && \quad-\big(1+2 k_1+2 k_3+4 k_4+4 k_5-8 k_7\big){}^2,\\
b_7 & = & 2 \big(2-7 k_0+4 k_1-2 k_2-4 k_4+8 k_6-\big(1+2 k_1+2 k_3+4 k_4+4 k_5-8 k_7\big){}^2\big).
\end{eqnarray*}

Notice that $b_0$ is a difference of squares.  As $2$ is not a quadratic residue we see that  $b_0\neq \pm 2$ and hence must be $\pm 1$,  by Lemma \ref{*lemma}.

\begin{lemma}  If $a^2-b^2=\pm 1$,  then 
\[ (a,b)\in \{(0,\pm1),(\pm 1,0) \}\]
\end{lemma}
\noindent{\bf Proof.} $a^2-b^2= (a+b)(a-b)=\pm 1$ so both factors are $\pm 1$.  $a-b=1$ gives $a=1+b$ and $1+2b=\pm 1$ so $b=-1$, and $a=0$,  or $b=0$. 
$a-b=-1$ gives $a=-1+b$ and $-1+2b=\pm 1$ so $b=1$, and $a=0$,  or $b=0$.  \hfill $\Box$

\medskip
 
We can now make the following assertions about the coefficients for the minimal polynomial of $\alpha$.

\begin{lemma}  Let 
\[ \alpha ^8+b_7\alpha ^7+b_6\alpha ^6+b_5\alpha ^5+b_4\alpha ^4+b_3 \alpha ^3+b_2\alpha ^2+b_1\alpha +b_0 = 0 \]
be the minimal polynomial for $\alpha$.  Then
\begin{itemize}
\item $b_0=\pm 1$.
\item $b_7$, $b_5$, $b_3$ and $b_1$ are even.
\item $b_6$ has the same parity as $\frac{1}{2} b_7$.
\item $b_4$ has the opposite parity to $b_6$.
\end{itemize}
\end{lemma}

\subsection{Degree $8$ coefficient bounds.}

We want to use the information in the above lemma to prove that there are no such polynomials,  as in degree $7$.  As earlier,  the coefficients in the search are determined by the root structure.  But this time the search for $\lambda$ is quite a bit larger.  However we can use the excluded disk to improve search bounds as $4.49 \leq |\alpha| \leq 4+2\sqrt{2}$ and $\alpha+\bar\alpha>8.9$ as we recall that all the real roots of the minimal polynomial for $\alpha$ lie in $[-1,1]$. Then a few elementary estimates got the search runtime down to a day.  We found no candidates.

\subsection{The case of degree $9$}

In light of our first calculation,  the following lemma suffices to deal with this case.

\begin{lemma} Let $\alpha$ be a complex root of a degree $9$ polynomial with $\lambda=\sqrt{2\alpha+1}$ an algebraic integer and we have $\IQ(\lambda)=\IQ(\alpha)$.  Then $\alpha$ is not a unit.
\end{lemma} 
The proof is simply a very long calculation.  We repeatedly apply the process above to determine the parity of the coefficients for the minimal polynomial for $\lambda$,
\[\lambda ^4+a_8\lambda ^8+a_7\lambda ^7+a_6\lambda ^6+a_5\lambda ^5+a_4\lambda ^4+a_3 \lambda ^3+a_2\lambda ^2+a_1\lambda +a_0\]
As before we find a multiple of the minimal polynomial for $\alpha$,

\noindent {\tiny
\begin{tabular}{|l|} \hline 
$ \left(1+a_1+a_3+a_5+a_7\right){}^2-\left(a_0+a_2+a_4+a_6+a_8\right){}^2 $ \\ \hline
$ \left(2 \left(1+a_1+a_3+a_5+a_7\right){}^2+4 \left(1+a_1+a_3+a_5+a_7\right) \left(4+a_3+2 a_5+3 a_7\right)\right. $\\ $\left.-4 \left(a_0+a_2+a_4+a_6+a_8\right) \left(a_2+2 a_4+3 a_6+4 a_8\right)\right) \alpha $ \\  \hline
$\left(8 \left(1+a_1+a_3+a_5+a_7\right) \left(6+a_5+3 a_7\right)+8 \left(1+a_1+a_3+a_5+a_7\right) \left(4+a_3+2 a_5+3 a_7\right) \right.$ \\ $\left. +4 \left(4+a_3+2 a_5+3 a_7\right){}^2-4 \left(a_2+2 a_4+3 a_6+4 a_8\right){}^2-8 \left(a_0+a_2+a_4+a_6+a_8\right) \left(a_4+3 a_6+6 a_8\right)\right) \alpha ^2$ \\ \hline
$8 \left(84+a_3^2+70 a_5+10 a_5^2+112 a_7+40 a_5 a_7+35 a_7^2+2 a_1 \left(10+a_5+4 a_7\right)+4 a_3 \left(2 a_5+5 \left(2+a_7\right)\right) \right. $ \\ $\left.-2 \left(a_0+a_2+a_4+a_6+a_8\right) \left(a_6+4 a_8\right)-2 \left(a_2+2 a_4+3 a_6+4 a_8\right) \left(a_4+3 a_6+6 a_8\right)\right) \alpha ^3$ \\ \hline
$16 \left(126-a_4^2+70 a_5+5 a_5^2-2 a_2 a_6-10 a_4 a_6-15 a_6^2+140 a_7+30 a_5 a_7+35 a_7^2+2 a_1 \left(5+a_7\right)\right.$\\ $\left. +2 a_3 \left(15+a_5+5 a_7\right)-2 a_0 a_8-10 a_2 a_8-30 a_4 a_8-70 a_6 a_8-70 a_8^2\right) \alpha ^4$\\ \hline
$32 \left(126+2 a_1+42 a_5+a_5^2-2 a_4 a_6-6 a_6^2+112 a_7+12 a_5 a_7+21 a_7^2+2 a_3 \left(6+a_7\right)-2 a_2 a_8-12 a_4 a_8-42 a_6 a_8-56 a_8^2\right) \alpha ^5$\\ \hline
$64 \left(84+2 a_3-a_6^2+56 a_7+7 a_7^2+2 a_5 \left(7+a_7\right)-2 a_4 a_8-14 a_6 a_8-28 a_8^2\right) \alpha ^6$\\  \hline
$128 \left(36+2 a_5+16 a_7+a_7^2-2 a_6 a_8-8 a_8^2\right) \alpha ^7$\\  \hline
$256 \left(9+2 a_7-a_8^2\right) \alpha ^8 + 512\alpha ^9$ \\ \hline  
\end{tabular}
}

\bigskip
We see
\begin{itemize} 
\item $ 9+2 a_7-a_8^2$ even, so \underline{$a_8$ is odd.}
\item $2 a_5+a_7^2-2 a_6$ divisible by $4$, so \underline{$a_7$ is even} and \underline{$a_5-a_6$ is even.}
\item $2 a_3-2 a_4+6 a_5-6 a_6-a_6^2+4 a_5 k_7+4 k_7^2-4 a_4 k_8-4 a_6 k_8$ divisible by $8$,  so \underline{$a_5$ and $a_6$ are even}. Further,  $a_3-a_4+2 k_5-2 k_6-2 k_6^2+2 k_7^2-2 a_4 k_8$ is divisible by $4$.
\item $3+a_1-a_2+6 a_3-6 a_4+2 k_5+2 k_5^2-2 k_6-2 a_4 k_6-4 k_6^2+2 a_3 k_7+2 k_7^2-2 a_2 k_8-4 a_4 k_8-4 k_6 k_8$ is divisible by $8$. So $a_1-a_2$ is odd.
\item $24-2 a_0+10 a_1-10 a_2+30 a_3-30 a_4-a_4^2+12 k_5+4 a_3 k_5+20 k_5^2-12 k_6-4 a_2 k_6-20 a_4 k_6-28 k_6^2-8 k_7+4 a_1 k_7+20 a_3 k_7+24 k_5 k_7+12 k_7^2-24 k_8-4 a_0 k_8-20 a_2 k_8-28 a_4 k_8-24k_6 k_8-24 k_8^2$ is divisible by $32$, so \underline{$a_4$ is even}, and $a_0+a_3$ is odd.
\item $84+a_3^2+140 k_5+40 k_5^2+224 k_7+160 k_5 k_7+140 k_7^2+4 a_1 \left(5+k_5+4 k_7\right)+8 a_3 \left(5+2 k_5+5 k_7\right)-4 \left(1+a_0+a_2+2 k_4+2 k_6+2 k_8\right) \left(2+k_6+4 k_8\right)-4 \left(3+k_4+3 k_6+6 k_8\right) \left(4+a_2+4 k_4+6 k_6+8 k_8\right)$ is divisible by $64$ and so \underline{$a_3$ is even, and $a_0$ is odd.}
\item $-4+20 a_1-20 a_2-a_2^2-24 k_0+60 k_3+4 a_1 k_3+12 k_3^2-64 k_4-12 a_2 k_4-8 k_0 k_4-24 k_4^2+84 k_5+12 a_1 k_5+48 k_3 k_5+40 k_5^2-96 k_6-24 a_2 k_6-24 k_0 k_6-80 k_4 k_6-60 k_6^2+112 k_7+24 a_1 k_7+80 k_3 k_7+120 k_5 k_7+84 k_7^2-136 k_8-40 a_2 k_8-48 k_0 k_8-120 k_4 k_8-168 k_6 k_8-112 k_8^2$ is divisible by $128$,  so \underline{$a_2$ is even and $a_1$ is odd}.
\end{itemize}
Subsequently we reduce the coefficients modulo $2$ to see that
\begin{itemize}
\item $k_1+k_2$ is even 
\item $k_1+k_3$ is even
\item $k_1+k_4$ is odd
\item $k_1+k_5$ is odd
\item $k_1+k_6$ is even
\item $k_1+ k_7$ is even
\item $k_0+k_1+ k_8$ is even
\end{itemize}
 Thus we substitute 
 \[ k_2=2m_2+k_1;k_3=2m_3+k_1;k_4=2m_4+1+k_1;k_5=2m_5+1+k_1; \]\[k_6=2m_6+k_1;k_7=2m_7+k_1;k_8=2m_8+1+k_1;k_0=2m_0+1\]
 We expand out the coefficients of the polynomials once again to see
 \begin{itemize}
\item $1+m_2+m_5+ m_6$ is even 
\item $1+ m_3+ m_5+ m_7$ is even
\item $k_1+m_0+m_4+ m_8$ is even
\item $k_1+m_3+m_7$ is even
\item $k_1+m_2+m_6$ is even
\end{itemize}
We then make substitutions once again of the sort $m_8=2n_8-k_1-m_0-m_4$, $m_7=2n_7-m_3-k_1$, and $m_6=2n_6-m_2-k_1$ and simplify the coefficients once again checking parity.  We find
 \begin{itemize}
\item $k_1+n_5+n_6+n_7+n_8$ is even 
\item $1+n_6+n_8$ is even
\item $1+n_6+n_8+m_2+ m_3+ m_4$ is even
\end{itemize}
make the substitutions
\[ n_7= 2r_7-\left(1+k_1+ n_5\right); n_8= 2r_8-\left(1+ n_6\right); \] \[ m_4=2 n_4-\left(1+n_6+n_8+m_2+ m_3\right); n_6=\left(2r_6-\left(1+ 2r_8-\left(k_1+n_5\right)\right)\right) \]

At this point every coefficient is an integral multiple of $512$ except the constant term coefficient which is $256 \left(r_7^2-r_8^2\right)$.  However $2$ is not a difference of squares so the constant term is at least $1024$ and $\alpha$ is not a unit.  

\medskip

We remark that a similar situation seems to arise in all odd degrees,  we also checked it for $3,5$ and $7$.

\section{The case that $p=2$ and $p=4$.}

This case is quite different to our previous searches in that much of the work has been done in early work of Flammang and Rhin \cite{FR}.  Motivated by a question of ours in addressing the identification of all arithmetic Kleinian groups generated by an elliptic of order $2$ and an elliptic of order $3$ they proved the following.

\begin{theorem}  There are $15909$ algebraic integers $\alpha$ whose minimal polynomial has exactly one complex conjugate pair of roots in the ellipse
\[ {\cal E} = \{z\in \IC: |z+1|+|z-2|\geq 5 \} \]
and such that all the real roots of this polynomial lie in the interval $[-1,2]$.  The degree of such an $\alpha$ is less than $10$,  and 
\begin{itemize}
\item $22$ polynomials in degree $2$,
\item $206$ polynomials in degree $3$,
\item $918$ polynomials in degree $4$,
\item $2524$  polynomials in degree $5$,
\item $4401$  polynomials in degree $6$,
\item $4260$  polynomials in degree $7$,
\item $2792$  polynomials in degree $8$,
\item $600$ polynomials in degree $9$, and
\item $186$  polynomials in degree $10$.
\end{itemize}
\end{theorem}

In order to use this list there are a few things we must establish.  First that the exterior of the moduli space of groups freely generated be elements of order $2$ and $4$, $\IC\setminus {\cal M}_{2,4}$ lies within this ellipse.  In fact it does not,  however this space has a $\IZ_2$ symmetry in the line $\{\Re e(z)=-1\}$ which respects discreteness (however the symmetric pair of discrete groups are not necessarily isomorphic --- only having a common index two subgroup.)  Thus we need only show $\{z:\Re(z)\geq -1\} \setminus {\cal M}_{2,4}$ lies in the ellipise.  Again it does not,  but an integer translate by $1$ of it does!  This observation was made by Zhang in her PhD thesis \cite{Zhang}.  At that time we did not have the technology to decide the question of whether a group had finite co-volume or not though Coopers PhD thesis \cite{Cooper} implemented a modified version of Week's Snappea code to decide for many of them. 

Actually we additionally searched these spaces,  using Flammang and Rhin's degree bounds,  just to check that we had found all the possibilities.

\medskip

A rough description of the space $\IC\setminus {\cal M}_{2,4}$ can be found as all the roots of all Farey polynomial lie in it.  Such a description is illustrated below.

 \scalebox{0.75}{\includegraphics[viewport=50 460 470 780]{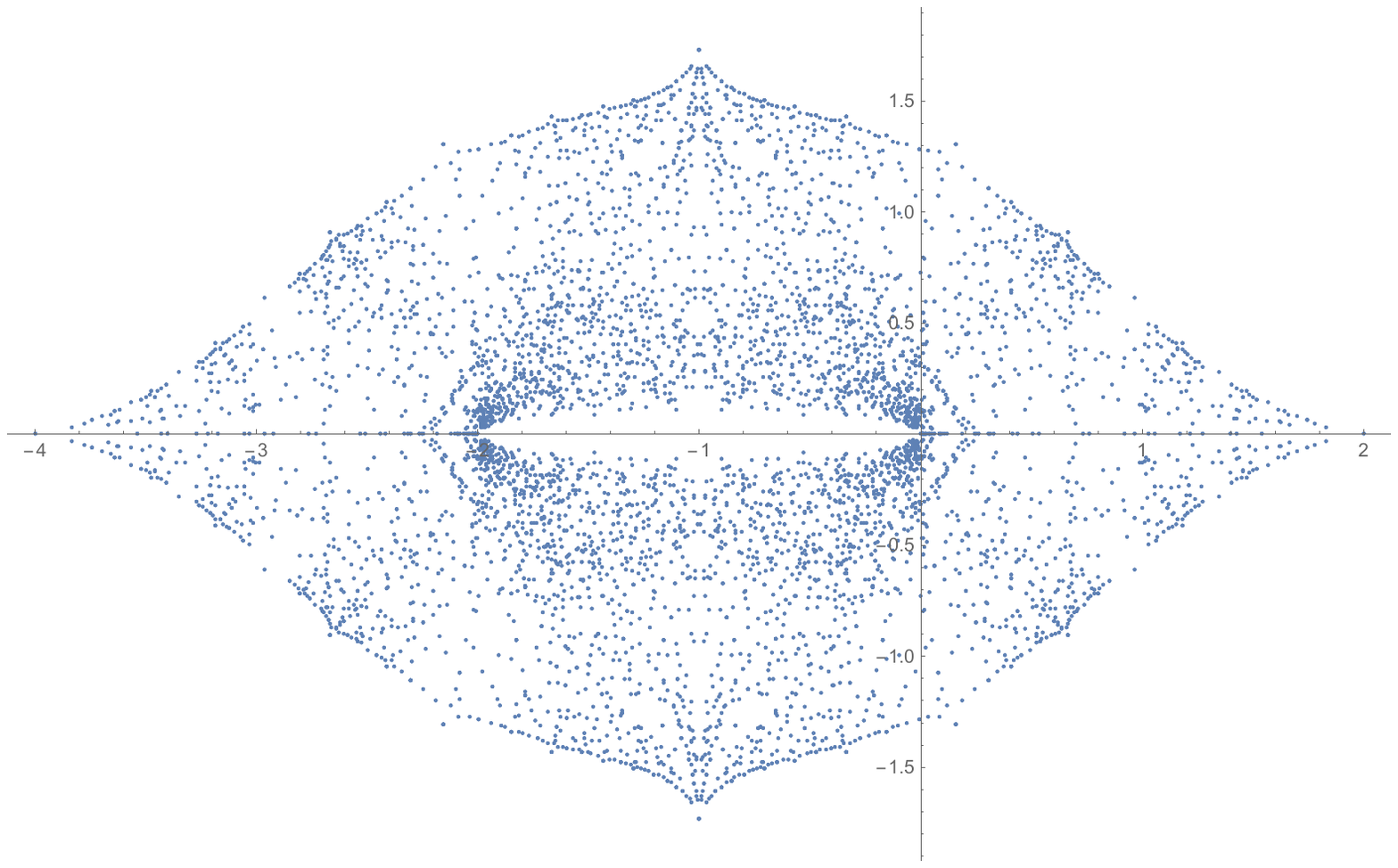}}\\
{\bf Figure 12.} {\em  The space $\IC\setminus {\cal M}_{2,4}$ found by roots of Farey polynomials.  Notice the symmetry in $\Re e(z)=-1$ }

\medskip

We will now give a provable description of this space in order to get reasonable degree bounds for $\alpha=\gamma+1$ to limit the actual number of polynomials we have to consider.  As previously we enumerate the ``first'' 129 Farey words.  We compute the pleating rays corresponding to the following $14$ slopes.
\[ F_{rat}=\{\frac{1}{2},\frac{2}{3},\frac{3}{4},\frac{3}{5},\frac{4}{5},\frac{5}{6},\frac{4}{7},\frac{5}{7},\frac{6}{7},\frac{5}{8},\frac{7}{8},\frac{5}{9},\frac{7}{9},\frac{8}{9}\}.\]
The associated Farey polynomials $p_{r/s}(z)$ are monic,  with integer coefficients and have degree no more than $9$.  The pleating ray is a particular branch of $p_{r/s}^{-1}((-\infty,-2])$, actually the branch making angle $r\pi/s$ at $\infty$ with the positive real axis.  We then compute $D_{r/s}=p_{r/s}(\{\Re e(z)\leq -2)$,  the $r/s$ pleating ray neighbourhood,  and the results of \cite{EMS} show that for all $r/s$,  $D_{r/s}\subset \overline{{\cal M}_{2,4}}$ touching $\partial {\cal M}_{2,4}$ at a single point -- the $r/s$-cusp group. This is illustrated below.

 \scalebox{0.5}{\includegraphics[angle=-90,viewport=-40 50 600 620]{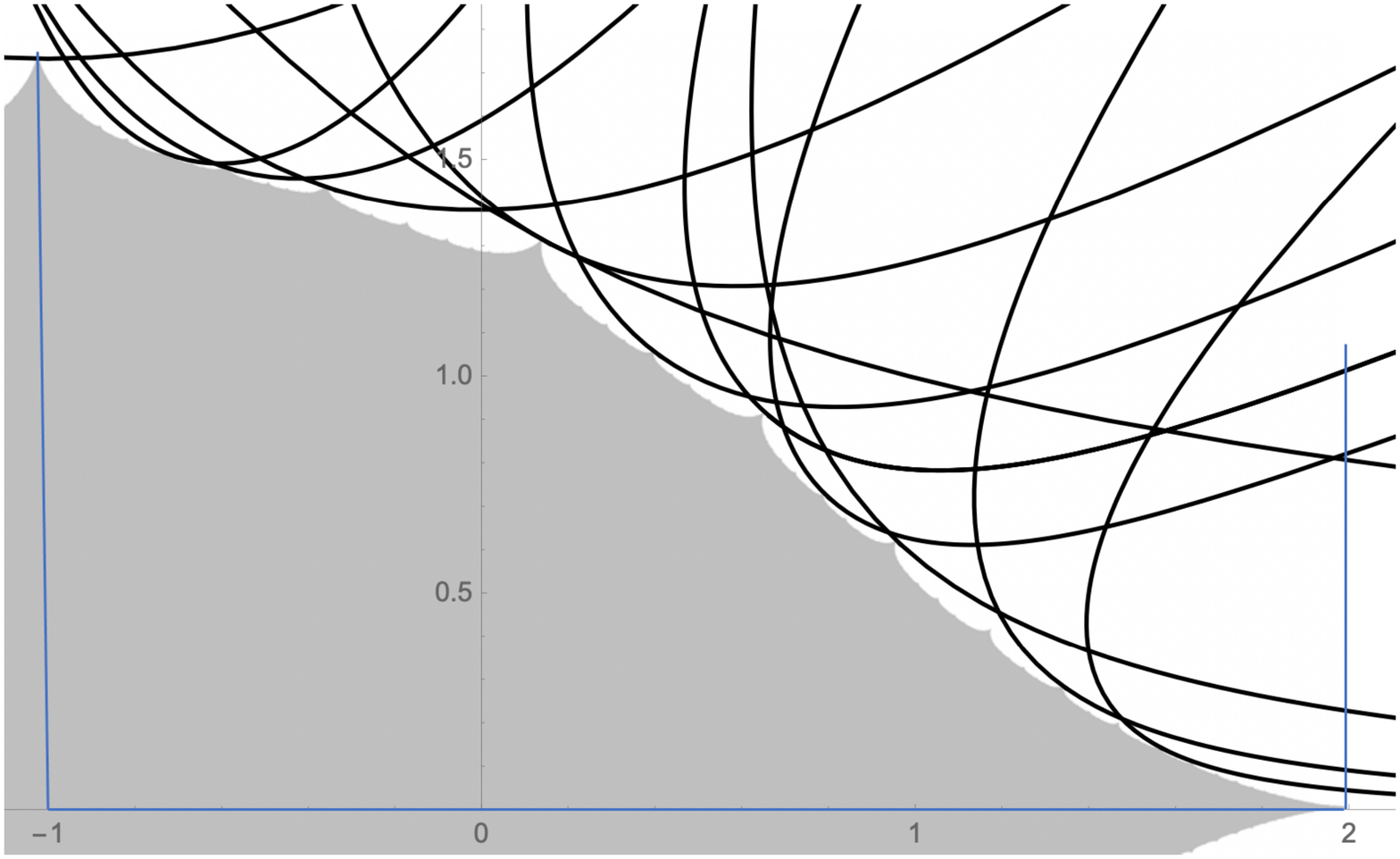}}\\
\noindent {\bf Figure 14.} {\em  The space $\IC\setminus {\cal M}_{2,4}$ (grey) and pleating ray neighbourhoods $D_{r/s}$.  The grey region is found by a conjectural description of the space of groups free on two generators given by Bowditch, \cite{Bowditch}.  The bounded parabaloid shaped regions are proved to lie completely in ${\cal M}_{2,4}$.}

\medskip

This gives us an approximation to the set ${\cal M}_{2,4}$. Should we find a value $\gamma\in \{z:\Re e(z)>-1,\Im m(z)>0\}$ which does not lie in the bounded region

\[ \IC \setminus \Big\{ \{-1 < \Re e(z)<2 \}\cup \bigcup_{\frac{r}{s}\in F_{rat}} D_{r/s} \Big\},    \]
then we know that the group generated by elements of order two and four with this particular commutator value is free on the two generators.

\medskip We remark that it is quite straightforward to test if a point is in some $D_{r/s}$.  We simply evaluate the Farey polynomial $p_{r/s}$ on it to see the image has real part less than $-2$ (we can guess the particular $r/s$ by inspection).  Then we check that this point lies in the same branch of $p_{r/s}^{-1}(\{\Re e(z)>-2\})$ as that determined by the pleating ray from an elementary path lifting which we can also implement computationally.

\bigskip

\subsection{The symmetry of the space ${\cal M}_{2,4}$}

First we describe the symmetry of the space ${\cal M}_{2,4}$.

\begin{theorem} Let $\langle f,g\rangle$ be a Kleinian group generated by elements $g$ of order $2$ and $f$ of order $4$ and set $\gamma=\gamma(f,g)$.  Then there is $\langle f,h\rangle$,  a Kleinian group generated by elements of order $2$ and $4$ and such that $\gamma(f,h)=-2-\gamma$. The groups  $\langle f,g\rangle$ and $\langle f,h\rangle$ have a common index two subgroup $\langle f,gfg^{-1}\rangle$ and so both are simultaneously lattices (or not).
\end{theorem}
\noindent{\bf Proof.}  The axis of $g$ bisects (and is perpendicular to) the common perpendicular $\ell$ between the axes of $f$ and $gfg^{-1}$.  Let $h$ be the elliptic of order two whose axis is perpendicular to both that of $g$ and $\ell$.  Then evidently $hfh^{-1}$ is an element of order the same as $f$ and it shares the same axis of $f$ and has the same trace.  Thus $hfh^{-1}=f^{\pm1}$.  We calculate that $\gamma(f,hfh^{-1})=\gamma(f,h)(\gamma(f,h)+2) = \gamma(f,g)(\gamma(f,g)+2)$ and a moments analysis shows $\gamma(f,g)\neq \gamma(f,h)$ so $\gamma(f,g)=-2+\gamma(f,h)$.  The statement regarding index two follows immediately by looking at the coset decomposition as $g$ and $h$ both have order two. \hfill $\Box$

 The figure below clearly indicates,  and it is not too difficult to prove just using the combinatorics of the isometric circles to construct a fundamental domain,  that 
\begin{equation} 
\IC \setminus {\cal M}_{2,4} \subset  \{-2+{\cal E}\} \cup \{-1+{\cal E}\}
\end{equation}
  \scalebox{0.73}{\includegraphics[viewport=40 420 770 750]{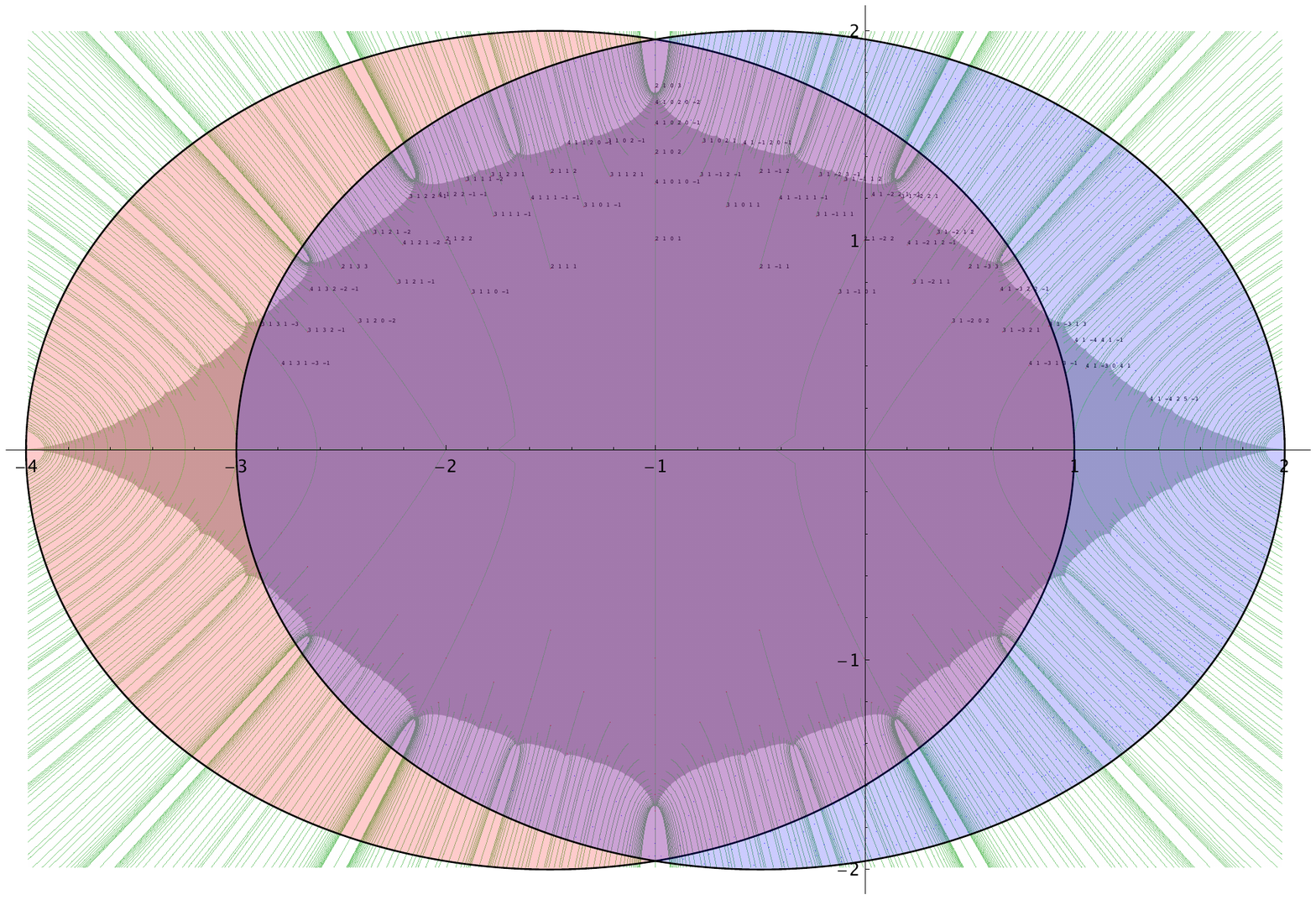}}\\
    
Notice that the interval $[-2,0]$ lies in both $\{-2+{\cal E}\}$ and $\{-1+{\cal E}\}$.  Thus we find that the values of $\gamma$ that we seek are those points in $-1+{\cal E}$ whose minimal polynomial has real roots which actually lie inside the interval $[-2,0]$ (recall they do lie in $[-2,1]$).  Once we have this list of points we go through the process of deciding if they are in or out of $\IC \setminus {\cal M}_{2,4}$.  

\medskip

Now using this symmetry we may assume $0 \leq \Re e(\alpha) \leq 3$.  We also have the elementary bounds $|\alpha|\leq 3$,  $|\alpha+1|\leq 4$ and $|\alpha-1|\leq 2$.  Now we obtain degree bounds as earlier following (\ref{alphabounds}).

Let 
\[ q_\alpha(z)=z^n+a_{n-1}z^{n-1}+\cdots+a_1z+a_0 = (z-\alpha)(z-\bar\alpha)(z-r_1) \cdots(z-r_{n-2}) \]
be the minimal polynomial for $\alpha$.  Then,  as $r_i\in [-1,1]$,
\begin{eqnarray} |q_\alpha(-1)||q_\alpha(0)||q_\alpha(1)| &=& |\alpha|^2|\alpha^2-1|^2\prod_{i=1}^{n-2}  |r_i| |1-r_{i}|\nonumber  \\
& \leq & 576 (0.3849)^{n-2}.  
\end{eqnarray}
as $\max \{x(x^2-1):x\in [-1,1]\} = 0.3849...$.
As the left hand side here is greater than one,  we deduce that $n\leq 8$.  Indeed if $p(0)\neq \pm 1$,  then $n\leq 7$.  

We therefore run through the $12331$ polynomials in the list above, along with the degree 8 polynomials with last coefficient 1.  We are only interested in those whose real roots lie in $[-1,1]$ and complex root has real part of absolute value no more than $\sqrt{3}$.

\subsection{Degree $2$} There are $12$ possible values,  of which $6$ lie in the region $\IC \setminus \Big\{ \{-1 < \Re e(z)<2 \}\cup \bigcup_{\frac{r}{s}\in F_{rat}} D_{r/s} \Big\}$.\\
 
 \scalebox{0.5}{\includegraphics[angle=-90,viewport=40 50 600 660]{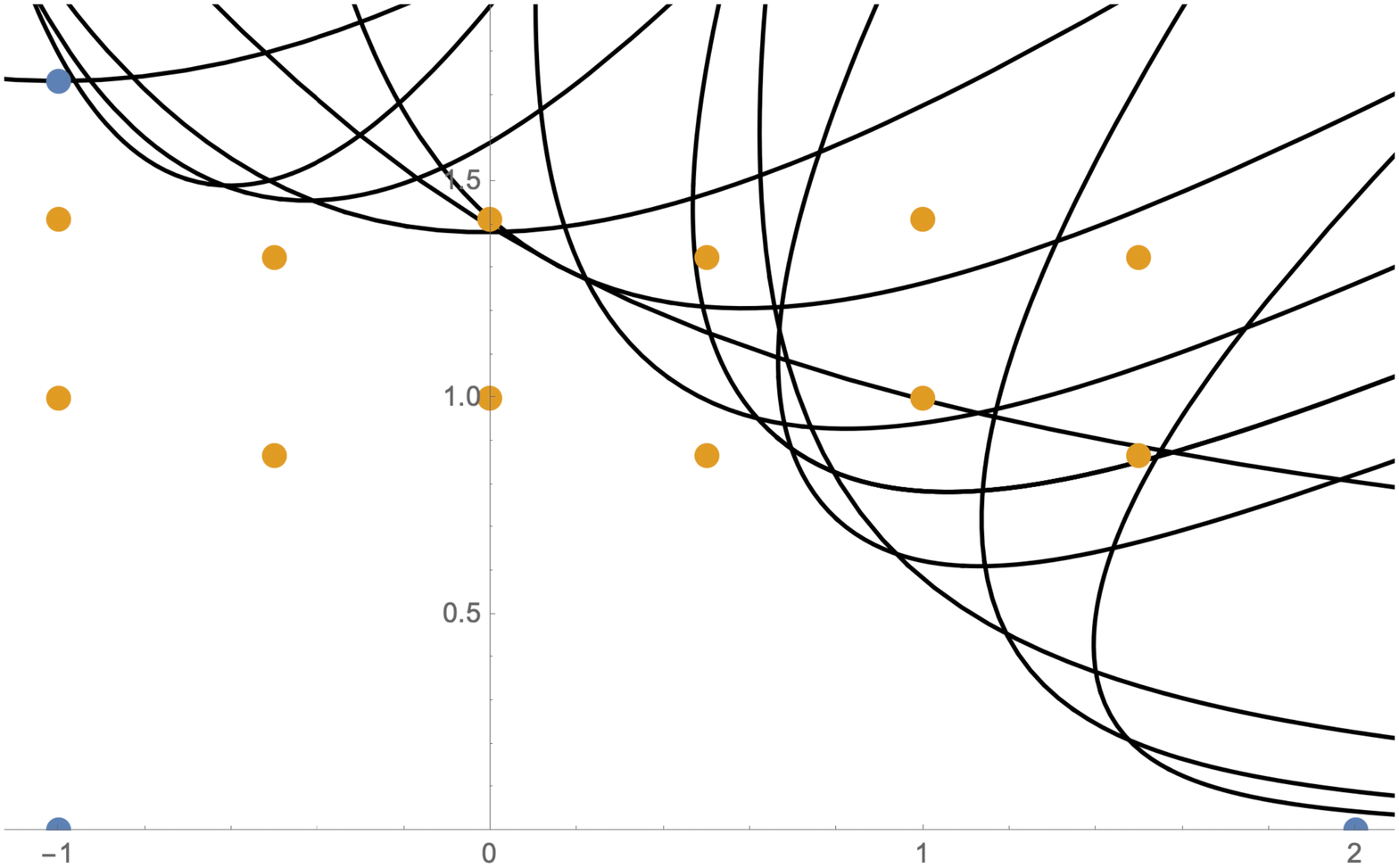}}\\

\subsection{Degree $3$} There are $175$ possible values,  of which $13$ lie in the region $\IC \setminus \Big\{ \{-1 < \Re e(z)<2 \}\cup \bigcup_{\frac{r}{s}\in F_{rat}} D_{r/s} \Big\}$.\\
 
 \scalebox{0.5}{\includegraphics[angle=-90,viewport=40 50 600 660]{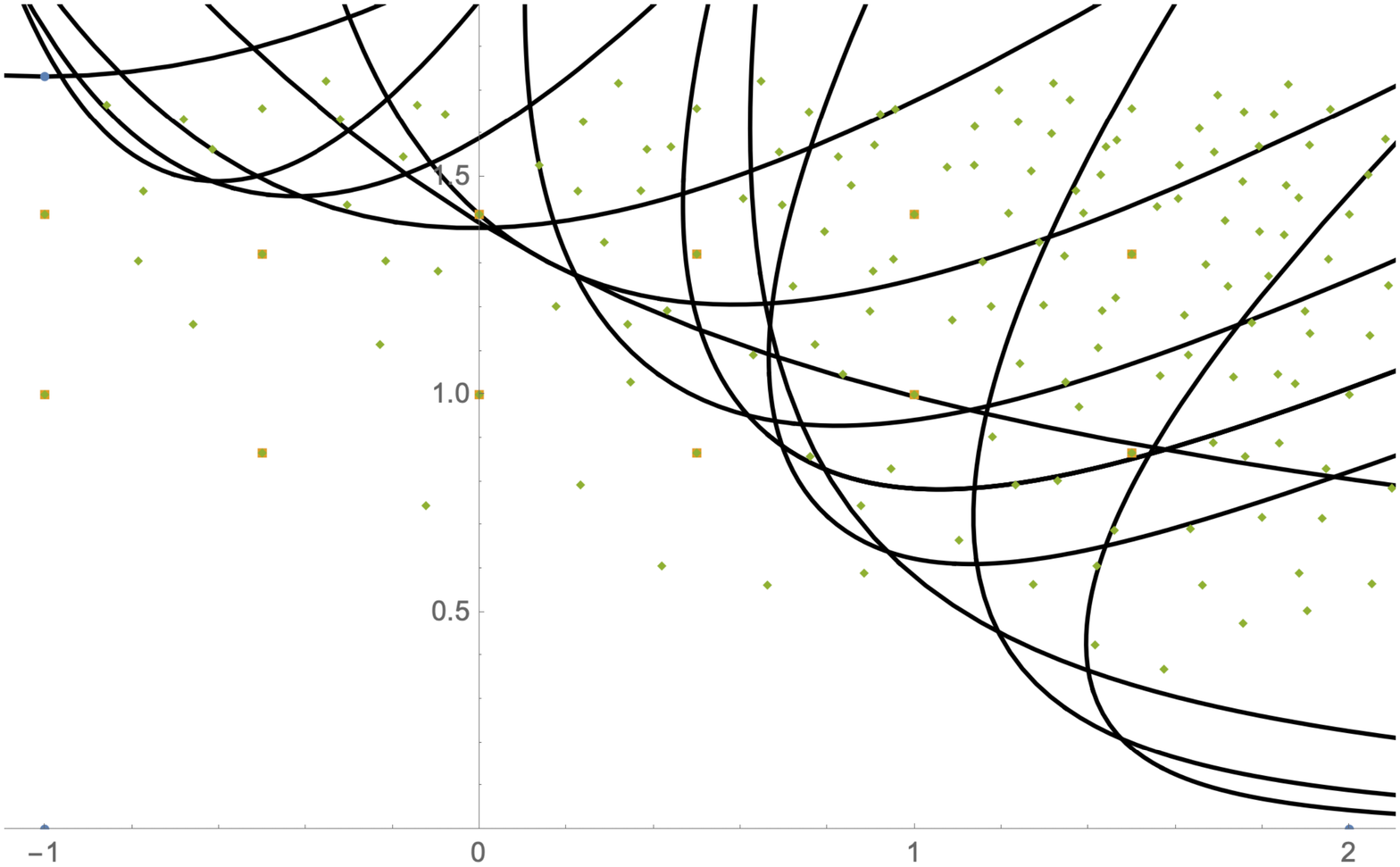}}\\

\subsection{Degree $4$} With rather coarse search bounds 
\[  -8\leq a \leq 1, -12\leq  b \leq 22,  -23\leq c \leq 23, -8\leq d \leq 8 \]
  
  we found $572$ possible values,  of which $23$ lie in the region $\IC \setminus \Big\{ \{-1 < \Re e(z)<2 \}\cup \bigcup_{\frac{r}{s}\in F_{rat}} D_{r/s} \Big\}$.  \\
  
 \scalebox{0.75}{\includegraphics[viewport=50 500 750 800]{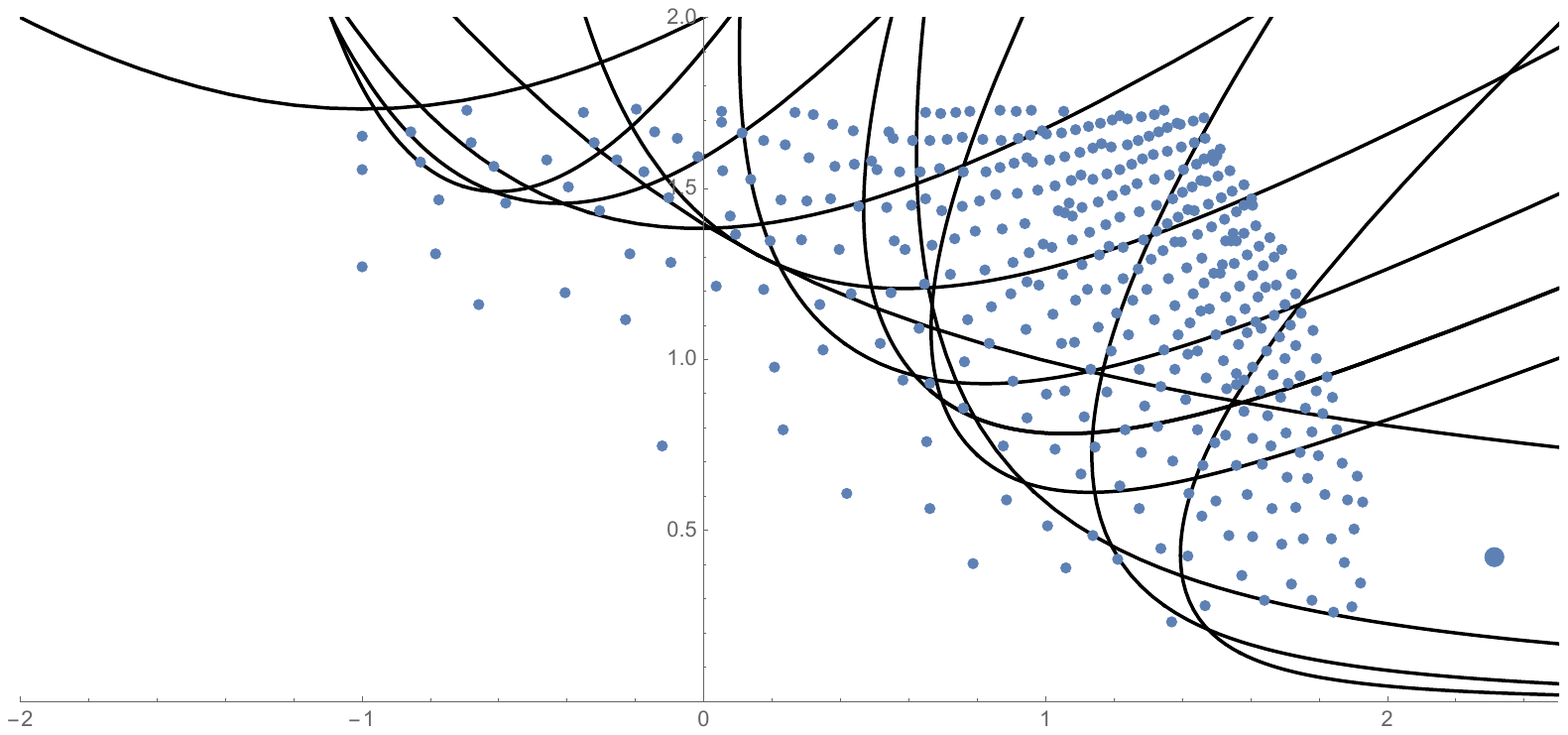}}\\
 
  We then checked irreducibility of the polynomial and this left $13$ possibilities in the region $\{\Re e(z)\geq -1\}$ and one further point was eliminated by adding another pleating ray neighbourhood.  This leaves the list we have presented.

\subsection{The case $p=4$.}  In this case we have a group generated by two elements of order $4$.  The following theorem shows that if $\langle f,g\rangle$ is such a group,  then there is a $\IZ_2$ extension to a group $\langle f,\Phi \rangle$ generated by elliptics of orders $2$ and $4$.  Of course if $\langle f,g\rangle$ is a discrete subgroup of an arithmetic Kleinian group,  then so is $\langle f,\Phi \rangle$.  The only issue is the calculation of the relevant commutators,  the arithmetic data, and identifying the slope.

\begin{theorem}  Let $\langle f,g\rangle$, generated by two elements of order $4$,  be a discrete subgroup of an arithmetic Kleinian group and $\gamma=\gamma(f,g)$ with $\gamma$ complex.  Suppose $\tilde{\gamma}(\tilde{\gamma}+2)=\gamma$.  Then there is $\Phi$ of order two and such that $\langle f,\Phi\rangle$ is a discrete subgroup of an arithmetic Kleinian group and that $\tilde{\gamma}=\gamma(f,\Phi)$.
\end{theorem} 
\noindent{\bf Proof.} One of the two elliptics of order two whose axis bisects the common perpendicular of the axes of $f$ and of $g$ and interchanges these axes is the element $
\Phi$ we seek.  Then $\Phi f\Phi^{-1}=g^{\pm 1}$ and $\langle f,\Phi\rangle$ is discrete as it contains a discrete group of index two. Similarly these groups are commensurable, and have the same invariant trace field. Arithmeticity is preserved by finite extensions.  To see this more concretely,  the factorisation condition must hold for $\gamma$,  that condition is precisely that $\IQ(\tilde{\gamma})=\IQ(\gamma)$. We have already seen the trace identity
\[ \gamma=\gamma(f,g)=\gamma(f,\Phi f\Phi)=\gamma(f,\Phi)(\gamma(f,\Phi)+2) =\tilde{\gamma}(\tilde{\gamma}+2) \]
 This completes the proof. \hfill $\Box$
 
 \medskip
 
 This result now shows us the can find all the groups in the case $p=4$ from our analysis of the previous case $p=2$.
 
 \section{Finding the groups.}
 
 The recent solution of Agol's conjecture \cite{ALSS,AOPSY} identifying all the Kleinian groups generated by two parabolic elements as two bridge knot and link groups and associated Heckoid groups suggests a method for identifying our groups.
 
 We  conjecture that the same is true for groups generated by two elements of finite order.  This suggests the following approach to identifying these groups.  Let $\Gamma=\langle f,g\rangle$ with $o(f)=4$ and $o(g)=p$.   We only sketch the ideas here.
 
 Let $m/n$ and $r/s$ be two rational numbers such that $m$, $n$ (resp. $r$, $s$) are coprime. Define an operation $\oplus$ such that $m/n\oplus r/s=(m + r)/(n + s)$; we call $\oplus$ Farey addition. Suppose $m/n < r/s$. Then $m/n < (m + r)/(n + s) < r/s$. These Farey fractions enumerate the simple closed curves on the four times punctured sphere,  and hence words in the free group of rank two.  These words correspond to pleating rays.  As examples: We start with a pair of rationals and associated rational words (in generators $X$ and $Y$,  with $x=X^{-1}$ and $y=Y^{-1}$): 
 \[ \{1/2, 1/1\} \mapsto \{\{X, Y, x, y\}, \{X, y\}\}\]
 and inductively create Farey words by addition;
\begin{eqnarray*} \left\{\frac{1}{2},\frac{2}{3},1\right\} &\mapsto& \{\{X,Y,x,y\},\{X,Y,x,Y,X,y\},\{X,y\}\} \\
 \left\{\frac{1}{2},\frac{3}{5},\frac{2}{3},\frac{3}{4},1\right\}&\mapsto&\{\{X,Y,x,y\},\{X,Y,x,y,X,y,x,Y,X,y\},\\  && \{X,Y,x,Y,X,y\},\{X,Y,x,Y,x,y,X,y\},\{X,y\}\} \end{eqnarray*}

Here is some Mathematica code which does this (from Zhang \cite{Z}).

\medskip
{\tiny \noindent L1 = {1/2, 1/1};\\
L2 = {{X, Y, x, y}, {X, y}} (*list of rational words*)\\
NF = 3 (*determines how many polynomials to make *)\\
For[j = 1, j$\leq$ NF,  j++,(* determines the number of rational words made*)] \\
d = Length[L1];\\
For[i=1, i$\leq$(2 d - 1), i=i + 2,\\
 N1 = (Numerator[L1[[i]]] + Numerator[L1[[i + 1]]] )/(Denominator[L1[[i]]] + Denominator[L1[[i + 1]]]);\\
 (* Farey addition*)\\
 L1 = Insert[L1, N1, i + 1]; (* orders the fractions *)\\
 T2 = Join[L2[[i]], L2[[i + 1]]]; d1 = Denominator[N1] + 1; \\
 R = Switch[T2[[d1]], x, X, X, x, y, Y, Y, y];\\
 T2 = ReplacePart[T2, R, d1]; L2 = Insert[L2, T2, i + 1]]]}

\medskip
Each entry in the second list corresponds to a word $W_{r/s}$.  We take a representation of a group generated by elliptics of order $4$ and $p$ by assigning
\[ X=\left(\begin{array}{cc}\frac{1+i}{\sqrt{2}}& 1 \\ 0 &\frac{1+i}{\sqrt{2}} \end{array} \right), \;\;\;\;\; 
Y=\left(\begin{array}{cc}\cos(\frac{\pi}{p})+i\sin(\frac{\pi}{p}) & 0 \\ -\mu  &\cos(\frac{\pi}{p})-i\sin(\frac{\pi}{p}) \end{array} \right),\]
Then 
\begin{equation} p_{r/s}(\mu)={\rm Trace}(W_{r/s})\end{equation}
Some examples (with $p=2$ and the first $9$ fractions)
\begin{center}
\begin{tabular}{|c|c|} \hline
Farey &  polynomial \\ \hline
$\frac{1}{2}$ &  $2-2 \sqrt{2} \mu +\mu ^2$ \\ \hline 
$\frac{4}{7}$ & $ -\sqrt{2}+25 \mu -58 \sqrt{2} \mu ^2+118 \mu ^3-65 \sqrt{2} \mu ^4+41 \mu ^5-7 \sqrt{2} \mu ^6+\mu ^7 $\\ \hline 
$\frac{3}{5}$ &  $\sqrt{2}-13 \mu +17 \sqrt{2} \mu ^2-19 \mu ^3+5 \sqrt{2} \mu ^4-\mu ^5$ \\ \hline 
$\frac{5}{8}$ &  $2-16 \sqrt{2} \mu +104 \mu ^2-144 \sqrt{2} \mu ^3+220 \mu ^4-100 \sqrt{2} \mu ^5+54 \mu ^6-8 \sqrt{2} \mu ^7+\mu ^8$ \\ \hline 
$\frac{2}{3}$ &  $-\sqrt{2}+5 \mu -3 \sqrt{2} \mu ^2+\mu ^3$ \\ \hline 
$\frac{5}{7}$ &  $\sqrt{2}-17 \mu +36 \sqrt{2} \mu ^2-84 \mu ^3+55 \sqrt{2} \mu ^4-39 \mu ^5+7 \sqrt{2} \mu ^6-\mu ^7$ \\ \hline 
$\frac{3}{4}$ &  $2-4 \sqrt{2} \mu +10 \mu ^2-4 \sqrt{2} \mu ^3+\mu ^4$ \\ \hline 
$\frac{4}{5}$ &  $-\sqrt{2}+5 \mu -11 \sqrt{2} \mu ^2+17 \mu ^3-5 \sqrt{2} \mu ^4+\mu ^5$ \\ \hline 
$\frac{1}{1}$ &  $\sqrt{2}-\mu$ \\ \hline 
\end{tabular} \\
\end{center}
Notice that the degree of the polynomial is the denominator or the fraction and that $W_{1/2}=XYxy=[X,y]$ and so $p_{1/2}(\mu)-2=\gamma(X,Y)$.

To identify our groups we used the following procedure. From $p$ we construct a list of about $100$ fractions and polynomials $p_{r/s}(\mu)$.  From an monic polynomial $z^3 +5z^2 +7z+1$ candidate we identify $\gamma=-2.41964+0.60629i$ and $\mu=\sqrt{2}+\sqrt{2+\gamma }=1.60761 + 1.5675 i$, from $-2 \sqrt{2} \mu +\mu ^2=\gamma$.  We then evaluate all these polynomials on $\mu$.

\begin{center}
\begin{tabular}{|c|c|} \hline
Farey $r/s$&  polynomial $P_{r/s}$\\ \hline
$\frac{1}{2}$ &  $-0.419643 - 0.606291 i$ \\ \hline 
$\frac{4}{7}$ & $-1.42553 + 1.59871 i$\\ \hline 
$\frac{3}{5}$ &  $0.540536 + 1.03152 i$ \\ \hline 
$\frac{5}{8}$ &  $1. - 2.84217 \times 10^{-14} i $ \\ \hline 
$\frac{2}{3}$ &  $1.02696 - 0.838121 i$ \\ \hline 
$\frac{5}{7}$ &  $-4.56052 + 1.21192 i$ \\ \hline 
$\frac{3}{4}$ &  $2.6478 + 1.72143 i$ \\ \hline 
$\frac{4}{5}$ &  $-3.39159 + 2.16591 i$ \\ \hline 
$\frac{1}{1}$ &  $0.398566 - 0.760591 i$ \\ \hline 
\end{tabular} \\
\end{center}
This list suggests that ${\rm Trace}(W_{5/8})=1$ and that therefore $W_{5/8}^3=1$.  We could deduce that this point $\gamma$ gives us the Heckoid group $(5/8;3)$. In this case the two bridge link with slope $5/8$ one component surgered by $(4,0)$ Dehn surgery and the other by $(2,0)$ Dehn surgery with the tunnelling word $W_{5/8}$ elliptic of order three.

We presented this data with round off error,   and we must decide if it is real or not.  There are two fortunate things here.  First we know apriori from the Identification Theorem that this group generated by elements of order $2$ and $4$ with $\gamma$ as the commutator {\em is} discrete.  We are just trying to find out what it is.  Secondly $\gamma$ is presented as an algebraic integer so we can use integer arithmetic,  identify the minimal polynomial of $P_{r/s}(\gamma)$ and its roots.  We can then assure ourselves that $1$ is indeed a root and that no others are close.  Or we could go back and enter $\gamma$ and compute $W_{r/s}$ symbolically.  These amount to the same thing of course.

In the case at hand $\mu = \sqrt{2}+ \sqrt{2 + \gamma}$ and $\mu$ is a root of
\[ 1 - 298 z^2 + 251 z^4 - 200 z^6 + 71 z^8 - 14 z^{10} + z^{12} \]
The minimal polynomial of $P_{5/8}(\mu)$ is then found to be $z-1$.

 This procedure works well in almost all cases.  However some of the cases required us searching through a few thousand slopes.  Doing this using integer arithmetic takes a very long time,  while using numerical approximation introduces round-off errors as there are approximately four times as many matrices to multiply together as the denominator of the slope (so sometimes around $500$ matrices).  We found that working with 40 digit precision (and hoping for a result accurate to a couple of decimal places) gave us an answer in a reasonable time - a couple of hours.  Of course once you have a possible answer  it is much easier and quicker to verify it is in fact correct.

\medskip
\noindent G.J. Martin, Institute for Advanced Study,  Massey University, New Zealand.\\ (g.j.martin@massey.ac.nz)\\
K. Salehi, Institute for Advanced Study,  Massey University, New Zealand.\\
Y. Yamishita, Department of Mathematics, Nara Women's University, Japan.

\end{document}